\documentclass[10pt, draft]{article}
\usepackage{amsmath}
\usepackage{amssymb}
\usepackage{latexsym}
\allowdisplaybreaks[4]
\date{}
\renewcommand{\theequation}{\arabic{section}.\arabic{equation}}
\newtheorem{Thm}{Theorem}[section]
\newtheorem{Prop}[Thm]{Proposition}
\newtheorem{Lemma}[Thm]{Lemma}

\newtheorem{Remark}[Thm]{Remark}
\renewcommand{\v}[1]{\vert #1 \vert}
\renewcommand{\Bbb}{\mathbb}
\newcommand{\V}[1]{\Vert #1 \Vert}

\newcommand{\R}{{\mathbb R}}
\newcommand{\C}{{\mathbb C}}
\title{An inverse problem for a three-dimensional heat equation 
in bounded regions with several convex cavities}
\author{Mishio Kawashita\thanks{Partly supported by JSPS KAKENHI Grant Number JP16K05232.}
\thanks{kawasita@hiroshima-u.ac.jp}
}
\date{}
\begin{document}
\maketitle
\begin{abstract}
\par

\vskip1pc
\end{abstract}
\par
In this paper, an inverse initial-boundary value problem for 
the heat equation in three dimensions is studied. Assume that
a three-dimensional heat conductive body contains several cavities of strictly convex. 
In the outside boundary of this body, a single pair of the temperature 
and heat flux is given as an observation datum for the inverse problem.  
It is found the minimum length 
of broken paths connecting arbitrary fixed point in the outside, 
a point on the boundary of the cavities and 
a point on the outside boundary in this order, if the minimum path is not line segment. 

%


\vskip 1truepc
\setcounter{section}{0}
\section{Introduction}\label{Introduction}
\vskip-6pt\noindent

Let $\Omega$ be a bounded domain of $\R^3$ with $C^{2}$ boundary. 
Let $D$ be an open subset of $\Omega$ with 
$C^{2}$ boundary and satisfy that
$\overline D\subset\Omega$ and $\Omega\setminus\overline D$ is connected.
We denote by $\nu_\xi$ and $\nu_y$ the unit outward normal vectors at
$\xi \in\partial D$ and
$y \in \partial\Omega$ on $\partial D$ and $\partial\Omega$ respectively.
\par
Let $T > 0$ be a fixed constant and $\rho$ be a continuous function on
$\partial{D}$. Consider the following initial boundary value problem
of the usual heat equation:
\begin{equation}
\left\{\hskip-12pt
\begin{array}{lll}
&(\partial_t - \triangle)u(t, x) = 0 
&\mbox{in } (0, T)\times\Omega\setminus\overline{D}, \\
&
\partial_{\nu}u(t, x) = f(t, x)
 & \mbox{on } (0, T)\times\partial\Omega,  \\
&(\partial_\nu+\rho(x))u(t, x) = 0
 & \mbox{on } (0, T)\times\partial{D},  \\
&u(0, x) = 0 
&\mbox{on } \Omega\setminus\overline{D}, \\
\end{array}
\right.
\label{Eq in inside for the measurement}
\end{equation} 
where $\partial_\nu = \sum_{j = 1}^3(\nu_x)_j\partial_{x_j}$ for 
$x \in \partial{D}\cup\partial\Omega$. 

\par
Mathematical studies on inverse problems arising thermal imaging are formulated 
as the boundary inverse problems for the usual heat equation.
In this inverse problem, pairs of the measurement data $(u, f)$ on the outside 
boundary, i.e. the temperature $u$ and the heat flux $f$ on 
$(0, T)\times\partial\Omega$, are given as observation data. 
The problem is to understand what information on $\partial{D}$
can be extracted by using these data on the outside boundary.
\par
Elayyan and Isakov \cite{EI} investigated the uniqueness problem corresponding 
to this type of inverse problem, which determine $D$ and $\rho$ uniquely 
by using infinitely many observation data.
In this paper, the case that 
infinitely many observation data are used to obtain inside information
is called by \lq\lq infinite measurement". The completely opposite case to 
infinite measurement is \lq\lq one measurement". This is the case that
only one pair of the observation data $(u, f)$ is allowed to use as the observation
data for inverse problem. In one measurement case, as in 
Bryan and Caudill \cite{BC}, the uniqueness results fail 
if the initial condition does not vanish. Hence, to handle
one measurement case, as in (\ref{Eq in inside for the measurement}), 
we need to assume that $u(0, x) = 0$ on $\Omega$.
\par
Other important problems are for stability and reconstruction. The stability
problem is to show continuous properties between the observation data and the unknown objects
($D$ and $\rho$). For stability problems, see Vessella \cite{Vessella1}, 
and references therein. 

\par
In this paper, the problems concerning reconstruction procedure, that is to find information on 
$D$ or $\rho$ from the observation data, are treated. For this problem, several methods are 
proposed by many authors.
For a one space dimensional case, Daido, Kang, and Nakamura \cite{DKN} gives 
an approach for this type of inverse problem by using
an analogue of the probe method introduced by Ikehata \cite{IP1}. This 
procedure is numerically simulated by Daido, Lei, Liu and Nakamura \cite{DLLN}.
\par
Various approaches for reconstruction procedures are proposed, however, 
there are relationships among them although the formulations are not 
similar to each other. 
These relations are found by Honda, Nakamura, Potthast and Sini \cite{HGPS}.
In the author's best knowledge, 
only the enclosure method is different from many of other approaches.
Thus, it is worth investigating inverse problems by the enclosure method.
\par

The enclosure method is originally developed by Ikehata \cite{I1} and \cite{I2}
for static problems formulated by elliptic boundary value problems. 
About boundary inverse problems for the heat equations, 
infinite measurement cases are treated by \cite{IK0}.
In \cite{IK00}, the case of inclusions (i.e. the case that $D$ is 
filled by other medium) is considered.  
In the case that the inside boundary of the inclusion may depend on time variable 
$t$ and is strictly convex for all $t$, Gaitan, Isozaki, Poisson, Siltanen and 
Tamminen \cite{Gaitan Isozaki et.al. 3} also investigated by using a similar
approach to the enclosure method. 
\par
Usually, to give reconstruction procedure, functions called \lq\lq indicator"
defined by using the observation data are introduced. 
From asymptotic behaviors of indicator functions, 
one can obtain informations for the inside. 
In \cite{IK0} and \cite{IK00},  
$h_D(\omega) = \sup_{x \in \partial{D}}x\cdot\omega$ ($\omega \in S^2$ for the 
three dimensional case), $d_D(p) = \inf_{x \in D}\v{x - p}$ 
($p \in \R^3\setminus\overline{\Omega}$), and
$R_D(q) = \sup_{x \in D}\v{x - q}$ ($q \in \R^3$) are extracted.  
Hence $D$ is enclosed by the sets such 
$\cap_{\omega \in S^2}\{ x \in \R^3 \vert x\cdot\omega < h_D(\omega)\}$
as $\cap_{p \in \R^3\setminus\overline{\Omega}}\{ x \in \R^3 \vert \v{x - p} > d_D(p) \}$
and $\cap_{q \in \R^3}\{ x \in \R^3 \vert \v{x - q} < R_D(q)\}$, which are
the origin of the word \lq\lq the enclosure method" as introduced in \cite{I1} and \cite{I2}.

\par
As stated in \cite{Ikehata and Kawashita 2}, infinite observation cases are
different from a one measurement case. This comes from how to choose
the indicator functions $I_\tau$, which contain a (real or complex) large 
parameter $\tau$ from the observation data $(u, f)$
on $(0, T)\times\partial\Omega$. 
We take functions $v_\tau(t, x)$ with large 
parameter $\tau$ satisfying $(\partial_t+\triangle)v = 0$
in $(0, T)\times\Omega$. From these functions, $I_\tau$ is defined by
\begin{align}
I_\tau = \int_{\partial\Omega}\int_0^T
\big(\partial_{\nu_y}v_{\tau}(t, y)u(t, y) - v_\tau(t, y)f(t, y)\big)dtdS_y.
\label{def of indicator 0}
\end{align}
For infinite measurement cases, 
the boundary data $f$ in (\ref{Eq in inside for the measurement}) can be changed
as it suits for $\partial_{\nu}v_\tau$ on $(0, T)\times\partial\Omega$.
Thus, the observation data $(u, f)$ are designed to obtain information of 
$\partial{D}$. Hence, as above, various amounts related 
to $\partial{D}$ are obtained. Note that most of the works stated above 
are for infinite measurement cases. 

\par
For one measurement case, only one pair of the observation data $(u, f)$ is given.
This means that we can not design the indicator functions like as infinite
measurement cases. Only we can do is to choose $v_\tau(t, x)$ for given $f$. 
One possibility of a choice of $v_\tau(t, x)$ is 
to take $v_\tau(t, x) = e^{-\tau^2t}q(x; \tau)$, where
$q(x, \tau)$ is the solution of 
$$
\left\{\hskip-12pt
\begin{array}{lll}
&(\triangle - \tau^2)q(x; \tau) = 0 
&\mbox{in } \Omega, \\
&
\partial_{\nu_y}q(y; \tau) = \int_{0}^{T}e^{-\tau^2t}f(t, y)dt
 & \mbox{on } \partial\Omega.
\end{array}
\right.
$$
In \cite{IK00}, by using $I_\tau$ with this choice of $v_\tau(t, x)$, 
${\rm dist}(D, \partial\Omega) = \inf\{\v{x - y} \vert x \in D, y \in \partial\Omega\}$
is extracted. 
\par
Another idea choosing $v_\tau(t, x)$ is to put 
$v_\tau(t, x) = e^{-\tau^2t}q(x; \tau)$
with functions $q(x; \tau)$ independent of $f$ and satisfying 
$(\triangle - \tau^2)q(x; \tau) = 0$ in $\Omega$. 
For a fixed $p \in \R^3\setminus\overline{\Omega}$ arbitrary taken, we put
$q(x; \tau) = e^{-\tau\vert x - p \vert}/(2\pi\v{x - p})$. 
This is a good example of $q(x; \tau)$. In \cite{Ikehata and Kawashita 2}, 
it is shown that asymptotic behaviors of the indicator function $I_\tau $ defined in (\ref{def of indicator 0}) 
by using this function $q(x; \tau)$ give 
\begin{equation*}
l(p, D) = \inf_{(\xi, y)\in\partial D\times \partial\Omega}l_p(\xi, y),
\end{equation*}
where
$$
l_p(\xi, y) = \vert p - \xi \vert + \vert \xi - y \vert, \quad
(\xi, y) \in \,\R^3\times\R^3.
$$
In \cite{Ikehata and Kawashita 2}, however, we need to assume strict convexity
of $\partial{D}$. In this paper, we treat the case that $D$ consists of
several strictly convex cavities.
\par
To describe main theorems, we introduce notations. 
Take an arbitrary point $p$ outside $\Omega$, and define $I(\lambda, p)$ by
replacing $v_\tau(t, x)$ in (\ref{def of indicator 0}) for defining $I_\tau$
with $v_\lambda(t, x)=e^{-\lambda^2t}E_{\lambda}(x, p)$, 
where $u(t, x)$ is the solution of (\ref{Eq in inside for the measurement}) 
and $E_{\lambda}(x, y)$ is given by
$$
E_{\lambda}(x, y) = \frac{e^{-\lambda\vert x - y \vert}}
{2\pi\vert x - y \vert},\,x \neq y,\, \vert \text{arg}\,\lambda \vert < \frac{\pi}{4}.
$$
Note that $(\triangle_x-\lambda^2)E_\lambda(x, y) + 2\delta(x - y)=0$ holds
in $\R^3$ in the sense of distribution, and the indicator function $I(\lambda, p)$
under consideration in this paper is of the form
\begin{equation}
I(\lambda, p) = \int_0^T\int_{\partial\Omega}e^{-\lambda^2t}\big(\partial_{\nu_y}E_{\lambda}(y, p)u(t, y)
-E_{\lambda}(y, p)\partial_{\nu_y}u(t, y)\big)dS_ydt.
\label{def of I(p, lambda)}
\end{equation}

\par
For a given $0 < \delta_0 < 1$, we denote by 
${\rm \C}_{\delta_0}$ the set of all complex numbers 
$\lambda$ such that 
$\text{Re}\,\lambda\ge\delta_0\vert\text{Im}\,\lambda\vert$. 
We also define $\Lambda_{\delta_1}$ by 
$$
\Lambda_{\delta_1} = \{ \lambda \in 
\C, 
\vert{\rm Im }\,\lambda\vert 
\leq \delta_1\frac{{\rm Re }\,\lambda}{\log{\rm Re }\,\lambda},
{\rm Re }\,\lambda \geq e\,\}.
$$
For $p \in \R^3\setminus\overline{\Omega}$, 
define
\begin{align*}
{\cal G}(p) &=b\{\xi \in \partial D \,\vert\, \nu_\xi\cdot(p - \xi) = 0 \}, 
\quad
{\cal G}^{\pm}(p) = \{\xi \in \partial D \,\vert\, \pm\nu_\xi\cdot(p - \xi) > 0 \}, 
\\
{\cal M}(p) &= \{(\xi, y) \in \partial D\times\partial\Omega \,\vert
\, l(p, D) = l_p(\xi,y) \}, 
\\
{\cal M}_1(p) &= \{(\xi,y) \in {\cal M}(p) \,\vert\, \xi  \in {\cal G}^+(p), 
\nu_\xi\cdot(y - \xi) > 0 \},
\\
{\cal M}_2^{\pm}(p) &= \{(\xi, y)\in {\cal M}(p) \,\vert\, \xi \in {\cal G}^{\pm}(p),
\pm\nu_\xi\cdot(y - \xi) < 0 \}, \,  
\intertext{and}
{\cal M}_g(p) &= \{(\xi, y) \in {\cal M}(p) \,\vert\, \xi \in {\cal G}(p) \}.
\end{align*}
These notations are the same as in \cite{Ikehata and Kawashita 2}.

\par

Throughout this paper, we put the following assumptions on $\partial{D}$ and $\Omega$:
\begin{enumerate}
\item[(I.1)] $D$ consists of several disjoint convex domains, namely   
$D = \cup_{j = 1}^{N}D_j$, where $D_j$ ($j = 1, 2, \ldots, N$) are 
disjoint domains of strictly convex with the boundaries $\partial{D}_j$ of 
class $C^{2}$.
\item[(I.2)] ${\cal M}_g(p){\cup}{\cal M}_2^{-}(p) = \emptyset$.
\end{enumerate}

Now we state what is obtained from the indicator function $I(\lambda, p)$.
We set $g(y; \lambda)$
\begin{equation}
g(y; \lambda) = \int_{0}^{T}e^{-\lambda^2t}f(t, y)dt
\quad(y \in \partial\Omega, \lambda \in \C_{\delta_0}). 
\label{LT-data}
\end{equation}
\begin{Thm}\label{main result no.2 for the paper}
Assume that $\partial\Omega$ and $\partial{D}$ are class of $C^2$
satisfying (I.1) and (I.2). Assume also that 
$f \in L^2((0, T)\times\partial\Omega)$, and there exists a constant 
$\beta_0 \in \R$ such that the function $g(y; \lambda)$ defined by (\ref{LT-data}) 
belongs to $C(\partial\Omega)$ for all 
$\lambda \in \C_{\delta_0}$ with large $\v{\lambda}$ and satisfies
\begin{align}
0 < \inf_{y \in \partial\Omega}
\liminf_{\v{\lambda} \to \infty}
{\rm Re}\,[\lambda^{\beta_0}g(y; \lambda)]
\leq \limsup_{\v{\lambda} \to \infty}
&\v{\lambda}^{\beta_0}
\Vert g(\cdot; \lambda) \Vert_{C(\partial\Omega)} 
< \infty 
\label{LT-data works as a source term}
\\
&\quad(\text{uniformly in } \lambda \in \C_{\delta_0}).
\nonumber
\end{align}
Then there exists a sufficiently
small $\delta_1 > 0$ such that
\begin{equation}
\lim_{\v{\lambda} \to \infty}
\frac{1}{\lambda}\log\vert I(\lambda, p) \vert = -l(p, D)
\quad(\text{uniformly in } \lambda \in \Lambda_{\delta_1}).
\label{the formula recovering the length 0}
\end{equation}
\end{Thm}
\begin{Remark}\label{the remark for main result for the paper}
\par
(1) There exist many $f \in L^2((0, T)\times\partial\Omega)$
satisfying (\ref{LT-data works as a source term}). Take $f \in C^1([0, T]; C(\partial\Omega))$
with $\inf_{y \in \partial\Omega}f(0, y) > 0$. 
As is in Remark 1 of \cite{Ikehata and Kawashita 2}, 
integration by parts implies that
$$
\V{\lambda^2g(\cdot; \lambda) - f(0, \cdot)}_{C(\partial\Omega)}
\leq \frac{\max_{0 \leq t \leq T}
\V{\partial_tf(t, \cdot)}_{C(\partial\Omega)}}
{\mu^{2}(1-\delta_0^2)}
\quad(\lambda \in \C_{\delta_0}).
$$
Note that $\delta_0$ is taken as $0 < \delta_0 < 1$.
Thus, this $f$ satisfies (\ref{LT-data works as a source term}) with $\beta_0 = 2$.
\par\noindent
(2) As is in (4) of Proposition 2 in p.1090 
of \cite{{Ikehata and Kawashita 2}}, if 
$(\xi_0, y_0) \in {\cal M}_g(p){\cup}{\cal M}_2^{+}(p)
{\cup}{\cal M}_2^{-}(p)$,
the points $p$, $\xi_0$ and $y_0$ consist of a line, and 
the point $\xi_0$ is on the line segment connecting $p$ and $y_0$.
Hence, for $(\xi_0, y_0) \in {\cal M}_2^{\pm}(p)$, there exists 
a point $(\xi_1, y_1) \in {\cal M}_2^{\mp}(p)$ respectively. 
If $\partial{D}$ itself is strictly convex, this point $ (\xi_1, y_1)$ is
uniquely determined, however, for non-strictly convex $\partial{D}$,
it is possible to be several points satisfying 
$(\xi_1, y_1) \in {\cal M}_2^{\mp}(p)$. 
In any case, ${\cal M}_2^{+}(p) = \emptyset$ if and only if
${\cal M}_2^{-}(p) = \emptyset$. 
Thus, ${\cal M}_g(p){\cup}{\cal M}_2^{+}(p){\cup}{\cal M}_2^{-}(p) = \emptyset$ if
(I.2) is assumed.
\par\noindent
(3) A sufficient condition that (I.2) holds is given in
Proposition 4 of \cite{Ikehata and Kawashita 2}. Note that
strict convexity of $\partial{D}$ does not used 
to show Proposition 4 of \cite{Ikehata and Kawashita 2}.
\par\noindent
(4) In Theorem \ref{main result no.2 for the paper}, it does not assume that
$l_p(\xi, y)$ is non-degenerate at ${\mathcal M}_1(p)$ (see (I.3) below for the precise description).
In this sense, Theorem \ref{main result no.2 for the paper} is better than the main
result in \cite{Ikehata and Kawashita 2}, since in \cite{Ikehata and Kawashita 2}, 
non-degenerate assumptions are also assumed even if
the case that $\partial{D}$ consists of one strictly convex surface.
\end{Remark}

\par

Formula (\ref{the formula recovering the length 0}) holds only for $\lambda \in \Lambda_{\delta_1}$.
This can be improved for $\lambda \in \C_{\delta_0}$ if we add the following assumption:
\begin{enumerate}
\item[(I.3)] Every point $(\xi_0, y_0) \in \partial{D}\times\partial\Omega$ attaining
$l(p, D)$ is non-degenerate critical point 
of $l_p(\xi, y)$.
\end{enumerate}
Note that as introduced in (4) of Remark \ref{the remark for main result for the paper},
(I.3) and strict convexity of $\partial{D}$ are always assumed in \cite{Ikehata and Kawashita 2}.

\begin{Thm}\label{main result for the paper}
Assume that $\partial\Omega$ and $\partial{D}$ satisfy (I.1), (I.2) and (I.3). 
Assume also that $f \in L^2((0, T)\times\partial\Omega)$, and 
there exists a constant 
$\beta_0 \in \R$ such that the function $g(y; \lambda)$ defined by (\ref{LT-data}) 
belongs to $C(\partial\Omega)$ for all 
$\lambda \in \C_{\delta_0}$ with large $\v{\lambda}$ and satisfies
(\ref{LT-data works as a source term}) for some $\beta_0 \in \R$. Further, assume that 
${\lambda}^{\beta_0}g(y; \lambda)$ is uniformly continuous in $y \in \partial\Omega$ with
respect to $\lambda \in \C_{\delta_0}$. Then, 
\begin{equation*}
\displaystyle
\lim_{\v{\lambda} \to \infty}
\frac{1}{\lambda}\log\vert I(\lambda, p) \vert = -l(p, D)
\quad(\text{uniformly in } \lambda \in \C_{\delta_0}).
\end{equation*}
\end{Thm}

\vskip1pc

\begin{Remark}\label{Remark for main result for the paper}
(1) Assumption (I.3) in Theorem \ref{main result for the paper}
is used to obtain an asymptotic behavior of $I(\lambda, p)$
as $\v{\lambda} \to \infty$ uniformly in $\lambda \in \C_{\delta_0}$.
In this sense, for non-degenerate case, we can say that 
the asymptotic behavior is better. 
\par\noindent
(2) If $\partial{D}$ and $\partial\Omega$ are $C^{2, \alpha_0}$ for some $0 < \alpha_0 < 1$, 
and $g(\cdot; \lambda) \in C^{0, \alpha_0}(\partial\Omega)$, 
$I(p, \lambda)$ has the following asymptotics:
\begin{align*}
I(\lambda, p) = \frac{1}{\lambda}e^{-\lambda l(p, D)}\left\{A(\lambda,p)g
+\Vert g(\,\cdot\, ;\lambda)\Vert_{C^{0,\alpha_0}(\partial\Omega)}O(\lambda^{-\alpha_0/2})
\right\}+O(\lambda^{-\frac{1}{2}}e^{-{\lambda}^2T})
\end{align*}
as $\vert\lambda\vert \to \infty$ uniformly with $\lambda \in \C_{\delta_0}$
for each $\delta_0>0$, where
$$\begin{array}{c}
\displaystyle
A(\lambda,p)g = \sum_{(\xi_0,y_0)\in {\mathcal M}_1(p)}C(\xi_0, y_0)H^+(\xi_0, y_0,p)g(y_0,\lambda). 
\end{array}
$$
In the above, 
$C(\xi_0, y_0)$ for each $(\xi_0, y_0)\in {\cal M}_1(p)$ is a positive constant
independent of $g$ and
\begin{align}
H^{+}(\xi, y,p)
= \frac{1}{\vert \xi - p \vert \vert \xi - y \vert}
\nu_\xi\cdot\left\{\frac{p - \xi}{\vert p - \xi \vert}
+\frac{y - \xi}{\vert y - \xi \vert}\right\},\,(\xi, y) \in \partial D\times\partial\Omega.
\label{definition of H^+}
\end{align}
This is the same asymptotic formula as in \cite{Ikehata and Kawashita 2} for the case of 
one strictly convex cavity.
Note that $(\xi, y) \in {\mathcal M}_1(p)$ means that $\nu_\xi\cdot(p - \xi) > 0$ and $\nu_y\cdot(y - \xi) > 0$, 
which yields
\begin{equation}
H^{+}(\xi, y, p) > 0 \quad((\xi, y) \in {\mathcal M}_1(p)).
\label{positivity of H^+ on M_1(p)}
\end{equation}
Thus, from (\ref{LT-data works as a source term}),  ${\rm Re} A(\lambda, p)g > 0$ holds.
\end{Remark}
\par

\par
Basic approaches for showing Theorem \ref{main result no.2 for the paper} 
and \ref{main result for the paper} are similar to our previous work \cite{Ikehata and Kawashita 2}
for the strictly convex case. As is in Section 3 of \cite{Ikehata and Kawashita 2}, 
a decomposition of $I(\lambda, p)$ and the representation formula of the main
term $I_0(\lambda, p)$ of $I(\lambda, p)$ are deduced by
using usual potential theory.  In Section \ref{Preliminaries},
a brief review of this decomposition and the formula of $I_0(\lambda, p)$ are given
(cf. Proposition \ref{integral representation of I_0}).
Note that the formula is of the form of Laplace integrals
on $\partial\Omega\times\partial{D}$ 
with exponential term $e^{-\lambda{l_p(\xi, y)}}$.

\par

\par
The amplitude functions of the Laplace integrals contain 
the inverse of the form 
$K(\lambda)(I -K(\lambda))^{-1}$ deduced from an integral operator 
$K(\lambda)$ on $\partial{D}$ with the integral kernel 
$K(\xi, \zeta; \lambda)$ estimated by
$$
\v{K(\xi, \zeta; \lambda)} \leq Ce^{-{\rm Re}\lambda\v{\xi - \zeta}}
\qquad(\xi, \zeta \in \partial{D}, \lambda \in \C_{\delta_0}).
$$
To obtain asymptotic behavior of $I_0(\lambda, p)$, it is 
crucial to get an estimate for the integral kernel 
$K^\infty(\xi, \zeta; \lambda)$ of $K(\lambda)(I -K(\lambda))^{-1}$ 
with the same exponential term $e^{-{\rm Re}\lambda\v{\xi - \zeta}}$
as for $K(\xi, \zeta; \lambda)$.
For the case $N = 1$, i.e. $\partial{D}$ is strictly convex,
such type of estimates is given in \cite{Ikehata and Kawashita}, 
and applied for an approach to an inverse problem via 
the enclosure method, which is the main subject of the previous work 
\cite{Ikehata and Kawashita 2}.

\par
For arbitrary $\partial{D}$, it seems to be hard to obtain
good estimates described above for $K^\infty(\xi, \zeta; \lambda)$.
For the case that $\partial{D}$ consists of several 
components $\partial{D}_j$ $(j = 1, 2, \ldots, N)$, however,
contributions to the estimates of the integral kernel 
$K^\infty(\xi, \zeta; \lambda)$ from the different components,
i.e. the case $\xi \in \partial{D}_j$ and $\zeta \in \partial{D}_k$ with
$j \neq k$ are weaker than the same components, i.e.
the case $\xi, \zeta \in \partial{D}_j$. In this paper, we call
the parts coming from the different components and the parts 
coming from the same components off-diagonal parts and diagonal parts 
respectively. Since the dominant part is given by
the diagonal parts, in Section \ref{Estimates of integral kernels}, 
we introduce the estimates of the integral kernels for the diagonal parts and
the amplitude functions of $I_0(\lambda, p)$. 
To control off-diagonal parts, we need to 
give additional argument, which is handled in Section \ref{eikyou-hyouka for off-diagonal parts}.

\par

In Section \ref{Proofs of the main theorems}, proofs of the main theorems are given. 
The main contributions for these Laplace integrals come from the points in 
${\mathcal M}(p)$. To pick up the main terms, we need to study on structures of
${\mathcal M}(p)$. Here, we use assumption (I.2), i.e. ${\mathcal M}_2^\pm(p) = {\mathcal M}_g(p)
= \emptyset$. By using local coordinate systems near ${\mathcal M}(p)$, eventually, the problems 
are reduced to investigating the asymptotic behaviors of Laplace integrals.
Since the appeared integrals seems to be different from usual and
typical ones, we give a brief outline to handle these integrals in 
Section \ref{Estimate of some Laplace integrals} 
for the paper to be self-contained. 

\par

To obtain the estimates of the diagonal parts, the kernel estimates 
for the case of one strictly convex cavity is essentially used. 
These estimates are given in
\cite{Ikehata and Kawashita} and \cite{Ikehata and Kawashita 2} 
by assuming that $\partial{D}$ is $C^{2, \alpha_0}$ for some $0 < \alpha_0 < 1$.
Note that this regularities assumption can be reduced to $C^2$ regularity, 
however, additional arguments are needed. This is performed in
\cite{Ikehata and Kawashita 3} by using strict convexity. Here, we can give a different
approach showing equi-continuous properties for a class of
local coordinate systems of $\partial{D}$, which is handled in Appendix.

\par

In the last of Introduction, we explain why assumption
(I.2) is needed. As is in (2) of 
Remark \ref{the remark for main result for the paper}, for
the points $(\xi_0, y_0) \in {\cal M}_g(p){\cup}{\cal M}_2^{+}(p)
{\cup}{\cal M}_2^{-}(p)$, attaining the minimum length 
$l(p, D)$, the point $\xi_0$ places on the line segment $py_0$. Hence,
contributions from the off-diagonal parts are the same levels 
as that from the diagonal parts. Thus, in this case, 
the off-diagonal parts can not be negligible.
This is essentially different from the case 
$(\xi_0, y_0) \in {\cal M}_1(p)$. Hence, the approach 
picking up the diagonal parts works only the case
$(\xi_0, y_0) \in {\cal M}_1(p)$, which is why assumption (I.2) is needed.

\setcounter{equation}{0}
\section{Decomposition of the indicator functions}\label{Preliminaries}

In the beginning, we give a remark on the class of the solutions of 
(\ref{Eq in inside for the measurement}) to make sure the meaning of 
integrals in $I(\lambda, p)$. We denote by $L^2(0, T; H)$ the space of $H$-valued $L^2$
functions in $t \in [0, T]$. 
For a Hilbert space $V$ with $V \subset L^2(\Omega\setminus\overline{D}) \subset V'$, we introduce the space
$W(0, T; V, V') = \{\, u \,\vert\, u \in L^2(0, T; V),
u' \in L^2(0, T; V')\,\}$, where $V'$ is the dual space of the 
Hilbert space $V$, and $u'$ means
the (weak) derivative in $t \in [0, T]$.  
Throughout this paper, we always assume
that the heat flux $f(t, y)$ belongs to
the space $L^2((0,\,T){\times}\partial\Omega)$.
Note that for any $ f \in L^2((0,\,T){\times}\partial\Omega)$, 
the weak solution $u \in W(0,\,T;H^1(\Omega\setminus\overline{D}), (H^1(\Omega\setminus\overline{D}))')$
of (\ref{Eq in inside for the measurement}) uniquely exists.
For the weak solutions see Section 1.5 of \cite{Ikehata and Kawashita 2}
and the references in it. Hence, the indicator function $I(\lambda, p)$ defined by (\ref{def of I(p, lambda)})
is well-defined.

\par
To show Theorem \ref{main result no.2 for the paper} and  \ref{main result for the paper}, 
we need to pick up the main term $I_0(\lambda, p)$ of 
the original indicator function $I(\lambda, p)$,
Define
$$\displaystyle
w(x; \lambda)=\int_0^T e^{-\lambda^2 t}u(t, x)dt,\,\,x\in\Omega\setminus\overline D, 
$$
which satisfies
$$
\left\{
\begin{array}{lll}
(\triangle-\lambda^2)w=u(T, x)e^{-\lambda^2T}\,\,
&\text{in}\,\Omega\setminus\overline D,\\
\big({\partial_\nu}+\rho)w=0\,\,
&\text{on}\,\partial D.
\end{array}
\right.
$$
Let $w_0(x; \lambda)$ be the solution of 
$$
\left\{
\begin{array}{lll}
(\triangle-\lambda^2)w_0 = 0\,\,
&\text{in}\,\Omega\setminus\overline D,
\\
\big({\partial_\nu}+\rho)w_0 = 0\,\,
&\text{on}\,\partial D,
\\
{\partial_\nu}w_0 = g \,\,
&\text{on}\,\partial\Omega,
\end{array}
\right.
$$
where $g(x; \lambda)$ is defined by (\ref{LT-data}).
As in Section 2 of \cite{IK0} or Appendix C of \cite{Ikehata and Kawashita 2}, 
$ w(x; \lambda)=w_0(x; \lambda)+O(e^{-\lambda^2 T})$ in $H^1(\Omega\setminus\overline{D})$
as weak sense, integration by parts implies 
\begin{equation*}
I(\lambda, p)=I_0(\lambda, p)+O(\lambda^{-\frac{1}{2}} e^{-\lambda^2 T})
\qquad(\text{as } \v{\lambda} \to \infty \text{ uniformly in } 
\lambda \in \C_{\delta_0}),
\end{equation*}
where
$$\displaystyle
I_0(\lambda, p) = \int_{\partial\Omega}\big(\partial_{\nu_y}E_{\lambda}(y,p)w_0(y; \lambda)
-E_{\lambda}(y,p)\partial_{\nu_y}w_0(y; \lambda)
\big)dS_y.
$$

\par
As is in \cite{Ikehata and Kawashita 2}, we use the layer potentials 
to construct $w_0(x; \lambda)$. From the layer potentials and 
the density functions, we can get the integral representation of 
$I_0(\lambda, p)$. The procedure is the same as in Section 3.1 
of \cite{Ikehata and Kawashita 2}. We give a brief review for it.


\par
Given $g\in C(\partial\Omega)$ and $h\in C(\partial D)$ define
\begin{align*}
\displaystyle
V_{\Omega}(\lambda)g(x)&=\int_{\partial\Omega}E_{\lambda}(x, y)g(y)dS_y,\,\,
x\in{\R}^3\setminus\partial\Omega
\\
\intertext{and}
\displaystyle
V_D(\lambda)h(x)&=\int_{\partial D}E_{\lambda}(x, \zeta)h(\zeta)dS_\zeta,\,\,
x\in{\R}^3\setminus\partial D.
\end{align*}
We construct $w_0$ in the form
\begin{equation}
w_0(x; \lambda)=V_{\Omega}(\lambda)\varphi(x; \lambda)+V_{D}(\lambda)\psi(x; \lambda), 
\label{potential expression of w_0}
\end{equation}
where $\varphi(\,\cdot\, ;\lambda)\in C(\partial\Omega)$ and 
$\psi(\,\cdot\, ;\lambda)\in C(\partial D)$ are called the density functions
satisfying the following equations in $C(\partial\Omega){\times}C(\partial{D})$:
\begin{equation}\begin{array}{c}
\displaystyle
\varphi(y; \lambda)-Y_{11}(\lambda)\varphi(y; \lambda)-Y_{12}(\lambda)\psi(y; \lambda) = g(y; \lambda)
\,\,\text{on}\,\partial\Omega, 
\\
\\
\displaystyle
\psi(\xi; \lambda)-Y_{21}(\lambda)\varphi(\xi; \lambda) - Y_{22}(\lambda)\psi(\xi; \lambda)=0
\,\,\text{on}\,\partial D.
\end{array}
\label{equations of varphi and psi}
\end{equation}
In (\ref{equations of varphi and psi}), $Y_{ij}(\lambda )$ $(i, j = 1, 2)$ are defined by
\begin{align*}
Y_{11}(\lambda)\varphi(y; \lambda) &= -\int_{\partial\Omega}\partial_{\nu_y}E_\lambda(y, z)
\varphi(z; \lambda)dS_z \quad (y \in \partial\Omega), 
\\
Y_{12}(\lambda)\psi(y; \lambda) &= -\int_{\partial{D}}\partial_{\nu_y}E_\lambda(y, \zeta)
\psi(\zeta; \lambda)dS_\zeta \quad (y \in \partial\Omega), 
\\
Y_{21}(\lambda)\varphi(\xi; \lambda) &= \int_{\partial\Omega}\big(\partial_{\nu_\xi}E_\lambda(\xi, z)
+\rho(\xi)E_\lambda(\xi, z)\big)\varphi(z; \lambda)dS_z  \quad (\xi \in \partial{D}), 
\intertext{and}
Y_{22}(\lambda)\psi(\xi; \lambda) &= \int_{\partial{D}}\big(\partial_{\nu_\xi}E_\lambda(\xi, \zeta)
+\rho(\xi)E_\lambda(\xi, \zeta)\big)
\psi(\zeta; \lambda)dS_\zeta  \quad (\xi \in \partial{D}).
\end{align*}
Note that for $\lambda \in \C$, ${\rm Re}\,\lambda > 0$, 
$Y_{11}(\lambda) \in B(C(\partial\Omega))$, $Y_{22}(\lambda) \in B(C(\partial{D}))$,
$Y_{12}(\lambda) \in B(C(\partial\Omega),C(\partial{D}))$ and 
$Y_{21}(\lambda) \in B(C(\partial{D}),C(\partial\Omega))$, 
where 
$B(E, F)$ is the set of bounded linear operators from $E$ to $F$, and $B(E) = B(E, E)$.
The operator norms of $Y_{ij}(\lambda)$ are estimated by 
\begin{align}
\V{Y_{11}(\lambda)}_{B(C(\partial\Omega))}&+\V{Y_{22}(\lambda)}_{B(C(\partial{D}))}
+ \V{Y_{12}(\lambda)}_{B(C(\partial{D}), C(\partial\Omega))}
\label{estimates of integral operators Y_{ij}}
\\&
+ \V{Y_{21}(\lambda)}_{B(C(\partial\Omega), C(\partial{D}))}
\leq C({\rm Re}\, \lambda)^{-1}
\qquad( \lambda \in \C, {\rm Re}\, \lambda > 0).
\nonumber
\end{align}
Hence, for $\lambda \in \C_{\delta_0}$ with sufficiently large ${\rm Re}\,\lambda$, 
equation (\ref{equations of varphi and psi}) can be solved by using the Neumann series. 
Since the inverse $(I - Y_{22}(\lambda))^{-1}$ is also constructed, 
from the second equation of (\ref{equations of varphi and psi}), it follows that
$\psi(\xi; \lambda) = (I - Y_{22}(\lambda))^{-1}Y_{21}(\lambda)\varphi(\xi; \lambda)$, 
which yields 
$$
\varphi(y; \lambda) 
= \{I - Y_{11}(\lambda) - Y_{12}(\lambda)(I - Y_{22})^{-1}Y_{21}(\lambda)\}^{-1}
g(y; \lambda). 
$$ 
From this and (\ref{estimates of integral operators Y_{ij}}), 
we obtain
\begin{align}
\varphi(y; \lambda) = 
g(y; \lambda)+&O(\lambda^{-1})
\Vert g(\,\cdot\, ;\lambda)\Vert_{C(\partial\Omega)}
\label{estimate of the density function varphi}
\\&
\quad(\text{uniformly in  $y \in \partial\Omega$, 
$\lambda \in \C_{\delta_0}$ as $\vert \lambda \vert \to \infty$}),
\nonumber
\end{align}
which is used in Section \ref{eikyou-hyouka for off-diagonal parts}.
\par

From (\ref{equations of varphi and psi}), 
(\ref{potential expression of w_0}) and the equality 
$$\displaystyle
I_0(\lambda, p)
=\int_{\partial D}\left(\frac{\partial E_{\lambda}}{\partial\nu}(\xi, p)
+\rho E_{\lambda}(\xi, p)\right)w_0(\xi; \lambda)dS_\xi
$$
given by integration by parts, we can write $I_0(\lambda, p)$
using only $\varphi$. This is given in Section 3.1 in \cite{Ikehata and Kawashita 2} 
for strictly convex case. Note that
this argument works even if $\partial{D}$ is not strictly convex.
Thus, we can obtain the same formula of $I_0(\lambda, p)$ as 
given in Proposition 1 of \cite{Ikehata and Kawashita 2}. 

\par

To obtain the formula of $I_0(\lambda, p)$, the transposed operator 
${}^tY_{22}(\lambda)$ of $Y_{22}(\lambda)$ is frequently used. 
Note that the operator ${}^tY_{22}(\lambda)$ is given by 
\begin{align}
{}^tY_{22}(\lambda)h(\zeta)
=\frac{1}{2\pi}
\int_{\partial D}e^{-\lambda\vert \xi - \zeta\vert}H(\xi, \zeta; \lambda)h(\xi)dS_\xi,
\,\, h \in C(\partial D), \zeta \in \partial D
\label{the transposed operator of Y_22(lambda)}
\end{align}
with the kernel $H(\xi, \zeta; \lambda) = \lambda H_0(\xi, \zeta)+H_1(\xi, \zeta)$,
where
$$\begin{array}{c}
\displaystyle
H_0(\xi, \zeta) = \frac{\nu_\xi\cdot (\zeta - \xi)}{\vert \xi - \zeta \vert^2}
\quad\text{ and }\quad
\displaystyle
H_1(\xi, \zeta) = \frac{1}{\vert \xi - \zeta \vert}
\left(\frac{\nu_\xi\cdot(\zeta - \xi)}{\vert \xi - \zeta \vert^2}
+\rho(\xi)\right).
\end{array}
$$
For $H_0(\xi, \zeta)$ and $H_1(\xi, \zeta)$, we define the operators 
$M^{(0)}(\lambda)$ and $\tilde{M}(\lambda)$ by 
\begin{align}
\displaystyle
M^{(0)}(\lambda)h(\zeta)
&=\frac{\lambda}{2\pi}\int_{\partial D}e^{-\lambda\vert \xi - \zeta\vert}H_0(\xi, \zeta)h(\xi)dS_\xi 
\label{definition of the kernel of {M}^{(0)}(lambda)}
\\
\intertext{and}
\displaystyle
\tilde{M}(\lambda)h(\zeta)
&=\frac{1}{2\pi}\int_{\partial D}e^{-\lambda\vert \xi - \zeta\vert}H_1(\xi, \zeta)h(\xi)dS_\xi,
\label{definition of the kernel of tilde{M}(lambda)}
\end{align}
respectively. 
Note that ${}^tY_{22}(\lambda)$ is decomposed into 
${}^tY_{22}(\lambda) = M^{(0)}(\lambda)+\tilde{M}(\lambda)$.

\par

Using ${}^tY_{22}(\lambda)$, we can represent $I_0(\lambda, p)$ as follows:
\begin{align*}
\displaystyle
I_0(\lambda, p) = \frac{1}{(2\pi)^2}\int_{\partial\Omega}dS_y
\varphi(y; \lambda)\int_{\partial D}
&e^{-\lambda\,l_p(\xi, y)}
\Big\{\frac{H(\xi, p; \lambda)}{\vert \xi - y\vert}
\nonumber\\
-\frac{H(\xi, y; \lambda)}{\vert \xi - p\vert}
&+2H(\xi, y; \lambda)F(\xi, p; \lambda)\Big\}dS_\xi,
\end{align*}
where
$$\displaystyle
F(\xi,p; \lambda) = e^{\lambda\vert \xi - p\vert}\left((I-{}^tY_{22}(\lambda))^{-1}
\frac{e^{-\lambda\vert\,\cdot\,-p\vert}}{\vert\,\cdot\,-p\vert}\right)(\xi).
$$
This is just (35) in p.1088 of \cite{Ikehata and Kawashita 2}.
\par

Next we decompose $F(\xi, p; \lambda)$ to pick up the main term of $I_0(\lambda, p)$.
We put $M(\lambda) = {}^tY_{22}(\lambda)(I - {}^tY_{22}(\lambda))^{-1}$ and 
$M^{(1)}(\lambda) = \tilde{M}(\lambda) + {}^tY_{22}(\lambda)M(\lambda)$.
From 
$
M(\lambda)
= {}^tY_{22}(\lambda)+({}^tY_{22}(\lambda))^2(I - {}^tY_{22}(\lambda))^{-1},
$
it follows that
\begin{equation}
M(\lambda) 
= M^{(0)}(\lambda)+M^{(1)}(\lambda), \qquad
M^{(1)}(\lambda) = \tilde{M}(\lambda)
+ {}^tY_{22}(\lambda)M(\lambda).
\label{decomposition of M(lambda)}
\end{equation}
Using these $M^{(j)}(\lambda)$, we set
\begin{align}
F^{(j)}(\xi, p; \lambda) = e^{\lambda\vert \xi - p \vert}\left(M^{(j)}(\lambda)
\left(\frac{e^{-\lambda\vert \,\cdot\,-p\vert}}{\vert\,\cdot\,-p\vert}\right)\right)(\xi)
\quad j = 0, 1.
\label{definition of F^{(j)}(xi, p; lambda)}
\end{align}
Since $(I - {}^tY_{22}(\lambda))^{-1}=I+M^{(0)}(\lambda)+M^{(1)}(\lambda)$, 
$F(\xi, p; \lambda)$ can be decomposed into
$$\displaystyle
F(\xi,p; \lambda)
= \frac{1}{\vert p - \xi \vert} + F^{(0)}(\xi, p; \lambda) 
+ F^{(1)}(\xi, p; \lambda).
$$
\par

\par
Using these notations and the function $H^{+}(\xi, y, p)$ introduced in (\ref{definition of H^+}), 
we can give an integral representation of $ I_0(\lambda,p)$.

\begin{Prop}\label{integral representation of I_0}
The decomposition
$$\displaystyle
I_0(\lambda,p)=\lambda I_{0\,0}(\lambda,p)+I_{0\,1}(\lambda,p),
$$
is valid, where
$$\begin{array}{c}
\displaystyle
G_0(\xi, y, p; \lambda)
=H^{+}(\xi, y, p)+2H_0(\xi, y)(F^{(0)}(\xi, p; \lambda)
+ F^{(1)}(\xi, p; \lambda)),\\
\\
\displaystyle
G_1(\xi, y, p; \lambda) = \frac{H_1(\xi, p)}{\vert \xi - y \vert}
+ \frac{H_1(\xi, y)}{\vert \xi - p\vert}
+ 2H_1(\xi, y)(F^{(0)}(\xi, p; \lambda) + F^{(1)}(\xi, p; \lambda))
\end{array}
$$
and
$$\displaystyle
I_{0\,j}(\lambda,p)
=\left(\frac{1}{2\pi}\right)^2
\int_{\partial\Omega}dS_y\varphi(y; \lambda)
\int_{\partial D}e^{-\lambda\,l_p(\xi,y)}G_j(\xi, y, p; \lambda)dS_\xi,\,\,j=0,1.
$$
\end{Prop}

\setcounter{equation}{0}
\section{Estimates of integral kernels}
\label{Estimates of integral kernels}

To show Theorem \ref{main result no.2 for the paper} and \ref{main result for the paper},
we need to give estimates of $I_{0\, j}(\lambda,p)$, which is reduced to getting estimates of 
$F^{(k)}(\xi, p; \lambda)$ ($k = 0, 1$) defined by (\ref{definition of F^{(j)}(xi, p; lambda)}).
In this section, necessary estimates of functions $F^{(k)}(\xi, p; \lambda)$ ($k = 0, 1$) are given. 
\par 
Since $D = \cup_{j = 1}^{N}D_j$ and each $D_j$ is disjoint, by the map 
$$
C(\partial{D}) \ni f \mapsto (f\vert_{\partial{D}_1}, 
f\vert_{\partial{D}_2}, \ldots, f\vert_{\partial{D}_N}) \in
C(\partial{D}_1)\times{\cdots\times}C(\partial{D}_N),
$$
we can identify the space $C(\partial{D})$ to 
$C(\partial{D}_1)\times{\cdots\times}C(\partial{D}_N)$.
In what follows, we put $f_j = f\vert_{\partial{D_j}}$ $(j = 1, 2, \ldots, N)$
for $f \in C(\partial{D})$. 
From (\ref{the transposed operator of Y_22(lambda)}), 
the integral kernel ${}^tY_{22}(\xi, \zeta; \lambda)$ of ${}^tY_{22}(\lambda)$ is given by
\begin{align*}
{}^tY_{22}(\xi, \zeta; \lambda) = \frac{1}{2\pi}e^{-\lambda\vert \xi - \zeta\vert}H(\zeta, \xi; \lambda).
\end{align*}
Hence, ${}^tY_{22}(\xi, \zeta; \lambda)$ is a measurable function on $\partial{D}\times\partial{D}$
with parameter $\lambda \in {\rm \C}_{\delta_0}$, and continuous for $\xi \neq \eta$. 
From the well known estimate 
\begin{align}
\vert \nu_{\xi}\cdot(\xi - \zeta ) \vert
\leq C\vert \xi - \zeta \vert^2,
\qquad(\xi, \zeta \in \partial{D} = \cup_{j = 1}^N\partial{D}_j)
\label{well known estimate for surface 1}
\end{align}
for $C^2$ surface, it follows that 
there exists a constant $C > 0$ such that
\begin{equation}
\vert {}^tY_{22}(\xi, \zeta; \lambda) \vert \leq C
\left({\rm Re }\lambda+\frac{1}{\vert \xi - \zeta \vert}\right)
e^{-{\rm Re }\lambda\vert \xi - \zeta \vert}\,\,
\quad(\xi, \zeta \in \partial D, \xi \neq \zeta, \lambda \in {\rm \C}_{\delta_0}).
\label{estimate of the kernel of Y(lambda)}
\end{equation}
For this integral kernel ${}^tY_{22}(\xi, \zeta; \lambda)$, we put
$$
{}^tY_{22}^{ij}(\xi, \zeta; \lambda) = {}^tY_{22}(\xi, \zeta; \lambda)
\qquad (\xi \in \partial{D_i}, \zeta \in \partial{D_j}, 
\lambda \in {\rm \C}_{\delta_0}),
$$
and define 
$$
{}^tY_{22}^{ij}(\lambda)f_j(\xi) = \int_{\partial{D}_j}{}^tY_{22}^{ij}(\xi, \zeta; \lambda)
f_j(\zeta)dS_\zeta.
$$
Note that for $f \in C(\partial{D})$, ${}^tY_{22}(\lambda)f(\xi)$ for $\xi \in \partial{D_i}$
can be written by 
$$
{}^tY_{22}(\lambda)f(\xi) = \sum_{j = 1}^N{}^tY_{22}^{ij}(\lambda)f_j(\xi)
\qquad(\xi \in \partial{D}_i, i = 1, 2, \ldots, N,
f \in C(\partial{D})).
$$

\par
In what follows, for simplicity, we write $\mu = \text{Re}\,\lambda$. 
For $\mu > 0$, ${}^tY_{22}(\lambda)$ is a bounded linear operator on 
$C(\partial{D})$, namely, each ${}^tY_{22}^{ij}(\lambda)$ is a bounded linear operator from
$C(\partial{D_j})$ to $C(\partial{D_i})$. From (\ref{estimates of integral operators Y_{ij}}), 
it follows that 
$$
\Vert {}^tY_{22}^{ij}(\lambda) \Vert_{B(C(\partial{D_j}), C(\partial{D_i}))} 
\leq C\mu^{-1}
\quad(i, j = 1, 2, \ldots, N)
$$
for some constant $C > 0$. 

\par

For each integral operator ${}^tY_{22}^{jj}(\lambda) \in B(C(\partial{D_j}))$ on 
$\partial{D_j}$ with the integral kernels ${}^tY_{22}^{jj}(\xi, \zeta; \lambda)$, 
${}^tY_{22}^{jj}(\lambda)(I - {}^tY_{22}^{jj}(\lambda))^{-1} \in 
B(C(\partial{D_j}))$ exists for $\lambda \in \C_{\delta_0}$, $\mu \geq \mu_0$ by 
choosing $\mu_0 > 0$ larger if necessary. In what follows, we put
$M_{D_j}(\lambda) = {}^tY_{22}^{jj}(\lambda)(I - {}^tY_{22}^{jj}(\lambda))^{-1}$. Note that
$M_{D_j}(\lambda)$ corresponds to the operator $M(\lambda)$ for the case that 
$\partial{D}$ consists of $\partial{D_j}$ only.
According to one cavity case as in \cite{Ikehata and Kawashita 2}, we define 
$M^{(0)}_{D_j}(\lambda)$, $\tilde{M}_{D_j}(\lambda)$ and $M^{(1)}_{D_j}(\lambda)$ by
\begin{align*}
\displaystyle
M^{(0)}_{D_j}(\lambda)h_j(\zeta)
&=\frac{\lambda}{2\pi}\int_{\partial D_j}e^{-\lambda\vert \xi - \zeta\vert}H_0(\xi, \zeta)h_j(\xi)dS_\xi 
\quad(\zeta \in \partial{D}_j, h_j \in C(\partial{D_j})), 
\\
\tilde{M}_{D_j}(\lambda)h_j(\zeta)
&=\frac{1}{2\pi}\int_{\partial D}e^{-\lambda\vert \xi - \zeta\vert}H_1(\xi, \zeta)h_j(\xi)dS_\xi
\quad(\zeta \in \partial{D}_j, h_j \in C(\partial{D_j}))
\end{align*}
and 
\begin{align}
M^{(1)}_{D_j}(\lambda) = \tilde{M}_{D_j}(\lambda) + {}^tY_{22}^{jj}(\lambda)M_{D_j}(\lambda).
\label{definition of M^{(1)}_{D_j}}
\end{align}
Note that the operators $M^{(0)}_{D_j}(\lambda)$, $\tilde{M}_{D_j}(\lambda)$ and $M^{(1)}_{D_j}(\lambda)$
correspond to $M^{(0)}(\lambda)$, $\tilde{M}(\lambda)$ and $M^{(1)}(\lambda)$
for one cavity case, respectively, and the relation 
$
M_{D_j}(\lambda) = M^{(0)}_{D_j}(\lambda) + M^{(1)}_{D_j}(\lambda)
$
holds.

\par
Since each $\partial{D_j}$ is strictly convex, as in 
Theorem 6.1 of \cite{Ikehata and Kawashita}, the integral kernel 
$M^{(1)}_{D_j}(\xi, \eta; \lambda)$ of $M^{(1)}_{D_j}(\lambda)$ 
has the following estimates:
\begin{Prop}\label{s-convex case estimates}
Assume that $\partial{D_j}$ is bounded, $C^{2}$ and 
strictly convex. Then, 
there exist positive constants $C$ and $\mu_0\ge 1$ such that
for all $\lambda\in {\rm \C}_{\delta_0}$ with 
$\mu=\text{Re}\,\lambda\ge\mu_0$
the operator $M^{(1)}_{D_j}(\lambda)$ has an integral kernel 
$M^{(1)}_{D_j}(\xi, \zeta; \lambda)$ which is measurable for 
$(\xi, \zeta \in \partial D\times\partial D$, 
continuous for $\xi \neq \zeta$ and has the estimate
$$\displaystyle
\vert M^{(1)}_{D_j}(\xi, \zeta;  \lambda)\vert
\le C e^{-\mu\vert \xi - \zeta \vert}
\left(1+\frac{1}{\vert \xi - \zeta \vert}
+\min\,\left\{\mu(\mu\vert \xi - \zeta \vert^3)^{1/2},\,
\frac{1}{\vert \xi - \zeta \vert^3}\right\}\right).
$$
\end{Prop}
\begin{Remark}\label{C^2 boundaries are OK}
In \cite{Ikehata and Kawashita},   
the above estimate in Proposition \ref{s-convex case estimates} is obtained for
strictly convex $\partial{D}_j$ with $C^{2, \alpha_0}$ ($0 < \alpha_0 < 1$) 
regularities. As is in \cite{Ikehata and Kawashita 3}, this regularity assumption can be relaxed to $C^2$. 
A proof of this relaxation is given in \cite{Ikehata and Kawashita 3}, which uses strict convexity of 
$\partial{D}_j$. A different proof is given in Appendix for  
the paper to be self-contained. 
\end{Remark}

\par
Note that 
since $\min\,\{\sqrt{a}, a^{-1}\}\le 1$ for all $a>0$, from
Proposition \ref{s-convex case estimates}, we get
\begin{equation}
\vert M^{(1)}_{D_j}(\xi, \zeta; \lambda)\vert
\leq C\left(\mu+\frac{1}{\vert \xi - \zeta \vert}\right)
e^{-\mu\vert \xi - \zeta \vert}.
\label{estimate of M^{(1)}_{D_j}}
\end{equation}
From (\ref{estimate of M^{(1)}_{D_j}}), (\ref{well known estimate for surface 1}) 
and the form of $M^{(0)}_{D_j}(\xi, \zeta; \lambda)$, 
we obtain
\begin{align}
\vert M_{D_j}(\xi, \zeta; \lambda)\vert &\leq 
C\left(\mu+\frac{1}{\vert \xi - \zeta \vert}\right)
e^{-\mu\vert \xi - \zeta \vert}
\label{kernel estimate for s-convex case}
\\&
\qquad(\xi, \zeta \in \partial{D}_j, \lambda \in \C_{\delta_0}, 
\mu \geq \mu_0, j = 1, 2, \ldots, N),
\nonumber
\end{align}
where $M_{D_j}(\xi, \zeta; \lambda)$ is the integral kernel of
$M_{D_j}(\lambda)$.

\par

Next, we introduce estimates of the integral kernel of $M^{(k)}(\lambda)$ ($k = 0, 1$). 
We denote by $M^{ij}(\xi, \zeta; \lambda)$, $M^{(0), ij}(\xi, \zeta; \lambda)$ 
and $M^{(1), ij}(\xi, \zeta; \lambda)$ 
($\xi \in \partial{D_i}, \zeta \in \partial{D_j}$) 
the $(i, j)$-components of the 
integral kernel of $M(\lambda)$, $M^{(0)}(\lambda)$ and $M^{(1)}(\lambda)$, respectively.
In what follows we put 
\begin{equation}
d_1 = \frac{1}{2}\min_{i \neq j}{\rm dist}(\partial{D_i}, \partial{D_j}) > 0, 
\label{definition of d_1}
\end{equation}
where
${\rm dist}(\partial{D_i}, \partial{D_j})
= \inf\{\, \vert \xi - \zeta \vert \,\vert\, \xi \in \partial{D_i}, \zeta \in \partial{D_j}\,\}$.
Note that from (\ref{estimate of the kernel of Y(lambda)}), 
for $\xi \in \partial{D_i}$, $\zeta \in \partial{D_j}$, $i \neq j$
\begin{align}
\big\vert {}^tY_{22}^{ij}(\xi, \zeta; \lambda) \big\vert 
&\leq C\Big(\mu+\frac{1}{\vert \xi - \zeta \vert}\Big)e^{-\mu\vert \xi - \zeta  \vert}
\leq C\Big(\mu+\frac{1}{2d_1}\Big)e^{-\mu\delta(2d_1)}
e^{-(1-\delta)\mu\vert \xi - \zeta  \vert}
\nonumber
\\&
\leq C\Big(\frac{1}{d_1\delta}+\frac{1}{2d_1}\Big)e^{-\mu\delta{d_1}}
e^{-(1-\delta)\mu\vert \xi - \zeta  \vert}
\label{estimate of Y_{ij}(x, y, lambda) for i neq j}
\\&\hskip20mm
\quad
(\xi \in \partial{D_i}, \zeta \in \partial{D_j}, i \neq j, 0 < \delta \leq 1).
\nonumber
\end{align}
From (\ref{definition of the kernel of {M}^{(0)}(lambda)}), 
the similar argument to getting (\ref{estimate of Y_{ij}(x, y, lambda) for i neq j})
implies that there exists a constant $C > 0$ such that
for all $i, j = 1, 2, \ldots, N$, $i \neq j$ and $0 < \delta \leq 1$, 
\begin{align}
\v{M^{(0), ij}(\xi, \zeta; \lambda)} \leq C\delta^{-1}e^{-{\delta}d_1\mu}
e^{-(1-\delta)\mu\vert \xi - \zeta \vert}
\quad(\xi \in \partial{D_i}, \zeta \in \partial{D_j}, 
\lambda \in \C_{\delta_0}, \mu \geq 1).
\label{estimate of M^{(0), ij}(xi, zeta; lambda), i neq j}
\end{align}
For the diagonal parts, (\ref{well known estimate for surface 1}) implies that
there exists a constant $C > 0$ such that
\begin{align}
\v{M^{(0), ii}(\xi, \zeta; \lambda)} \leq C{\mu}e^{-\mu\vert \xi - \zeta \vert}
\quad(\xi, \zeta \in \partial{D_i}, \xi \neq \zeta, \lambda \in \C_{\delta_0}, \mu \geq 1)
\label{estimate of M^{(0), ii}(xi, zeta; lambda)}
\end{align}
for all $i = 1, 2, \ldots, N$.

\par
The problem is to give estimates for 
$M^{(1), ij}(\xi, \zeta; \lambda)$.
\begin{Prop}\label{integral kernel of M^{(1)}}
There exist constants $C > 0 $, $\mu_1 > 0$ and $0 < \delta_1 \leq 1$ such that 
\par\noindent
(1) the integral kernel $M^{(1), ij}(\xi, \zeta; \lambda)$ is estimated by
$$
\big\vert M^{(1), ij}(\xi, \zeta; \lambda) \big\vert 
\leq C\delta^{-4}e^{-{\delta}d_1\mu}
e^{-(1-\delta)\mu\vert \xi - \zeta \vert}
\quad(\xi \in \partial{D_i}, \zeta \in \partial{D_j}, 
\lambda \in \C_{\delta_0}, \mu \geq \mu_1\delta^{-3})
$$
for all $i, j = 1, 2, \ldots, N$, $i \neq j$ and $0 < \delta \leq \delta_1$;
\par\noindent
(2) the integral kernel $M^{(1), jj}(\xi, \zeta; \lambda)$ is estimated by
$$
\big\vert M^{(1), jj}(\xi, \zeta; \lambda) 
- M^{(1)}_{D_j}(\xi, \zeta; \lambda) \big\vert 
\leq C\delta^{-4}e^{-\delta{d_1}\mu}
e^{-\mu\vert \xi - \zeta \vert}
\quad(\xi, \zeta \in \partial{D_j}, 
\lambda \in \C_{\delta_0}, \mu \geq \mu_1\delta^{-3})
$$
for all $j = 1, 2, \ldots, N$ and $0 < \delta \leq \delta_1$.
\end{Prop}
\begin{Remark}\label{remark of integral kernel of M^{(1)}}
(2) of Proposition \ref{integral kernel of M^{(1)}} and
(\ref{estimate of M^{(1)}_{D_j}}) imply that
\begin{align*}
\big\vert M^{(1), jj}(\xi, \zeta; \lambda) \big\vert 
&\leq C\Big(\mu+\delta^{-4}e^{-\delta{d_1}\mu}
+\frac{1}{\vert \xi - \zeta \vert}\Big)
e^{-\mu\vert \xi - \zeta\vert}
\\&\hskip10mm
\quad(\xi, \zeta \in \partial{D_j}, 
\lambda \in \C_{\delta_0}, \mu \geq \mu_1\delta^{-3}, 
0 < \delta \leq \delta_1).
\end{align*}
Choosing $\delta = \delta_1$ in the above, and replacing
$\mu_1$ with $ \mu_1\delta_1^{-3} $
denoted by
$\mu_1$ again, we obtain
$$
\big\vert M^{(1), jj}(\xi, \zeta; \lambda) \big\vert 
\leq C'\Big(\mu+\frac{1}{\vert \xi - \zeta \vert}\Big)
e^{-\mu\vert \xi - \zeta \vert}
\quad(\xi, \zeta \in \partial{D_j}, 
\lambda \in \C_{\delta_0}, \mu \geq \mu_1).
$$
\end{Remark}

\par
Proposition \ref{integral kernel of M^{(1)}} can be obtained by decomposing
off-diagonal parts of the integral kernels. These procedures and a proof of 
Proposition \ref{integral kernel of M^{(1)}} are given in 
Section \ref{eikyou-hyouka for off-diagonal parts}. Here, we proceed to introduce 
the estimates of $F^{(k)}(\xi, p; \lambda)$ ($k = 0, 1$) given by 
Proposition \ref{integral kernel of M^{(1)}}.

\par

For given $\varepsilon > 0$, we define
\begin{align*}
{\cal G}_{\varepsilon}(p) &= \{\xi \in\partial D \,\vert\,
\text{dist}(\xi,{\cal G}(p)) \geq \varepsilon \,\},\,\,
\quad
{\cal G}_{\varepsilon}^+(p) = {\cal G}_{\varepsilon}(p){\cap}{\cal G}^+(p).
\end{align*}
We also put
${\cal G}^{+, 0}(p) = 
\{\xi \in {\cal G}^{+}(p) \,\vert\, tp+(1-t)\xi \notin \partial{D} (0 < t \leq 1)\,\}
$
and
${\cal G}^{+, 0}_\varepsilon(p) = {\cal G}_{\varepsilon}(p){\cap}{\cal G}^{+, 0}(p)$.
The definitions (\ref{definition of F^{(j)}(xi, p; lambda)}) of $F^{(k)}(\xi, p; \lambda)$ $(k = 0, 1)$ 
in Section \ref{Preliminaries} imply that 
\begin{equation}
F^{(k)}(\xi, p; \lambda) = \sum_{j = 1}^NF^{(k), ij}(\xi, p; \lambda)
\qquad(\xi \in \partial{D_i}, i = 1, 2, \ldots, N, k = 0, 1),
\label{Form of F^{(k)}(x, p, lambda)}
\end{equation}
where 
\begin{align}
F^{(k), ij}(\xi, p; \lambda) = e^{\lambda\vert \xi - p \vert}\int_{\partial{D_j}}
M^{(k), ij}(\xi, \zeta; \lambda)\frac{e^{-\lambda\vert \zeta - p \vert}}{\vert \zeta  - p \vert}
dS_\zeta
\quad(\xi \in \partial{D_i})
\label{formula of F^{(k)}_{ij}(x, p, lambda)}.
\end{align}
To obtain the estimates of $I_{0\, j}(\lambda, p)$, the following estimates of
$F^{(k), ij}(\xi, p; \lambda)$ are necessary:
\begin{Prop}\label{estimates of $F^{(j)}$}
\par
There exists a positive constant $\mu_1$ such that 
the following assertions hold:
\par
\noindent
(1)  There exist positive constants $C$ and $d_+$ 
such that if $i, j = 1, 2, \ldots, N$, $i \neq j$, 
$\lambda \in \,{\rm \C}_{\delta_0}$,  
$\mu = \text{Re}\,\lambda\ge\mu_1\delta^{-3}$ and  
$0 < \delta \leq 1$, then
\begin{align*}
\vert F^{(k), ij}(\xi, p; \lambda)\vert &\leq C\delta^{-4}
e^{\delta{d_+}\mu}
\quad (\xi \in \partial{D_i}, i \neq j, k = 0,1).
\end{align*}
\noindent
(2) There exists a positive constant $C$ 
such that if $\lambda\in\,{\rm \C}_{\delta_0}$
and $\mu = \text{Re}\,\lambda\ge\mu_1$, then
\begin{align*}
\vert F^{(k), jj}(\xi, p; \lambda)\vert &\leq C\mu
\quad (\xi \in \partial{D_j}, k = 0,1, j = 1, 2, \ldots, N).
\end{align*}
\noindent
(3) Given $\varepsilon >0$ and an open set $U \subset \partial{D}$ satisfying 
$\overline{U} \subset {\cal G}^{+, 0}_\varepsilon(p){\cap}\partial{D}_i$ for some 
$i = 1, 2, \ldots, N$, there exist positive constants $C$, $d_2$, 
$0 < \delta_2 \leq 1$ such that if $j = 1, 2, \ldots, N$, $j \neq i$, 
$\lambda \in \,{\rm \C}_{\delta_0}$, 
$\mu = \text{Re}\,\lambda\ge\mu_1\delta^{-3}$ and 
$0 < \delta \leq \delta_2$, then
\begin{align*}
\vert F^{(k), ij}(\xi, p; \lambda)\vert &\leq C\delta^{-4}e^{-2\delta{d_2}\mu}
\quad(\xi \in \overline{U}, j \neq i, k = 0, 1).
\end{align*}
\par
\noindent
(4)  Given $\varepsilon > 0$ there exists a positive constant 
$C_{\varepsilon}$ such that if 
$\lambda\in\,{\rm \C}_{\delta_0}$ and 
$\text{Re}\,\lambda\ge\mu_0$, then
\begin{align*}
\vert F^{(k), jj}(\xi, p; \lambda)\vert &\leq C_\varepsilon\mu^{-1}
\quad (\xi \in {\cal G}_{\varepsilon}^{+}(p){\cap}\partial{D_j}, 
j = 1, 2, \ldots, N, k = 0,1).
\end{align*}
\end{Prop}

\par

To show Proposition \ref{estimates of $F^{(j)}$} and various estimates of 
the integral kernels, we need to use local coordinate systems of $\partial{D}$. 
For $a \in \R^3$ and $r > 0$, we put $B(a, r) = \{ x \in \R^3 \vert
\v{x - a} < r \}$. 
We denote by ${\mathcal B}^{\,2}(\R^2)$ the set of $C^2$ functions $f$ in $\R^2$ 
such that the norm $\V{f}_{{\mathcal B}^{\,2}(\R^2)} 
= \max_{\v{\alpha} \leq 2}\sup_{x \in \R^2}\v{\partial_x^{\alpha}f(x)}$ is finite. 
Since $\partial{D}$ is compact, we can take the following coordinate systems:

\begin{Lemma}\label{standard coordinates}
There exists $0 < r_0$ such that, for all $\xi \in \partial D$, 
$\partial D\cap B(\xi,2r_0)$ can be represented as a graph of a function on the tangent plane
of $\partial D$ at $\xi$, that is,
there exist an open neighborhood $U_\xi$ of $(0,0)$ in $\R^2$
and a function $g = g_\xi \in {\mathcal B}^{2}(\R^2)$ with $g(0,0) = 0$ and 
$\nabla g(0,0) = 0$ such that the map
$$
U_\xi \ni \,\sigma=(\sigma_1,\sigma_2)\mapsto
\xi + \sigma_1 e_1+\sigma_2e_2-g(\sigma_1,\sigma_2)\nu_\xi \in \partial D\cap B(\xi,2r_0)
$$
gives a system of local coordinates around $\xi$, where 
$\{e_1, e_2\}$ is an orthogonal basis for $T_{\xi}(\partial D)$.  Moreover
the norm $\Vert g \Vert_{{\mathcal B}^{2}(\R^2)}$ has an
upper bound independent of $\xi \in \partial D$. 
\end{Lemma}
In what follows, we call this system of coordinates the standard system 
of local coordinates around $\xi$.

As is given in Lemma 3.1 of \cite{Ikehata and Kawashita} or
Lemma 5.3 of \cite{Ikehata and Kawashita 2}, 
the following estimates, which are frequently used, 
are shown by the standard system of local coordinates:

\begin{Lemma}\label{estimates frequently used}
 Let $r_0$ be the same constant as that of Lemma \ref{standard coordinates}.
 There exists a positive constant $C$ depending only on $\partial D$ such that
\par\noindent
 (i)  for all $\xi \in \partial D$, $0 < \rho'_0 \leq r_0$, $\mu>0$, $0 \leq k < 2$
$$
\int_{B(\xi, \rho'_0)\,\cap\,\partial D}\frac{e^{-\mu\vert \xi -\zeta \vert}}{\vert \xi -\zeta \vert^k}\,dS_\zeta
\le\frac{C}{2-k}\,\min\,\{\mu^{-2+k},\,(\rho'_0)^{2-k}\};
$$
\par\noindent
(ii)  for all $\xi \in \partial D$, $\mu > 0$, $0 \leq k < 2$
$$
\int_{\partial D}\frac{e^{-\mu\vert \xi -\zeta \vert}}{\vert \xi -\zeta \vert^k}\,dS_\zeta
\le\frac{C}{2-k}\mu^{-(2-k)}\left(1+\frac{\mu^{2-k}e^{-\mu r_0}}{r_0^k}\right).
$$
\end{Lemma}
Although $C^{2, \alpha_0}$ regularities for
$\partial{D}$ is assumed in \cite{Ikehata and Kawashita}
and \cite{Ikehata and Kawashita 2}, the proofs given in 
\cite{Ikehata and Kawashita} and \cite{Ikehata and Kawashita 2}
work even if $\partial{D}$ is $C^2$. Hence, 
Lemma \ref{estimates frequently used} holds for $C^2$ boundary case.

\par
Take $\mu = 1$ and $k = 1$ in (ii) of Lemma \ref{estimates frequently used},
it follows that
$$
e^{-\inf_{\xi, \zeta \in \partial{D}}\v{\xi - \zeta}}
\int_{\partial{D_p}}\frac{dS_\zeta}{\vert \xi - \zeta \vert}
\leq 
\int_{\partial D}\frac{e^{-\vert \xi - \zeta \vert}}{\vert \xi - \zeta \vert}\,dS_\zeta
\leq C,
$$
which yields 
\begin{equation}
\int_{\partial{D_p}}\frac{dS_\zeta}{\vert \xi - \zeta \vert} \leq C
\qquad(\xi \in \partial{D}_p, p = 1, 2, \ldots, N).
\label{estimate of integral of 1/v{x - z} on partial{D}_p}
\end{equation}

Now we are in the position to give a proof of 
Proposition \ref{estimates of $F^{(j)}$} assuming Proposition \ref{integral kernel of M^{(1)}} holds.
\par\noindent
Proof of Proposition \ref{estimates of $F^{(j)}$}: 
From (\ref{estimate of M^{(0), ij}(xi, zeta; lambda), i neq j}), 
(\ref{estimate of M^{(0), ii}(xi, zeta; lambda)}), 
the estimate for $M^{(1), ij}(\xi, \zeta; \lambda)$
given in (1) of Proposition \ref{integral kernel of M^{(1)}} and the forms of
$F^{(k), ij}(\xi, p; \lambda)$ given in 
(\ref{formula of F^{(k)}_{ij}(x, p, lambda)}),  
it follows that
for any $i, j = 1, 2, \ldots, N$, $i \neq j$ and $0 < \delta \leq 1$, 
\begin{align}
\big\vert F^{(k), ij}(\xi, p; \lambda)  \big\vert &\leq
e^{\mu\vert \xi - p \vert}\int_{\partial{D_j}}C\delta^{-4}e^{-{\delta}d_1\mu}
e^{-(1-\delta)\mu\vert \xi - \zeta \vert}
\frac{e^{-\mu\vert \zeta - p \vert}}{\vert \zeta  - p \vert}dS_\zeta
\nonumber\\&
\leq \frac{C\delta^{-4}}{{\rm dist}(p, \partial{D})}
\int_{\partial{D_j}}e^{\mu(\vert \xi - p \vert-(1-\delta)\vert \xi - \zeta \vert
-\vert \zeta - p \vert)}dS_\zeta.
\label{estimate of F^{(k)}_{ij} for i neq j}
\end{align}

\par
Put $d_+ = \max\{\,\vert \xi - \zeta \vert \,\vert\, \xi, \zeta \in \partial{D}\,\} > 0$.
Noting
\begin{align*}
e^{\mu(\vert \xi - p \vert-(1-\delta)\vert \xi - \zeta \vert
-\vert \zeta - p \vert)}
\leq e^{\mu(\vert \xi - p \vert-\vert \xi - \zeta \vert
-\vert \zeta - p \vert)}e^{\mu\delta\vert \xi - \zeta \vert}
\leq e^{\mu\delta{d_+}}\quad
(\xi, \zeta \in \partial{D}),
\end{align*}
we obtain (1) of Proposition \ref{estimates of $F^{(j)}$}.

\par
When $i = j$, the second estimate in Remark \ref{remark of integral kernel of M^{(1)}}
and estimate (\ref{estimate of integral of 1/v{x - z} on partial{D}_p}) imply that
\begin{align}
\big\vert F^{(k), jj}(\xi, p; \lambda)  \big\vert &\leq
e^{\mu\vert \xi - p \vert}\int_{\partial{D_j}}C'\Big(\mu+\frac{1}{\vert \xi - \zeta \vert}\Big)
e^{-\mu\vert \xi - \zeta \vert}
\frac{e^{-\mu\vert \zeta - p \vert}}{\vert \zeta  - p \vert}dS_\zeta
\nonumber
\\&
\leq \frac{C'}{{\rm dist}(p, \partial{D})}
\int_{\partial{D_j}}e^{\mu(\vert \xi - p \vert - \vert \xi - \zeta \vert
-\vert \zeta - p \vert)}
\Big(\mu+\frac{1}{\vert \xi - \zeta \vert}\Big)dS_\zeta 
\leq C'\mu, 
\label{for estimate of F^{(k), jj}}
\end{align}
which follows (2) of Proposition \ref{estimates of $F^{(j)}$}.

\par
\par 
Next is for the proof of (3) in Proposition \ref{estimates of $F^{(j)}$}.
We need the following lemma:
\begin{Lemma}\label{L1 for estimates of F^{(k)}_{ij} i neq j}
Given $\varepsilon >0$ and an open set $U \subset \partial{D}$ satisfying 
$\overline{U} \subset {\cal G}^{+, 0}_\varepsilon(p){\cap}\partial{D}_i$ for some 
$i = 1, 2, \ldots, N$, there exist $0 < \delta_2 \leq 1$ and $d_2 > 0$ such that
$$
\vert \zeta - p \vert+(1-\delta)\vert \xi - \zeta \vert \geq
\vert \xi - p \vert + 2d_2
\quad(\xi \in \overline{U}, \zeta \in \partial{D_j}, 0 < \delta \leq \delta_2)
$$
for any $j = 1, 2, \ldots, N$ with $j \neq i$.
\end{Lemma}
Proof: From the definition of ${\cal G}^{+, 0}_\varepsilon(p)$, it follows that 
for any $\xi \in {\cal G}^{+, 0}_\varepsilon(p){\cap}\partial{D_i}$ and 
$\zeta \in \partial{D_j}$ with $j \neq i$,
$\vert p - \zeta \vert + \vert \zeta - \xi \vert > \vert p - \xi \vert$, which yields
$$
\vert p - \zeta \vert + \vert \zeta - \xi \vert > \vert p - \xi \vert
\quad(\xi \in \overline{U}, \zeta \in \cup_{j \neq i}\partial{D}_j).
$$
Since 
$\overline{U}\times\cup_{j \neq i}\partial{D_j}$
is a bounded closed set, from the above estimate, there exists a 
constant $d_2 > 0$ such that
$$
\vert p - \zeta \vert + \vert \zeta - \xi \vert \geq \vert p - \xi \vert
+3d_2
\quad(\xi \in \overline{U}, \zeta \in \partial{D_j}, i \neq j).
$$
We put 
$\varphi(\xi, \zeta, \delta) = \v{p - \zeta} +(1-\delta)\v{\zeta - \xi} - \v{p - \xi}$ and
$d_+' = \max\{\, \vert \xi - \zeta \vert \,\vert\, \xi \in \partial{D_i}, 
\zeta \in \partial{D_j}, i, j = 1, 2, \ldots, N, i \neq j \,\} > 0$.
Note that the above estimate implies that
$\varphi(\xi, \zeta, 0) = \v{p - \zeta} +\v{\zeta - \xi} - \v{p - \xi} \geq 3d_2$ 
$(\xi \in \overline{U}, \zeta \in \cup_{j \neq i}\partial{D}_j)$.
We define $\delta_2 = \min\{1, d_2/d'_+ \}$ so that 
$0 < \delta_2 \leq 1$.
\par
Noting 
\begin{align*}
\vert \varphi(\xi, \zeta, \delta) - \varphi(\xi, \zeta, 0) \vert
&= \delta\vert \xi - \zeta \vert \leq d_+'\delta
\quad(0 < \delta \leq \delta_2, (\xi, \zeta) \in \overline{U}\times\cup_{j \neq i}\partial{D_j}),
\end{align*}
we obtain 
\begin{align*}
\varphi(\xi, \zeta, \delta) &\geq \varphi(\xi, \zeta, 0) - 
\vert \varphi(\xi, \zeta, \delta) - \varphi(\xi, \zeta, 0) \vert
\geq 3d_2-d_+'\delta \geq 2d_2, 
\end{align*}
which completes the proof of Lemma \ref{L1 for estimates of F^{(k)}_{ij} i neq j}.
\hfill$\blacksquare$
\vskip1pc

\par
Take any $\varepsilon > 0$ and an open set $U \subset \partial{D}$ satisfying 
$\overline{U} \subset {\cal G}^{+, 0}_\varepsilon(p){\cap}\partial{D}_i$.  
Lemma \ref{L1 for estimates of F^{(k)}_{ij} i neq j} yields 
\begin{align*}
e^{\mu(\vert \xi - p \vert-(1-\delta)\vert \xi - \zeta \vert
-\vert \zeta - p \vert)}
\leq e^{\mu(\vert \xi - p \vert-\vert \xi - p \vert -2d_2)}
=  e^{-2\mu{d_2}}
\quad(\xi \in \overline{U}, \zeta \in \cup_{j \neq i}\partial{D_j},
0 < \delta \leq \delta_2).
\end{align*}
This estimate and (\ref{estimate of F^{(k)}_{ij} for i neq j}) imply that
\begin{align*}
\big\vert F^{(k), ij}(\xi, p; \lambda)  \big\vert 
\leq \frac{C\delta^{-4}}{{\rm dist}(p, \partial{D})}
{\rm Vol}(\partial{D_j}) e^{-2\mu{d_2}}
\quad(0 < \delta \leq \delta_2), 
\end{align*}
which shows (3) of Proposition \ref{estimates of $F^{(j)}$}.

\par
Last, we show (4) of Proposition \ref{estimates of $F^{(j)}$}. 
Since $\partial{D}_j$ is strictly convex, as in (i) of Lemma 5.2 in 
\cite{Ikehata and Kawashita 2}, p.1095, for any $\varepsilon > 0$, there
exists a constant $C_\varepsilon> 0$ such that  
\begin{align*}
\vert \xi - \zeta \vert + \vert \zeta - p \vert \geq \v{\xi - p} + C_\varepsilon\v{\zeta - \xi}
\quad(\zeta \in \partial{D}_j, \xi \in {\cal G}_\varepsilon^+(p){\cap}{\partial{D}_j}).
\end{align*}
Hence, from (ii) of Lemma \ref{estimates frequently used}, it follows that
for $\xi \in {\cal G}_\varepsilon^+(p){\cap}{\partial{D}_j}$,  
\begin{align*}
\int_{\partial{D_j}}&e^{\mu(\vert \xi - p \vert - \vert \xi - \zeta \vert
-\vert \zeta - p \vert)}\Big(\mu+\frac{1}{\vert \xi - \zeta \vert}\Big)dS_\zeta
\leq \int_{\partial{D_j}}e^{-C_\varepsilon\mu\v{\zeta - \xi}}
\Big(\mu+\frac{1}{\vert \xi - \zeta \vert}\Big)dS_\zeta
\\&
\leq C\big(\mu(C_\varepsilon\mu)^{-2}+(C_\varepsilon\mu)^{-1}\big)
= CC_\varepsilon^{-2}(1+C_\varepsilon)\mu^{-1}.
\end{align*}
The above estimate and (\ref{for estimate of F^{(k), jj}})
give (4) of Proposition \ref{estimates of $F^{(j)}$}.
\hfill$\blacksquare$
\vskip1pc
\setcounter{equation}{0}
\section{Proofs of the main theorems}
\label{Proofs of the main theorems}

For $\xi^{(0)} \in \partial{D}$, $y^{(0)} \in \partial\Omega$ and $\varepsilon > 0$, we put 
$U_{\varepsilon}(\xi^{(0)})
=\{\xi \in \partial D \,\vert\, \vert \xi - \xi^{(0)} \vert < \varepsilon \}$ 
and $V_{\varepsilon}(y^{(0)})
=\{y \in\partial\Omega\,\vert\,\vert y - y^{(0)} \vert < \varepsilon \}$.
We need the following properties of points in ${\mathcal M}(p)$:
\begin{Lemma}\label{px^{(j)} does not cross partial{D}}
Assume that $C^2$ surface $\partial{D}$ satisfies (I.1) and (I.2). Then, 
the following properties hold:
\par\noindent
(1) For each $m = 1, 2, \ldots, N$, 
${\cal G}^{+, 0}(p){\cap}\partial{D_{m}}$ is an open set in $\partial{D}$.
\par\noindent
(2) For any point $(\xi^{(0)}, y^{(0)}) \in {\cal M}(p) = {\cal M}_1(p)$, there exist constants
$\varepsilon > 0$ and $\varepsilon' > 0$, and a number $m \in \{1, 2, \ldots, N\} $ 
satisfying
$\overline{U_{2\varepsilon}(\xi^{(0)})} 
\subset \partial{D_{m}}\cap{\cal G}^{+, 0}_{\varepsilon'}(p)$.
\end{Lemma}
Proof: Note that ${\cal G}^{+}(p)$ is an open set and ${\cal G}^{+}(p)\cup{\cal G}(p)$ is
a closed set since the function $x\longmapsto (p - \xi)\cdot\nu_\xi$ is continuous. 
For (1), it suffices to show that ${\cal G}^{+, 0}(p)$ is open.
Take any $\xi_0 \in {\cal G}^{+, 0}(p)$. Then $\xi_0 \in \partial{D}_j$ for some 
$j \in \{1, 2, \ldots, N\}$.
We can assume that $ j = 1$ without loss of generality. We denote by $l[p, \xi_0]$ the 
line segment $p\xi_0$. Since $l[p, \xi_0]$ does not intersect $\cup_{j = 2}^N\partial{D}_j$, 
there exists a constant $\delta > 0$ such that 
$U_\delta\cap(\cup_{j = 2}^N\partial{D}_j) = \emptyset$,
where $U_\delta = \{ z \in \R^3 \,\vert\, {\rm dist}(z, l[p, \xi_0]) < \delta \}$. 
We can take this $\delta > 0$ small enough that
$U_\delta{\cap}\partial{D}_1 \subset {\cal G}^+(p)$ since $\xi_0 \in {\cal G}^+(p)$.
Note that $U_\delta{\cap}\partial{D}_1$ is an open set in $\partial{D}_1$.
To obtain (1), it suffices to show 
$U_\delta{\cap}\partial{D}_1 \subset {\cal G}^{+, 0}(p) $.
\par

Take any $\xi \in U_\delta{\cap}\partial{D}_1$. Since $U_\delta$ is convex set, 
$l[p, \xi] \subset U_\delta$, which yields
$l[p, \xi]\cap(\cup_{j = 2}^N\partial{D}_j) = \emptyset$. 
From $U_\delta{\cap}\partial{D}_1 \subset {\cal G}^+(p)$, $\nu_\xi\cdot(p - \xi) > 0$,
which means $l[p, \xi]\cap{D_1} = \emptyset$ since $D_1$ is convex and 
$\nu_\xi$ is the unit outer normal of $\partial{D_1}$ at $\xi$.
Thus, $l[p, \xi]\cap\partial{D} = \{\xi\}$, i.e. $\xi \in {\cal G}^{+, 0}(p)$ is shown,
which implies (1) of Lemma \ref{px^{(j)} does not cross partial{D}}.

\par
Next, we show (2). Take any point $(\xi^{(0)}, \zeta^{(0)}) \in {\cal M}(p)$. 
Since $(\xi^{(0)}, \zeta^{(0)}) \in {\cal M}_1(p)$, $\xi^{(0)} \in \partial{D}_{m}$ holds for some 
$m \in \{1, 2, \ldots, N\}$.
Since $\partial{D}_i\cap\partial{D}_m = \emptyset$ if $i \neq m$, 
this $m$ is unique. For this $m$, we show $\xi^{(0)} \in {\mathcal G}^{+, 0}(p)$. If we obtain this, 
from (1) of Lemma \ref{px^{(j)} does not cross partial{D}}, 
there exists $\varepsilon > 0$ such that $ \overline{U_{2\varepsilon}(\xi^{(0)})} 
\subset \partial{D_{m}}\cap{\cal G}^{+, 0}(p)$. Since 
${\cal G}^{+, 0}_{\varepsilon'}(p) = {\cal G}_{\varepsilon'}(p){\cap}{\cal G}^{+, 0}(p)$ 
by the definition of ${\cal G}^{+, 0}_{\varepsilon'}(p)$, we obtain
$\overline{U_{2\varepsilon}(\xi^{(0)})} 
\subset \partial{D_{m}}\cap{\cal G}^{+, 0}_{\varepsilon'}(p)$ if we choose 
$\varepsilon' > 0$ sufficiently small enough.
\par
To obtain $\xi^{(0)} \in {\mathcal G}^{+, 0}(p)$, it suffices to show that
the line segment $p\xi^{(0)}$ crosses $\partial{D}$ at only $\xi^{(0)}$.
Assume that $p\xi^{(0)}$ crosses $\partial{D}$ at 
$\zeta = tp+(1-t)\xi^{(0)} \in \partial{D}$ for some $0 < t \leq 1$. 
If $y^{(0)} $ does not contain the line ${\xi^{(0)}}p$, 
it follows that
$\vert \zeta - \xi^{(0)} \vert + \vert \xi^{(0)} - y^{(0)} \vert 
> \vert \zeta - y^{(0)} \vert $.
If not, since $(\xi^{(0)}, y^{(0)}) \in {\mathcal M}_1(p)$ means that
$ \nu_{\xi^{(0)}}\cdot(p - \xi^{(0)}) > 0 $ and 
$ \nu_{\xi^{(0)}}\cdot(y^{(0)} - \xi^{(0)}) > 0 $, 
$y^{(0)}$ is on the line segment $l[\xi^{(0)}, p]$, 
which yields
$\vert \zeta - \xi^{(0)} \vert + \vert \xi^{(0)} - y^{(0)} \vert 
> \vert \zeta - y^{(0)} \vert $. 
In any case, we obtain
\begin{align*}
l(p, D) &= 
l_p(\xi^{(0)}, y^{(0)}) = \vert p - \xi^{(0)} \vert + \vert \xi^{(0)} - y^{(0)} \vert
= \vert p - \zeta \vert + \vert \zeta - \xi^{(0)} \vert + \vert \xi^{(0)} - y^{(0)} \vert
\\&
> \vert p - \zeta \vert + \vert \zeta - y^{(0)} \vert = l_p(\zeta, y^{(0)}) 
\geq l(p, D),
\end{align*}
which is a contradiction. This completes the proof of 
Lemma \ref{px^{(j)} does not cross partial{D}}.
\hfill$\blacksquare$
\vskip1pc\noindent
\par

\par\noindent
Proof of Theorem \ref{main result no.2 for the paper}: 
From the proof of (2) of Lemma \ref{px^{(j)} does not cross partial{D}}, we can
take $\varepsilon > 0$ and $\varepsilon' > 0$ in (2) of Lemma \ref{px^{(j)} does not cross partial{D}}
arbitrary small. 
From (\ref{positivity of H^+ on M_1(p)}), 
we can also assume that 
$\inf_{(\xi, y) \in \overline{U_{2\varepsilon}(\xi^{(0)})}{\times}\overline{V_{2\varepsilon}(y^{(0)})}}
H^{+}(\xi, y, p) > 0$.
Hence, compactness of ${\cal M}(p)$ implies
that there exist points $(\xi^{(j)}, y^{(j)}) \in {\cal M}(p)$, 
numbers $m_j \in \{1, 2, \ldots, N\} $, and constants
$\varepsilon_j > 0$ and $\varepsilon_j' > 0$ $(j = 1, 2, \ldots, N_1)$ such that
$\overline{U_{2\varepsilon_j}(\xi^{(j)})} \subset \partial{D_{m_j}}\cap{\cal G}^{+, 0}_{\varepsilon_j'}(p)$, 
${\cal M}(p) \subset \cup_{j = 1}^{N_1}U_{\varepsilon_j/3}(\xi^{(j)}){\times}V_{\varepsilon_j/3}(y^{(j)})$
and 
\begin{align}
\inf_{(\xi, y) \in \overline{U_{2\varepsilon_j}(\xi^{(j)})}{\times}
\overline{V_{2\varepsilon_j}(y^{(j)})}}H^{+}(\xi, y, p) > 0.
\label{positivity of H^{+}}
\end{align}

\par
Take cut-off functions $\Psi_j\in C_0^{2}(U_{\varepsilon_j}(\xi^{(j)}){\times}V_{\varepsilon_j}(y^{(j)}))$ 
with $\Psi_j(\xi, y) = 1$ in $U_{\varepsilon_j/2}(\xi^{(j)})\times V_{\eta_j/2}(y^{(j)})$ and
$\Psi_j(\xi, y)=0$ in $(U_{2\varepsilon_j/3}(\xi^{(j)}){\times}V_{2\varepsilon_j/3}(y^{(j)}))^c$, and put
$$
I_{0kj}(\lambda,p)
=\int_{V_{\varepsilon_j}(y^{(j)})}dS_y\lambda^{\beta_0}\varphi(y; \lambda)
\int_{U_{\varepsilon_j}(\xi^{(j)})}e^{-\lambda l_p(\xi, y)}
\Psi_j(\xi, y)G_k(\xi, y, p; \lambda)dS_\xi, 
$$
and define $I_{0}^{(j)}(\lambda, p)$ by 
$I_{0}^{(j)}(\lambda, p) = I_{00j}(\lambda, p)+{\lambda}^{-1}I_{01j}(\lambda, p)$. 
Note that there exists a positive constant $c_0$ such that
$$\displaystyle
l_p(\xi, y) \geq l(p, D)+c_0\,\,\text{for}\, (\xi, y)
\in (\partial D\times\partial\Omega)\setminus\left(\cup_{j=1}^{N_1}
U_{\varepsilon_j/3}(\xi^{(j)})\times V_{\varepsilon_j/3}(y^{(j)})\right)
$$
since $l_p(\xi, y) > l(p, D)$ for all $(\xi, y) \in 
(\partial D\times\partial\Omega)\setminus\left(\cup_{j=1}^{N_1}
U_{\varepsilon_j/3}(\xi^{(j)})\times V_{\varepsilon_j/3}(y^{(j)})\right)$
being a compact set. The above estimate, (1) and (2) of 
Proposition \ref{estimates of $F^{(j)}$},
and Proposition \ref{integral representation of I_0} imply that
there exist constants $C > 0$, $\mu_1 \geq 1$, $0 < \delta_1 \leq 1$ such that
\begin{align}
\Big\vert \lambda^{\beta_0-1}I_{0}(\lambda, p) 
- \frac{1}{(2\pi)^2}\sum_{j=1}^{N_1}I_{0}^{(j)}(\lambda, p) \Big\vert
\leq Ce^{-\mu l(p, D)}
e^{-c_0\mu}(\mu+\delta^{-4}e^{\delta{d_+}\mu})\Vert g\Vert_{C(\partial\Omega)}
\label{Picking up the main terms}
\end{align}
for any $0 < \delta \leq \delta_1$, $\lambda \in \C_{\delta_0}$, $\mu \geq \delta^{-3}\mu_1$.

\par
Take local coordinates $\xi = s^{(j)}(\sigma)$ and $y = \tilde{s}^{(j)}(\tilde\sigma)$
in $U_{2\varepsilon_j}(\xi^{(j)})$ and $V_{2\varepsilon_j}(y^{(j)})$ with
$\xi^{(j)} = s^{(j)}(0)$ and $y^{(j)} = \tilde{s}^{(j)}(0)$ respectively.
We put $\tilde{\Psi}_j(\sigma, \tilde{\sigma}) = 
\Psi_j(s^{(j)}(\sigma), \tilde{s}^{(j)}(\tilde{\sigma}))$ and 
$\tilde{l_p}^{(j)}(\sigma, \tilde{\sigma}) 
= l_p(s^{(j)}(\sigma), \tilde{s}^{(j)}(\tilde{\sigma}))$,
and write $J_j(\sigma, \tilde{\sigma})$ as 
the local coordinate expressions of the surface elements. 
Using these coordinates and notations, we obtain
\begin{align*}
I_{0}^{(j)}(\lambda, p)
=\int_{\Bbb R^4} e^{-\lambda\tilde{l_p}^{(j)}(\sigma, \tilde{\sigma})}
\tilde{\Psi}_j(\sigma, \tilde{\sigma})
\alpha^{(j)}(\sigma, \tilde\sigma; \lambda)d\sigma d\tilde{\sigma},
\end{align*}
where $\alpha^{(j)}(\sigma, \tilde\sigma; \lambda)$ is defined by
\begin{align*}
\alpha^{(j)}(\sigma, \tilde\sigma; \lambda) = J_j(\sigma,\tilde{\sigma})
\lambda^{\beta_0}\varphi(\tilde{s}^{(j)}(\tilde{\sigma}); \lambda)
\big(G_0(&s^{(j)}(\sigma), \tilde{s}^{(j)}(\tilde{\sigma}), p; \lambda)
\\&
+ \lambda^{-1}G_1(s^{(j)}(\sigma), \tilde{s}^{(j)}(\tilde{\sigma}), p; \lambda)
\big).
\end{align*}

\par

Since $\overline{U_{2\varepsilon_j}(\xi^{(j)})} 
\subset \partial{D_{m_j}}\cap{\cal G}^{+, 0}_{\varepsilon_j'}(p)$ for each $j = 1, 2, \ldots, N_1$, 
(\ref{Form of F^{(k)}(x, p, lambda)}), Proposition \ref{integral representation of I_0} 
and
(3) and (4) of Proposition \ref{estimates of $F^{(j)}$} imply that
there exist constants $C > 0$, $d_2 > 0$ and $0 < \delta_2 \leq 1$ such that
\begin{align}
\vert G_0(\xi, y, p; \lambda) &- H^+(\xi, y, p) \vert 
\leq C(\mu^{-1}+\delta_2^{-4}e^{-\mu{\delta_2}d_2}),
\quad
 \vert G_1(\xi, y, p; \lambda) \vert \leq C
\label{estimates of G_0 and G_1}
\\&
((\xi, y) \in \,\overline{U_{\varepsilon_j}(\xi^{(j)})}\times\partial\Omega, 
\lambda\in\,{\rm \C}_{\delta_0} \text{ with } 
\text{Re}\,\lambda \geq \mu_1\delta_2^{-3}).
\nonumber
\end{align}
From (\ref{estimate of the density function varphi}) and (\ref{LT-data works as a source term}),
it follows that
\begin{align}
\lambda^{\beta_0}\varphi(y; \lambda) = 
\lambda^{\beta_0}g(y; \lambda)+&O(\lambda^{-1})
\quad(\text{uniformly in  $y \in \partial\Omega$, 
$\lambda \in \C_{\delta_0}$ as $\vert \lambda \vert \to \infty$}),
\label{uniformity for varphi in lambda}
\end{align}
and there exist constants $C > 0$ and $\mu_2 > 0$ such that
$$
{\rm Re}\,\, [\lambda^{\beta}g(y; \lambda)] \geq C
\quad(y \in \partial\Omega, \lambda \in \C_{\delta_0}, \mu \geq \mu_2).
$$
Combining these estimates with (\ref{positivity of H^{+}}), 
we obtain the following decomposition of $\alpha^{(j)}$: 
\begin{align*}
&\alpha^{(j)}(\sigma, \tilde\sigma; \lambda) 
= \alpha^{(j)}_1(\sigma, \tilde\sigma; \lambda)
+\lambda^{-1}\tilde\alpha^{(j)}_1(\sigma, \tilde\sigma; \lambda),
\\&
{\rm Re}\,\,[\alpha^{(j)}_1(\sigma, \tilde\sigma; \lambda)] \geq C,  
\quad
\v{\alpha^{(j)}_1(\sigma, \tilde\sigma; \lambda)}+
\v{\tilde\alpha^{(j)}_1(\sigma, \tilde\sigma; \lambda)} \leq C'
\quad((\sigma, \tilde\sigma) \in {\rm supp}\,\tilde\Psi_j)
\end{align*}
for some constants $C > 0$ and $C' > 0$. Note that there exists a constant
$ C > 0$ such that
$$ 
l(p, D) \leq \tilde{l_p}^{(j)}(\sigma, \tilde{\sigma}) 
\leq l(p, D)+C(\v{\sigma}^2+\v{\tilde\sigma}^2)
\quad((\sigma, \tilde\sigma) \in {\rm supp}\,\tilde\Psi_j)
$$
since $\tilde{l_p}^{(j)}(\sigma, \tilde{\sigma})$ is $C^2$ and
$\nabla_{(\sigma, \tilde{\sigma})}\tilde{l_p}^{(j)}(0, 0) = 0$.
\par

From these properties of $\alpha^{(j)}$ and 
$\tilde{l_p}^{(j)}(\sigma, \tilde{\sigma})$, it easily follows that
there exists a constant $C > 0$ such that
$$
\v{e^{{\lambda}l(p, D)}I_{0}^{(j)}(\lambda, p)} \leq C
\qquad (\text{uniformly in $\lambda \in \Lambda_{\delta_0}$ as } \v{\lambda} \to \infty).
$$
For lower bounds, the arguments for the Laplace integrals of some type given 
in Section 7 of \cite{Ikehata and Kawashita 3} implies that there exist constants
$\delta_1 > 0$ and $C > 0$ such that
$$
{\rm Re}\,[e^{{\lambda}l(p, D)}I_{0}^{(j)}(\lambda, p)] \geq C\mu^{-1}
\quad(\text{uniformly in $\lambda \in \Lambda_{\delta_1}$ as } \v{\lambda} \to \infty).
$$ 
Note that the Laplace integrals appeared in \cite{Ikehata and Kawashita 3} are of the cases
that the principal part of the amplitude functions, corresponding to the part 
$ \alpha^{(j)}_1$ of $\alpha^{(j)}$ for our case, does not contain the parameter $\lambda$.
Thus, the types of the integrals are slightly different from each other.
From this reason and for the paper to be self-contained, a proof for the above estimate
is given in Section \ref{Estimate of some Laplace integrals} (cf. 
Proposition \ref{lower estimate for some Laplace integral}). 
Combining (\ref{Picking up the main terms}) with the above estimates, we obtain
Theorem \ref{main result no.2 for the paper}.
\hfill$\blacksquare$
\vskip1pc\noindent
\par

\par\noindent
Proof of Theorem \ref{main result for the paper}: 
In this case, since (I.1), (I.2) and (I.3) are assumed,  
${\cal M}_2^+(p)\cup{\cal M}_2^-(p)\cup{\cal M}_g(p) = \emptyset$, 
and each point in ${\cal M}(p)$ is non-degenerate critical point 
of $l_p(\xi, y)$. These imply that ${\cal M}(p)$ is discrete set, 
which is expressed by 
$
{\cal M}(p) = {\cal M}_1(p) = \{\, (\xi^{(j)}, y^{(j)}) \,\vert\,
j = 1, 2, \ldots, N_1 \,\}.
$
From (2) of Lemma \ref{px^{(j)} does not cross partial{D}}, for 
any $j = 1, 2, \ldots, N_1$, 
there exist constants $\varepsilon_j > 0$ and $\varepsilon_j' > 0$ such that 
$\left(U_{2\varepsilon_j}(\xi^{(j)})\times V_{2\varepsilon_j}(y^{(j)})\right)\cap {\cal M}(p)
=\{(\xi^{(j)}, y^{(j)})\}$ and 
$\overline{U_{2\varepsilon_j}(\xi^{(j)})} \subset \partial{D_{m_j}}{\cap}
{\cal G}^{+, 0}_{\varepsilon_j'}(p)$. 
In this case, we can also obtain (\ref{Picking up the main terms}) and 
(\ref{estimates of G_0 and G_1}). 

\par
Taking local 
coordinates $\xi = s^{(j)}(\sigma)$ and $y = \tilde{s}^{(j)}(\tilde\sigma)$
in $U_{\varepsilon_j}(\xi^{(j)})$ and $V_{\varepsilon_j}(y^{(j)})$ with
$\xi^{(j)} = s^{(j)}(0)$ and $y^{(j)} = \tilde{s}^{(j)}(0)$ respectively, 
in this case, we decompose $I_{0}^{(j)}(\lambda, p)$ into 
$I_{0}^{(j)}(\lambda, p) = \tilde{I}_{00j}(\lambda, p)+{\lambda}^{-1}\tilde{I}_{01j}(\lambda, p)$,
where for each $j = 1, 2, \ldots, N_1$, 
$$
\tilde{I}_{0kj}(\lambda, p)
= \int_{\R^4} e^{-\lambda\tilde{l_p}^{(j)}(\sigma, \tilde{\sigma})}
\tilde{\Psi}_j(\sigma, \tilde{\sigma})
\beta^{(j)}_k(\sigma; \lambda)d\tilde{\sigma}
\quad(k = 0, 1), 
$$
$\beta^{(j)}_0(\sigma; \lambda) = 
\lambda^{\beta_0}g(\tilde{s}^{(j)}(\tilde{\sigma}); \lambda)
H^+(s^{(j)}(\sigma), \tilde{s}^{(j)}(\tilde{\sigma}), p)
J_j(\sigma, \tilde{\sigma})$,  
and $\beta^{(j)}_1(\sigma; \lambda) $ is given by 
\begin{align*}
\beta^{(j)}_1(\sigma; \lambda) &= 
\lambda^{\beta_0+1}\big(\varphi(\tilde{s}^{(j)}(\tilde{\sigma}); \lambda)
-g(\tilde{s}^{(j)}(\tilde{\sigma}); \lambda)\big)
H^+(s^{(j)}(\sigma), \tilde{s}^{(j)}(\tilde{\sigma}), p)
J_j(\sigma, \tilde{\sigma})
\\&\hskip6mm
+
\lambda^{\beta_0}\varphi(\tilde{s}^{(j)}(\tilde{\sigma}); \lambda)\big\{
\lambda\big(G_0(s^{(j)}(\sigma), \tilde{s}^{(j)}(\tilde{\sigma}), p; \lambda)
- H^+(s^{(j)}(\sigma), \tilde{s}^{(j)}(\tilde{\sigma}), p))
\\&\hskip70mm
+ G_1(s^{(j)}(\sigma), \tilde{s}^{(j)}(\tilde{\sigma}), p; \lambda)\big\}
J_j(\sigma, \tilde{\sigma}).
\end{align*}
Each $(\xi^{(j)}, y^{(j)})$ is non-degenerate, 
$\text{Hess}\,(\tilde{l_p}^{(j)})(0,0) > 0$ holds. 
Since $\lambda^{\beta_0}g(s^{(j)}(\sigma); \lambda)$ is uniformly continuous in 
$\sigma \in U_{2\varepsilon_j}(\xi^{(j)})$ 
with respect to $\lambda \in \C_{\delta_0}$, 
$\lim_{\sigma \to 0}\beta^{(j)}_0(\sigma; \lambda) = \beta^{(j)}_0(0; \lambda)$ 
uniformly in $\lambda \in \C_{\delta_0}$. From (\ref{LT-data works as a source term}),
$\beta^{(j)}_0(0; \lambda)$ is bounded for $\lambda \in \C_{\delta_0}$. 
Further, (\ref{estimates of G_0 and G_1}) and (\ref{uniformity for varphi in lambda}) yield that 
$\beta^{(j)}_1(\sigma; \lambda)$ is uniformly bounded for $\sigma \in U_{2\varepsilon_j}(\xi^{(j)})$ and 
$\lambda \in \C_{\delta_0}$.
Hence, Laplace method (cf. 
Proposition \ref{lower estimate for some Laplace integral for non-degenerate case}) implies 
$$
\tilde{I}_{01j}(\lambda,p)
= e^{-\lambda l(p, D)} 
\Vert g(\,\cdot, \,\lambda)\Vert_{C(\partial\Omega)}O(\lambda^{\beta_0-2})
$$
and
$$
\begin{array}{c}
\displaystyle
\tilde{I}_{00j}(\lambda,p)
=\frac{J_j(0,0)e^{-\lambda l(p, D)}}
{\sqrt{\text{det}\,(\text{Hess}\,(\tilde{l_p}^{(j)})(0,0))}}
\left(\frac{2\pi}{\lambda}\right)^2
\left({\lambda}^{\beta_0}g(y^{(j)},\lambda)H^+(\xi^{(j)}, y^{(j)}, p)
+ o(1)\right)
\end{array}
$$
as $\v{\lambda} \to \infty$ uniformly for $\lambda \in \C_{\delta_0}$.
From (\ref{positivity of H^+ on M_1(p)}), it follows that
$H^+(\xi^{(j)}, y^{(j)}, p) > 0 $ holds since $(\xi^{(j)},y^{(j)}) \in {\mathcal M}_1(p)$. 
This completes the proof of Theorem \ref{main result for the paper}.
\hfill$\blacksquare$
\vskip1pc\noindent
\par

Note that if $\partial{D}$ and $\partial\Omega$ are $C^{2, \alpha_0}$ for some $0 < \alpha_0 < 1$, 
and $g(\cdot; \lambda) \in C^{0, \alpha_0}(\partial\Omega)$, it holds that 
$\beta^{(j)}_k(\cdot; \lambda) \in C^{0, \alpha_0}$ near $\sigma = 0$. Hence, 
from Remark \ref{the case of Horder continuous boundary} 
in Section \ref{Estimate of some Laplace integrals}, we obtain (2) of 
Remark \ref{Remark for main result for the paper}.

\setcounter{equation}{0}
\section{The influence from the off-diagonal parts}
\label{eikyou-hyouka for off-diagonal parts}

In this section, a proof of Proposition \ref{integral kernel of M^{(1)}}
is given. As in Proposition \ref{integral kernel of M^{(1)}} and 
estimate (\ref{kernel estimate for s-convex case}), the integral kernels of 
the operators $M^{(1)}_{D_j}(\lambda)$ and 
$M_{D_j}(\lambda) = {}^tY_{22}^{jj}(\lambda)(I - {}^tY_{22}^{jj}(\lambda))^{-1}$,
which are for the case that $\partial{D}$ consists
of only one strictly convex cavity $D_j$, are given. Hence, we need to evaluate
the influences among other cavities, which is performed by decomposing 
the whole operator $(I - {}^tY_{22}(\lambda))^{-1}$ into the diagonal parts and
the off-diagonal parts. 
\par

Before giving the decomposition, 
we introduce the following estimates used frequently:
\begin{Lemma}\label{estimates frequently used 2}
There exist constants $C > 0$ and $d_1 > 0$ such that
\begin{align*}
\int_{\partial{D_p}}\left(\mu+\frac{1}{\vert \xi - \zeta \vert}\right)
e^{-\mu\delta\vert \xi - \zeta \vert}dS_\zeta
&\leq C\delta^{-1}e^{-d_1\delta\mu}
\\
\quad(\xi \in \partial{D}_i, &i, p = 1, 2, \ldots, N, i \neq p, 0 < \delta \leq 1, \mu > 0),  
\\
\int_{\partial{D_i}}\left(\mu+\frac{1}{\vert \xi - \zeta \vert}\right)
e^{-\mu\delta\vert \xi - \zeta \vert}dS_\zeta
&\leq C\delta^{-1-q}\mu^{-q}
\\
\quad(\xi \in \partial{D}_i, &i = 1, 2, \ldots, N, q = 0, 1, 0 < \delta \leq 1, \mu > 0).  
\end{align*}
\end{Lemma}
Proof: Recalling (\ref{definition of d_1}), we obtain 
 $\vert \xi - \zeta \vert \geq 2d_1$ $(\xi \in \partial{D_i}, \zeta \in \partial{D_p})$
for $i \neq p$. This implies that 
there exists a constant $C > 0$ such that
\begin{align*}
\int_{\partial{D_p}}\left(\mu+\frac{1}{\vert \xi - \zeta \vert}\right)
e^{-\mu\delta\vert \xi - \zeta \vert}dS_\zeta
&\leq \left(\mu+\frac{1}{2d_1}\right)e^{-2\mu\delta{d_1}}{\rm Vol}(\partial{D_p})
\\&
\leq C\delta^{-1}e^{-\mu\delta{d_1}}
\quad(0 < \delta \leq 1, \xi \in \partial{D}_i). 
\end{align*}
For the case $i = p$, from 
$\mu\delta{\vert \xi - \zeta \vert}e^{-\mu\delta{\vert \xi - \zeta \vert}} \leq 1$,
it follows that
$$
\int_{\partial{D_p}}\left(\mu+\frac{1}{\vert \xi - \zeta \vert}\right)
e^{-\mu\delta\vert \xi - \zeta \vert}dS_\zeta
\leq 2\delta^{-1}\int_{\partial{D_p}}\frac{dS_\zeta}{\vert \xi - \zeta \vert}
\quad(0 < \delta \leq 1).
$$
The above estimate and (\ref{estimate of integral of 1/v{x - z} on partial{D}_p})
imply the estimate for the case $q = 0$.
For the case $q = 1$, from (ii) of Lemma \ref{estimates frequently used},
for $0 < \delta \leq 1$, it follow that
\begin{align*}
\int_{\partial{D_p}}\left(\mu+\frac{1}{\vert \xi - \zeta \vert}\right)
e^{-\mu\delta\vert \xi - \zeta \vert}dS_\zeta
\leq C\big(\mu(\delta\mu)^{-2}+(\delta\mu)^{-1}\big) \leq C\delta^{-2}\mu^{-1},
\end{align*}
which completes the proof of Lemma \ref{estimates frequently used 2}.
\hfill$\blacksquare$
\vskip1pc

We put
$
Y_{D}(\lambda) = {\rm diag}({}^tY_{22}^{11}(\lambda), {}^tY_{22}^{22}(\lambda),
\cdots, {}^tY_{22}^{NN}(\lambda))$, 
where ${\rm diag}(a_1, a_2, \cdots, a_N)$ is the diagonal matrix with 
$(p, q)$-component $a_p\delta_{pq}$, and $\delta_{pq}$ is Kronecker's delta.
Note that $Y_{D}(\lambda)(I-Y_{D}(\lambda))^{-1}$ is given by
$$
Y_{D}(\lambda)(I-Y_{D}(\lambda))^{-1} 
= {\rm diag}(M_{D_1}(\lambda), M_{D_2}(\lambda), \cdots, M_{D_N}(\lambda)).
$$
\par
To handle off-diagonal parts, we introduce 
$W(\lambda) = ({}^tY_{22}(\lambda) - Y_{D}(\lambda))(I-Y_{D}(\lambda))^{-1}$ and
$\tilde{W}(\lambda) 
= ({}^tY_{22}(\lambda) - Y_{D}(\lambda))Y_{D}(\lambda)(I-Y_{D}(\lambda))^{-1}$. 
Noting 
\begin{align*}
I - {}^tY_{22}(\lambda) &= I - Y_{D}(\lambda)-({}^tY_{22}(\lambda) - Y_{D}(\lambda))
\\
&= (I-({}^tY_{22}(\lambda) - Y_{D}(\lambda))(I-Y_{D}(\lambda))^{-1})(I-Y_{D}(\lambda)),
\end{align*}
we obtain
\begin{equation}
(I- {}^tY_{22}(\lambda))(I-Y_{D}(\lambda))^{-1} = I-W(\lambda).
\label{Expression of I-Y(lambda)}
\end{equation}
\par
We define the operators $W^{ij}(\lambda)$ and $\tilde{W}^{ij}(\lambda)$ by 
$$
W(\lambda)f(\xi) = \sum_{j = 1}^NW^{ij}(\lambda)f_j(\xi),
\qquad
\tilde{W}(\lambda)f(\xi) = \sum_{j = 1}^N
\tilde{W}^{ij}(\lambda)f_j(\xi)
\qquad(\xi \in \partial{D}_i) 
$$
for $f \in C(\partial{D})$ and $i = 1, 2, \ldots, N$.
Since each $(i, j)$-component of 
$$
({}^tY_{22}(\lambda) - Y_{D}(\lambda))Y_{D}(\lambda)(I-Y_{D}(\lambda))^{-1}
$$
with $i \neq j$ and $i, j = 1, 2, \ldots, N$ is given by 
${}^tY_{22}^{ij}(\lambda)M_{D_j}(\lambda) $, 
and
\begin{align*}
W(\lambda) &= ({}^tY_{22}(\lambda) - Y_{D}(\lambda))(I-Y_{D}(\lambda))^{-1}
\\&
= ({}^tY_{22}(\lambda) - Y_{D}(\lambda))
+ ({}^tY_{22}(\lambda) - Y_{D}(\lambda))Y_{D}(\lambda)(I-Y_{D}(\lambda))^{-1},
\end{align*}
we obtain the following relations:
\begin{align*}
W^{ii}(\lambda) &= \tilde{W}^{ii}(\lambda) = 0
\quad(i = 1, 2, \ldots, N), 
\\
W^{ij}(\lambda) &= {}^tY_{22}^{ij}(\lambda)+\tilde{W}^{ij}(\lambda),
\quad
\tilde{W}^{ij}(\lambda) = {}^tY_{22}^{ij}(\lambda)M_{D_j}(\lambda)
\quad(i, j = 1, 2, \ldots, N, i \neq j). 
\end{align*}

\par
From the definition of $W(\lambda)$, $(I-W(\lambda))^{-1}$ exists 
for $\lambda \in \C_{\delta_0}$, $\mu \geq \mu_0$ by 
choosing $\mu_0 > 0$ larger if necessary. 
In what follows, we put $W^\infty(\lambda) = W(\lambda)(I-W(\lambda))^{-1}$,
which can also be written as 
$$
W^\infty(\lambda)f(\xi) 
= \sum_{j = 1}^NW^{\infty, ij}(\lambda)f_j(\xi)
\qquad(\xi \in \partial{D}_i, f \in C(\partial{D})) 
$$
by using the operators $W^{\infty, ij}(\lambda) \in B(C(\partial{D}_j), C(\partial{D}_i))$.
We denote by $W^{ij}(\xi, \zeta; \lambda)$ and $W^{\infty, ij}(\xi, \zeta; \lambda)$
the integral kernel of $W^{ij}(\lambda)$ and $W^{\infty, ij}(\lambda)$
respectively.

\par
We need the following estimates of $W^{ij}(\xi, \zeta; \lambda)$:

\begin{Prop}\label{estimate of off-diagonal part 1}
There exist constants $d_1 > 0$ and $C_1 > 0$ such that for all 
$i, j = 1, 2, \ldots, N$ with $i \neq j$ and $0 < \delta \leq 1$, 
the integral kernel $W^{ij}(\xi, \zeta; \lambda)$ is estimated by
$$
\big\vert W^{ij}(\xi, \zeta; \lambda) \big\vert 
\leq C_1\delta^{-2}e^{-{\delta}d_1\mu}e^{-(1-\delta)\mu\vert \xi - \zeta \vert}
\quad(\xi \in \partial{D_i}, \zeta \in \partial{D_j}, 
\lambda \in \C_{\delta_0}, \mu \geq \mu_0).
$$
\end{Prop}
Note that from the definition of $W^{ij}(\lambda)$, 
$W^{jj}(\xi, \zeta; \lambda) = 0$ $(\xi, \zeta \in \partial{D_j}, 
j = 1, 2, \ldots, N)$. 
\par
\par\noindent
Proof of Proposition \ref{estimate of off-diagonal part 1}: 
Assume that $i \neq j$.
Since $\tilde{W}^{ij}(\lambda) = {}^tY_{22}^{ij}(\lambda)M_{D_j}(\lambda)$, 
the integral kernel $\tilde{W}^{ij}(\xi, \zeta; \lambda)$ of $\tilde{W}^{ij}(\lambda)$
has the following integral representation:
\begin{align}
\tilde{W}^{ij}(\xi, \zeta; \lambda) = \int_{\partial{D_j}}
{}^tY_{22}^{ij}(\xi, \eta; \lambda)
M_{D_j}(\eta, \zeta; \lambda)dS_\eta
\quad(\xi \in \partial{D_i}, \zeta \in \partial{D_j}).
\label{form of tilde{W}_{ij}(x, y, lambda)}
\end{align}
The above representation, estimates (\ref{estimate of the kernel of Y(lambda)})
and (\ref{kernel estimate for s-convex case}), and 
$\vert \xi - \eta \vert \geq 2d_1$ $(\xi \in \partial{D_i}, 
\eta \in \partial{D_j})$ imply that
\begin{align*}
\big\vert \tilde{W}^{ij}(\xi, \zeta; \lambda) \big\vert
&\leq C\int_{\partial{D_j}}
\left(\mu+\frac{1}{\vert \xi - \eta \vert}\right)
e^{-\mu\vert \xi - \eta \vert}
\left(\mu+\frac{1}{\vert \eta - \zeta \vert}\right)
e^{-\mu\vert \eta - \zeta \vert}dS_\eta
\\&
\leq C\mu\int_{\partial{D_j}}
\left(\mu+\frac{1}{\vert \eta - \zeta \vert}\right)
e^{-\mu(\vert \xi - \eta \vert+\vert \eta - \zeta \vert)}dS_\eta
\\&\hskip30mm
\quad(\xi \in \partial{D_i}, \zeta \in \partial{D_j}, 
\lambda \in \C_{\delta_0}, \mu \geq \mu_0).
\end{align*}
For $\xi \in \partial{D_i}$, $\eta, \zeta \in \partial{D_j}$, it follows that
\begin{align*}
\vert \xi - \eta \vert + \vert \eta - \zeta \vert
&\geq \delta\vert \xi - \eta \vert + \delta\vert \eta - \zeta \vert
+ (1-\delta)\vert \xi - \zeta \vert
\\&
\geq 2d_1{\delta} + \delta\vert \eta - \zeta \vert
+ (1-\delta)\vert \xi - \zeta \vert.
\end{align*}
From this estimate and Lemma \ref{estimates frequently used 2}, we obtain
\begin{align*}
\int_{\partial{D_j}}
\Big(\mu&+\frac{1}{\vert \eta - \zeta \vert}\Big)
e^{-\mu(\vert \xi - \eta \vert + \vert \eta - \zeta \vert)}dS_\eta
\\&
\leq 
e^{-2d_1{\delta}\mu}e^{-(1-\delta)\mu\vert \xi - \zeta \vert}
\int_{\partial{D_j}}\Big(\mu+\frac{1}{\vert \eta - \zeta \vert}\Big)
e^{-\mu\delta\vert \eta - \zeta \vert}dS_\eta
\\&
\leq C\delta^{-2}\mu^{-1}e^{-2d_1{\delta}\mu}e^{-(1-\delta)\mu\vert \xi - \zeta \vert}
\qquad
(0 < \delta \leq 1, \xi \in \partial{D_i}, \zeta \in \partial{D_j}).
\end{align*}
Thus, we can find a constant $C > 0$ satisfying
\begin{align}
\big\vert \tilde{W}^{ij}(\xi, \zeta; \lambda) \big\vert
\leq C\delta^{-2}e^{-2d_1{\delta}\mu}e^{-(1-\delta)\mu\vert \xi - \zeta \vert}
\quad
(\xi \in \partial{D_i}, \zeta \in \partial{D_j}, i \neq j, 0 < \delta \leq 1)
\label{estimate of tilde{W}_{ij}(x, y, lambda) (i neq j)}.
\end{align}

\par
Estimate (\ref{estimate of tilde{W}_{ij}(x, y, lambda) (i neq j)})
and (\ref{estimate of the kernel of Y(lambda)}), and 
$
W^{ij}(\xi, \zeta; \lambda) = {}^tY_{22}^{ij}(\xi, \zeta; \lambda) + \tilde{W}^{ij}(\xi, \zeta; \lambda)
$
for $i \neq j$ imply Proposition \ref{estimate of off-diagonal part 1}.
\hfill$\blacksquare$
\vskip1pc\noindent
\par

From Proposition \ref{estimate of off-diagonal part 1}, we can give estimates 
of $W^{\infty, ij}(\xi, \zeta; \lambda)$.

\begin{Prop}\label{estimate of off-diagonal part 2}
There exists a constant $\mu_1 > 0$ such that for all 
$i, j = 1, 2, \ldots, N$ and $0 < \delta \leq 1$, 
the integral kernel $W^{\infty, ij}(\xi, \zeta; \lambda)$ is estimated by
$$
\big\vert W^{\infty, ij}(\xi, \zeta; \lambda) \big\vert 
\leq 2C_1\delta^{-2}e^{-{\delta}d_1\mu}
e^{-(1-\delta)\mu\vert \xi - \zeta \vert}
\quad(\xi \in \partial{D_i}, \zeta \in \partial{D_j}, 
\lambda \in \C_{\delta_0}, \mu \geq \mu_1\delta^{-3}),
$$
where $C_1 > 0$ and $d_1 > 0$ are the constants given in 
Proposition \ref{estimate of off-diagonal part 1}. Further,
there also exist constants $0 < \delta_1 \leq 1$ and
$d_3 > 0$ such that for any $j = 1, 2, \ldots, N$ and
$0 < \delta \leq \delta_1$, 
$W^{\infty, jj}(\xi, \zeta, \lambda)$ is estimated by
$$
\big\vert W^{\infty, jj}(\xi, \zeta; \lambda) \big\vert 
\leq 2C_1\delta^{-2}e^{-{\delta}d_1\mu}e^{-2d_3\mu}
e^{-\mu\vert \xi - \zeta \vert}
\quad(\xi, \zeta \in \partial{D_j}, 
\lambda \in \C_{\delta_0}, \mu \geq \mu_1\delta^{-3}).
$$
\end{Prop}

\par\noindent
Proof:
We start to getting estimates of the repeated kernels of the integral operator
$W(\lambda)$. We put 
$W^{(n)}(\lambda) = (W(\lambda))^n$ $(n = 1, 2, \ldots)$, and denote 
by $W^{(n)}_{ij}(\lambda)$ the $(i, j)$-components of $W^{(n)}(\lambda)$,
and by $W^{(n)}_{ij}(\xi, \zeta; \lambda)$ ($\xi \in \partial{D_i}$, 
$\zeta \in \partial{D_j}$) the integral kernel of $W^{(n)}_{ij}(\lambda)$.
Then it follows that

$$
W^{(n)}_{ij}(\lambda)f_j(\xi) = \int_{\partial{D_j}}W^{(n)}_{ij}(\xi, \zeta; \lambda)f_j(\zeta)dS_\zeta
\qquad(\xi \in \partial{D_i}).
$$
By induction, we show
\begin{align}
\big\vert W^{(n)}_{ij}(\xi, \zeta; \lambda) \big\vert 
\leq C_2^{n-1}&(C_1\delta^{-2}e^{-\mu{\delta}d_1})^ne^{-(1-\delta)\mu\vert \xi - \zeta \vert}
\label{estimate of repeated kernel for off-diagonal part}
\\&\hskip-10mm
(\xi \in \partial{D_i}, \zeta \in \partial{D_j}, 0 < \delta \leq 1, i, j = 1, 2, \ldots, N,
n = 1, 2, \ldots),
\nonumber
\end{align}
where $C_2 = {\rm Vol}(\partial{D}) = \sum_{j = 1}^N{\rm Vol}(\partial{D_j}) > 0$.
From $W^{(1)}_{ij}(\xi, \zeta; \lambda) = W_{ij}(\xi, \zeta; \lambda)$, 
Proposition \ref{estimate of off-diagonal part 1} shows that the case $n = 1$ is true.
Assume that the case less than or equal to $n$ are true.
\par
Note that the kernel $W^{(n+1)}_{ij}(\xi, \zeta; \lambda)$ is given by
\begin{equation}
W^{(n+1)}_{ij}(\xi, \zeta; \lambda)
= \sum_{p = 1}^N\int_{\partial{D_p}}W_{ip}(\xi, \eta; \lambda)W^{(n)}_{pj}(\eta, \zeta; \lambda)dS_\eta.
\label{def of repeated kernels for off-diagonal part}
\end{equation}
Hence Proposition \ref{estimate of off-diagonal part 1} and the assumption of induction 
imply that
\begin{align*}
\big\vert W^{(n+1)}_{ij}(\xi, \zeta; \lambda) \big\vert
\leq C_2^{n-1}(C_1\delta^{-2}e^{-\mu{\delta}d_1})^{n+1}
\sum_{p = 1}^N\int_{\partial{D_p}}
e^{-(1-\delta)\mu(\vert \xi - \eta \vert + \vert \eta - \zeta \vert)}dS_\eta.
\end{align*}
From $e^{-(1-\delta)\mu(\vert \xi - \eta \vert + \vert \eta - \zeta \vert)}
\leq e^{-(1-\delta)\mu\vert \xi - \zeta \vert}$, 
it follows that 
$$
\int_{\partial{D_p}}
e^{-(1-\delta)\mu(\vert \xi - \eta \vert + \vert \eta - \zeta \vert)}dS_\eta
\leq e^{-(1-\delta)\mu\vert \xi - \zeta \vert}{\rm Vol}(\partial{D_p}).
$$
This implies
$$
\big\vert W^{(n+1)}_{ij}(\xi, \zeta; \lambda) \big\vert \leq
C_2^{n}(C_1\delta^{-2}e^{-\mu{\delta}d_1})^{n+1}e^{-(1-\delta)\mu\vert \xi - \zeta\vert}, 
$$
which means that the case $n+1$ is also true.
Thus, we obtain (\ref{estimate of repeated kernel for off-diagonal part}).
\par
For handling the diagonal parts $W^{(n)}_{jj}(\xi, \zeta, \lambda)$, we need 
the following lemma:
\begin{Lemma}\label{L1 for diagonal part of the repeated kernels for off-diagonal part}
There exist $0 < \delta_1 \leq 1$ and $d_3 > 0$ such that
$$
(1-\delta)(\vert \xi - \eta \vert+\vert \eta - \zeta \vert) \geq
\vert \xi - \zeta \vert + 2d_3
\quad(\xi, \zeta \in \partial{D_j}, \eta \in \partial{D_p}, j \neq p, 0 < \delta \leq \delta_1).
$$
\end{Lemma}
Proof: We put $\varphi(\xi, \zeta, \eta, \delta) = (1-\delta)(\vert \xi - \eta \vert+\vert \eta - \zeta \vert)
- \vert \xi - \zeta \vert$. The function $\varphi(\xi, \zeta, \eta, 0)$ is continuous on the compact set
$\cup_{j = 1}^N\partial{D_j}\times\partial{D_j}\times
(\partial{D}\setminus\partial{D_j})$, and $\varphi(\xi, \zeta, \eta, 0) > 0$ $((\xi, \zeta, \eta) \in 
\cup_{j = 1}^N\partial{D_j}\times\partial{D_j}\times(\partial{D}\setminus\partial{D_j}))$,
which yields
$$
d_3 = 3^{-1}\inf\{\, \varphi(\xi, \zeta, \eta, 0) \,\vert\, 
(\xi, \zeta, \eta) \in \cup_{j = 1}^N
\partial{D_j}\times\partial{D_j}\times(\partial{D}\setminus\partial{D_j})\,\} > 0.
$$
We put $d_+' = \max\{\, \vert \xi - \zeta \vert \,\vert\, \xi \in \partial{D_i}, 
\zeta \in \partial{D_j}, i, j = 1, 2, \ldots, N, i \neq j \,\} > 0$.
Note that
\begin{align*}
\vert \varphi(\xi, \zeta, \eta, \delta) - \varphi(\xi, \zeta, \eta, 0) \vert
&= \delta(\vert \xi - \eta \vert+\vert \eta - \zeta \vert) \leq 2d_+'\delta
\\&
\quad((\xi, \zeta, \eta) \in 
\cup_{j = 1}^N\partial{D_j}\times\partial{D_j}\times(\partial{D}\setminus\partial{D_j})),
\end{align*}
we put $\delta_1 = \min\{1, (2d_+')^{-1}d_3 \}$. Then $0 < \delta_1 \leq 1$, and 
for $0 < \delta \leq \delta_1$
and $(\xi, \zeta, \eta) \in 
\cup_{j = 1}^N\partial{D_j}\times\partial{D_j}\times(\partial{D}\setminus\partial{D_j})$,
\begin{align*}
\varphi(\xi, \zeta, \eta, \delta) &\geq \varphi(\xi, \zeta, \eta, 0) - 
\vert \varphi(\xi, \zeta, \eta, \delta) - \varphi(\xi, \zeta, \eta, 0) \vert
\geq 3d_2-2d_+'\delta \geq 2d_3, 
\end{align*}
which completes the proof of 
Lemma \ref{L1 for diagonal part of the repeated kernels for off-diagonal part}.
\hfill$\blacksquare$
\vskip1pc\noindent
\par

Now estimates of $W^{(n)}_{jj}(\xi, \zeta; \lambda)$ are given as follows:
Noting $W^{jj}(\xi, \zeta; \lambda) = 0$ $(j = 1, 2, \ldots, N)$ and 
(\ref{def of repeated kernels for off-diagonal part}), for $ n \geq 2$
we have 
\begin{align*}
W^{(n)}_{jj}(\xi, \zeta; \lambda)
= \sum_{p \neq j}\int_{\partial{D_p}}
W_{jp}(\xi, \eta; \lambda)W^{(n-1)}_{pj}(\eta, \zeta; \lambda)dS_\eta.
\end{align*}
The above equality, Proposition \ref{estimate of off-diagonal part 1}
and (\ref{estimate of repeated kernel for off-diagonal part}) imply that
for any $\xi, \zeta \in \partial{D_j}$, 
\begin{align*}
\big\vert W^{(n)}_{jj}(\xi, \zeta; \lambda) \big\vert 
\leq C_2^{n-2}(C_1\delta^{-2}e^{-\mu{\delta}d_1})^{n}
\sum_{p \neq j}\int_{\partial{D_p}}
e^{-(1-\delta)\mu(\vert \xi - \eta \vert + \vert \eta - \zeta \vert)}
dS_\eta.
\end{align*}
Since Lemma \ref{L1 for diagonal part of the repeated kernels for off-diagonal part}
yields that there exist $0 < \delta_1 \leq 1$ and $d_2 > 0$ such that for any
$0 < \delta \leq \delta_1$ and $j \neq p$,
$$
(1-\delta)(\vert \xi - \eta \vert+\vert \eta - \zeta \vert) \geq
\vert \xi - \zeta \vert + 2d_3
\quad(\xi, \zeta \in \partial{D_j}, \eta \in \partial{D_p}, 0 < \delta \leq \delta_1), 
$$
which yields
$$
\int_{\partial{D_p}}
e^{-(1-\delta)\mu(\vert \xi - \eta \vert + \vert \eta - \zeta \vert)}
dS_\eta
\leq {\rm Vol}(\partial{D_p})
e^{-\mu\vert \xi - \zeta \vert}e^{-2\mu{d_3}}
\qquad(\xi, \zeta \in \partial{D_j}, p \neq j).
$$
From these estimates, we obtain 
\begin{align}
\big\vert W^{(n)}_{jj}(\xi, \zeta; \lambda) \big\vert 
&\leq C_2^{n-1}(C_1\delta^{-2}e^{-\mu{\delta}d_1})^{n}
e^{-\mu\vert \xi - \zeta \vert}e^{-2\mu{d_3}}
\label{estimate of repeated kernel for diagonal part}
\\&\hskip10mm
(\xi, \zeta \in \partial{D_j}, 0 < \delta \leq \delta_1, 
j = 1, 2, \ldots, N).
\nonumber
\end{align}

\par

Now we put $\mu_1 = \max\{\mu_0, 2C_1C_2/d_1\,\} \geq \mu_0$.
For $\mu \geq \delta^{-3}\mu_1$, it follows that
$$
C_1C_2\delta^{-2}e^{-\mu{\delta}d_1}
\leq C_1C_2\delta^{-2}(\mu{\delta}d_1)^{-1}
(\mu{\delta}d_1)e^{-\mu{\delta}d_1}
\leq C_1C_2d_1^{-1}\delta^{-3}\mu^{-1}
\leq 1/2.
$$
This estimate, (\ref{estimate of repeated kernel for off-diagonal part})
and (\ref{estimate of repeated kernel for diagonal part}) imply
\begin{align*}
\big\vert W^{(n)}_{ij}(\xi, \zeta; \lambda) \big\vert &\leq
C_1\delta^{-2}e^{-\mu{\delta}d_1}\Big(\frac{1}{2}\Big)^{n-1}
e^{-(1-\delta)\mu\vert \xi - \zeta \vert}
\\&\hskip30mm
\quad(\xi \in \partial{D_i}, \zeta \in \partial{D_j}, 
0 < \delta \leq 1), 
\\
\big\vert W^{(n)}_{jj}(\xi, \zeta; \lambda) \big\vert &\leq
C_1\delta^{-2}e^{-\mu{\delta}d_1}\Big(\frac{1}{2}\Big)^{n-1}
e^{-\mu\vert \xi - \zeta \vert}e^{-2\mu{d_3}}
\quad(\xi, \zeta \in \partial{D_j}, 0 < \delta \leq \delta_1). 
\end{align*}
Noting that $W^{\infty, ij}(\xi, \zeta; \lambda) 
= \sum_{n = 1}^{\infty}W^{(n)}_{ij}(\xi, \zeta; \lambda)$, 
since $W^{\infty, ij}(\lambda) 
= \sum_{n = 1}^{\infty}W^{(n)}_{ij}(\lambda)$, 
we obtain
\begin{align*}
\big\vert W^{\infty, ij}(\xi, \zeta; \lambda) \big\vert &\leq
2C_1\delta^{-2}e^{-\mu{\delta}d_1}
e^{-(1-\delta)\mu\vert \xi - \zeta \vert}
\quad(\xi \in \partial{D_i}, \zeta \in \partial{D_j}, 
0 < \delta \leq 1), 
\\
\big\vert W^{\infty, jj}(\xi, \zeta; \lambda) \big\vert &\leq
2C_1\delta^{-2}e^{-\mu{\delta}d_1}
e^{-\mu\vert \xi - \zeta \vert}e^{-2\mu{d_3}}
\quad(\xi, \zeta \in \partial{D_j}, 0 < \delta \leq \delta_1), 
\end{align*}
which completes the proof of Proposition \ref{estimate of off-diagonal part 2}.
\hfill$\blacksquare$
\vskip1pc\noindent
\par

Now, we proceed to get
estimates for the integral kernel of $M(\lambda) = {}^tY_{22}(\lambda)(I - {}^tY_{22}(\lambda))^{-1}$.
From (\ref{Expression of I-Y(lambda)}), it follows that
\begin{align*}
(I - {}^tY_{22}(\lambda))^{-1} 
&= (I - Y_{D}(\lambda))^{-1}(I - W(\lambda))^{-1}
\\&
= I + Y_{D}(\lambda)(I - Y_{D}(\lambda))^{-1}
+ W(\lambda)(I-W(\lambda))^{-1}
\\&\hskip20mm
+ Y_{D}(\lambda)(I - Y_{D}(\lambda))^{-1}W(\lambda)(I - W(\lambda))^{-1},
\end{align*}
which yields
\begin{align*}
M(\lambda) = Y_{D}(\lambda)(I - Y_{D}(\lambda))^{-1}
&+ W(\lambda)(I - W(\lambda))^{-1}
\\&
+ Y_{D}(\lambda)(I - Y_{D}(\lambda))^{-1}W(\lambda)(I - W(\lambda))^{-1}
\end{align*}
since $M(\lambda) = {}^tY_{22}(\lambda)(I - {}^tY_{22}(\lambda))^{-1} = (I - {}^tY_{22}(\lambda))^{-1} - I$. 
We denote by $M^{ij}(\xi, \zeta; \lambda)$ the $(i, j)$-components of the
integral kernel of $M(\lambda)$. The above expression implies that
\begin{align}
M^{ij}(\xi, \zeta; \lambda) &= \delta_{ij}M_{D_j}(\xi, \zeta; \lambda)
+ W^{\infty, ij}(\xi, \zeta; \lambda)
\nonumber
\\&\hskip15mm
+ \int_{\partial{D_i}}
M_{D_i}(\xi, \eta; \lambda)
W^{\infty, ij}(\eta, \zeta; \lambda)dS_\eta
\quad(\xi \in \partial{D_i}, \zeta \in \partial{D_j})
\label{representation of M_{ij}(x, y, lambda)}.
\end{align}
From (\ref{decomposition of M(lambda)}) and 
(\ref{definition of the kernel of tilde{M}(lambda)}), 
for $\xi \in \partial{D_i}, \zeta \in \partial{D_j}$, 
the $(i, j)$-components of the integral kernel $M^{(1), ij}(\xi, \zeta; \lambda)$ of $M^{(1)}(\lambda)$
is given by 
\begin{align}
M^{(1), ij}(\xi, \zeta; \lambda) =\frac{1}{2\pi}e^{-\lambda\vert \xi - \zeta \vert}H_1(\xi, \zeta)
+ \sum_{p = 1}^N\int_{\partial{D_p}}
{}^tY_{22}^{ip}(\xi, \eta; \lambda)
M^{pj}(\eta, \zeta; \lambda)dS_\eta.
\label{representation of M^{(1)}_{ij}(x, y, lambda)}
\end{align}
Since $M^{(1)}_{D_j}$ is defined by (\ref{definition of M^{(1)}_{D_j}}), 
$M^{(1)}_{D_j}(\xi, \zeta; \lambda)$ are written by 
\begin{align}
M^{(1)}_{D_j}(\xi, \zeta; \lambda) = \frac{1}{2\pi}e^{-\lambda\vert \xi - \zeta \vert}H_1(\xi, \zeta)
+ \int_{\partial{D_j}}{}^tY_{22}^{jj}(\xi, \eta; \lambda)M_{D_j}(\eta, \zeta; \lambda)dS_\eta
\label{representation of M^{(1)}_{D_j}(x, y, lambda)}. 
\end{align}
\begin{Lemma}\label{integral kernel of Y(lambda)(I-Y(lambda))^{-1}}
There exist constants $C > 0 $ and $\mu_1 > 0$ such that for all 
$0 < \delta \leq 1$ and $i, j = 1, 2, \ldots, N$ with $i \neq j$, 
the integral kernel $M^{ij}(\xi, \zeta; \lambda)$ given by (\ref{representation of M_{ij}(x, y, lambda)}) 
is estimated by
$$
\big\vert M^{ij}(\xi, \zeta; \lambda) \big\vert 
\leq C\delta^{-3}e^{-{\delta}d_1\mu}
e^{-(1-\delta)\mu\vert \xi - \zeta \vert}
\quad(\xi \in \partial{D_i}, \zeta \in \partial{D_j}, 
\lambda \in \C_{\delta_0}, \mu \geq \mu_1\delta^{-3}).
$$
There also exist constants $C > 0 $, $\mu_1 > 0$ and 
$0 < \delta_1 \leq 1$ such that for all 
$j = 1, 2, \ldots, N$  and $0 < \delta \leq \delta_1$, 
the integral kernel $M^{jj}(\xi, \zeta; \lambda)$ is estimated by
\begin{align*}
\big\vert M^{jj}(\xi, \zeta; \lambda) 
- M_{D_j}(\xi, \zeta; \lambda) \big\vert 
&\leq C\delta^{-2}e^{-\delta{d_1}\mu}e^{-{d_3}\mu}
e^{-\mu\vert \xi - \zeta \vert}
\\&
\quad(\xi, \zeta \in \partial{D_j}, 
\lambda \in \C_{\delta_0}, \mu \geq \mu_1\delta^{-3}).
\end{align*}
\end{Lemma}
\begin{Remark}\label{remark of integral kernel of Y(lambda)(I-Y(lambda))^{-1}}
From Lemma \ref{integral kernel of Y(lambda)(I-Y(lambda))^{-1}} and
(\ref{kernel estimate for s-convex case}), we obtain
\begin{align*}
\big\vert M^{jj}(\xi, \zeta; \lambda) \big\vert 
&\leq C\Big(\mu+\delta^{-2}e^{-\delta{d_1}\mu}
+\frac{1}{\vert \xi - \zeta \vert}\Big)
e^{-\mu\vert \xi - \zeta \vert}
\\&\hskip10mm
\quad(\xi, \zeta \in \partial{D_j}, 
\lambda \in \C_{\delta_0}, \mu \geq \mu_1\delta^{-3}, 
0 < \delta \leq \delta_1),
\end{align*}
where $0 < \delta_1 \leq 1$ and $\mu_1 > 0$ are given in 
Lemma \ref{integral kernel of Y(lambda)(I-Y(lambda))^{-1}}.
\end{Remark}
Proof: We show the first estimate in Lemma \ref{integral kernel of Y(lambda)(I-Y(lambda))^{-1}}.
Assume that $i \neq j$. It suffices to show 
\begin{align}
\int_{\partial{D_i}}\big\vert
M_{D_i}(\xi, \eta; \lambda)&
W^{\infty, ij}(\eta, \zeta; \lambda)\big\vert
dS_\eta \leq 
C\delta^{-3}e^{-{\delta}d_1\mu}
e^{-(1-\delta)\mu\vert \xi - \zeta \vert}
\label{for estimates of M_{ij} no 1}
\\&\hskip20mm
\quad(\xi \in \partial{D_i}, \zeta \in \partial{D_j}, i \neq j,  
\lambda \in \C_{\delta_0}, \mu \geq \mu_1\delta^{-3})
\nonumber
\end{align}
since other terms in the representation 
(\ref{representation of M_{ij}(x, y, lambda)}) of
the integral kernel $M^{ij}(\xi, \zeta; \lambda)$ are given by 
Proposition \ref{estimate of off-diagonal part 2}.

\par
Keeping the case $i \neq j$ in mind,  
and using 
(\ref{kernel estimate for s-convex case}) and 
Proposition \ref{estimate of off-diagonal part 2}, for $0 < \delta \leq 1$, 
we obtain 
\begin{align*}
\int_{\partial{D_i}}\big\vert
&M_{D_i}(\xi, \eta; \lambda)
W^{\infty, ij}(\eta, \zeta; \lambda)\big\vert
dS_z 
\\&
\leq 2CC_1\delta^{-2}e^{-{\delta}d_1\mu}
\int_{\partial{D_i}}\Big(\mu+\frac{1}{\vert \xi - \eta \vert}\Big)
e^{-\mu\vert \xi - \eta \vert}
e^{-(1-\delta)\mu\vert \eta - \zeta \vert}dS_\eta.
\end{align*}
Since 
\begin{align*}
e^{-\mu\vert \xi - \eta \vert}e^{-(1-\delta)\mu\vert \eta - \zeta \vert}
\leq e^{-(1-\delta)\mu\vert \xi - \zeta \vert}e^{-\mu\delta\vert \xi - \eta \vert}
\quad(\xi \in \partial{D_i}, \zeta, \eta \in \partial{D_j}),
\end{align*}
$0 < \delta \leq 1$ and Lemma \ref{estimates frequently used 2} imply
\begin{align*}
\int_{\partial{D_i}}\big\vert
&M_{D_i}(\xi, \eta; \lambda)
W^{\infty, ij}(\eta, \zeta; \lambda)\big\vert
dS_\eta 
\\&
\leq 2CC_1\delta^{-2}e^{-{\delta}d_1\mu}
e^{-(1-\delta)\mu\vert \xi - \zeta \vert}
\int_{\partial{D_i}}\Big(\mu+\frac{1}{\vert \xi - \eta \vert}\Big)
e^{-\delta\mu\vert \xi - \eta \vert}dS_\eta
\\&
\leq C\delta^{-3}e^{-{\delta}d_1\mu}
e^{-(1-\delta)\mu\vert \xi - \zeta \vert},
\end{align*}
which shows (\ref{for estimates of M_{ij} no 1}).

\par
For the case $i = j$, 
Proposition \ref{estimate of off-diagonal part 2} and 
(\ref{representation of M_{ij}(x, y, lambda)}), it suffices to show
\begin{align}
\int_{\partial{D_j}}\big\vert
M_{D_j}(\xi, \eta; \lambda)&
W^{\infty, jj}(\eta, \zeta; \lambda)\big\vert
dS_\eta \leq 
C\delta^{-2}e^{-{\delta}d_1\mu}e^{-d_3\mu}
e^{-\mu\vert \xi - \zeta \vert}
\label{for estimates of M_{ij} no 2}
\\&\hskip20mm
\quad(\xi, \zeta \in \partial{D_j}, 
\lambda \in \C_{\delta_0}, \mu \geq \mu_1\delta^{-3},
0 < \delta \leq \delta_1),
\nonumber
\end{align}
where $0 < \delta_1 \leq 1$ is the constant given in 
Proposition \ref{estimate of off-diagonal part 2}.

\par
From the estimate of $W^{\infty, jj}(\xi, \zeta; \lambda)$
in Proposition \ref{estimate of off-diagonal part 2} and 
(\ref{kernel estimate for s-convex case}), it follows that for $0 < \delta_1 \leq 1$, 
\begin{align*}
\int_{\partial{D_j}}\big\vert
&M_{D_j}(\xi, \eta; \lambda)
W^{\infty, jj}(\eta, \zeta; \lambda)\big\vert
dS_\eta 
\\&
\leq 2CC_1\delta^{-2}e^{-{\delta}d_1\mu}e^{-2d_3\mu}
\int_{\partial{D_j}}\Big(\mu+\frac{1}{\vert \xi - \eta \vert}\Big)
e^{-\mu\vert \xi - \eta \vert}
e^{-\mu\vert \eta - \zeta \vert}dS_\eta.
\end{align*}
Since (\ref{estimate of integral of 1/v{x - z} on partial{D}_p}) implies 
$$
e^{-d_3\mu}\int_{\partial{D_j}}\Big(\mu+\frac{1}{\vert \xi - \eta \vert}\Big)dS_\eta
\leq Ce^{-d_3\mu}(\mu+1) \leq C
\quad(\mu \geq 1),
$$
we obtain (\ref{for estimates of M_{ij} no 2}), which completes the proof of
Lemma \ref{integral kernel of Y(lambda)(I-Y(lambda))^{-1}}.
\hfill$\blacksquare$
\vskip1pc\noindent
\par

Now we are in the position to show Proposition \ref{integral kernel of M^{(1)}}. 
\par
\noindent
Proof of Proposition \ref{integral kernel of M^{(1)}}:
We put 
$$
A^{ijp}(\xi, \zeta; \lambda) = \int_{\partial{D_p}}
{}^tY_{22}^{ip}(\xi, \eta; \lambda)
M^{pj}(\eta, \zeta; \lambda)dS_\eta
\quad(\xi \in \partial{D_i}, \zeta \in \partial{D_j}), 
$$
which are in the integral representation 
(\ref{representation of M^{(1)}_{ij}(x, y, lambda)})
of the integral kernel $ M^{(1)}_{ij}(\xi, \zeta; \lambda) $ of 
$M^{(1)}(\lambda)$.
Here we consider the following four cases (i)-(iv), though 
they do not correspond to the partition of the possible cases.

\par\noindent
(i) The case $j \neq p$: 
From the first estimate in Lemma \ref{integral kernel of Y(lambda)(I-Y(lambda))^{-1}}
and (\ref{estimate of the kernel of Y(lambda)}), there exists a constant $C > 0$ such that
for any $0 < \delta \leq 1$ and $\lambda \in \C_{\delta_0}, \mu \geq \mu_1\delta^{-3}$,
\begin{align*}
\vert A^{ijp}(\xi, \zeta; \lambda) \vert &\leq C
\delta^{-3}e^{-{\delta}d_1\mu}\int_{\partial{D_p}}
\left(\mu+\frac{1}{\vert \xi - \eta \vert}\right)
e^{-\mu\vert \xi - \eta \vert}e^{-(1-\delta)\mu\vert \eta - \zeta \vert}
dS_\eta
\\&
\leq C\delta^{-3}e^{-{\delta}d_1\mu}e^{-(1-\delta)\mu\vert \xi - \zeta \vert}
\int_{\partial{D_p}}\left(\mu+\frac{1}{\vert \xi - \eta \vert}\right)
e^{-\mu\delta\vert \xi - \eta \vert}dS_\eta.
\end{align*}
Hence, Lemma \ref{estimates frequently used 2} implies 
\begin{align}
\vert A^{ijp}(\xi, \zeta; \lambda) \vert &\leq 
C\delta^{-4}e^{-{\delta}d_1\mu}e^{-(1-\delta)\mu\vert \xi - \zeta \vert}
\label{for proof of integral kernel of Y(lambda)(I-Y(lambda))^{-1} no 1}
\\
\quad(\xi \in \partial{D_i}, \zeta \in \partial{D_j}, &i, j, p = 1, 2, \ldots, N, 
j \neq p, \lambda \in \C_{\delta_0}, \mu \geq \mu_1\delta^{-3}, 0 < \delta \leq 1).
\nonumber 
\end{align}

\vskip1pc\noindent
(ii) The case $j \neq p$ and $i = j$:
For these $i$ and $j$, as in the case (i), it follows that
\begin{align*}
\vert A^{jjp}(\xi, \zeta; \lambda) \vert &\leq 
C\delta^{-3}e^{-{\delta}d_1\mu}
\int_{\partial{D_p}}
\left(\mu+\frac{1}{\vert \xi - \eta \vert}\right)
e^{-(1 - \delta)\mu(\vert \eta - \zeta \vert+ \vert \xi - \eta \vert)}
e^{-\delta\mu\vert \xi - \eta \vert}dS_\eta
\\&
\hskip50mm
(\xi, \zeta \in \partial{D}_j, \lambda \in \C_{\delta_0}, 
\mu \geq \mu_1\delta^{-3}, 0 < \delta \leq 1).
\end{align*}
Hence, Lemma \ref{L1 for diagonal part of the repeated kernels for off-diagonal part}
and Lemma \ref{estimates frequently used 2} yield 
\begin{align}
\vert A^{jjp}(\xi, \zeta; \lambda) \vert &\leq 
C\delta^{-3}e^{-{\delta}d_1\mu}
e^{-\mu\vert \xi - \zeta \vert}e^{-2d_3\mu}
\int_{\partial{D_p}}
\left(\mu+\frac{1}{\vert \xi - \eta \vert}\right)
e^{-\delta\mu\vert \xi - \eta \vert}dS_\eta
\nonumber\\
&\leq
C\delta^{-4}e^{-{\delta}d_1\mu}e^{-\mu\vert \xi - \zeta \vert}e^{-2d_3\mu}
\label{for proof of integral kernel of Y(lambda)(I-Y(lambda))^{-1} no 1'}
\\
\quad(\xi, \zeta \in \partial{D_j}, &j, p = 1, 2, \ldots, N, 
j \neq p, 
\lambda \in \C_{\delta_0}, \mu \geq \mu_1\delta^{-3}, 0 < \delta \leq \delta_1).
\nonumber 
\end{align}

\vskip1pc\noindent
(iii) The case $j = p$:
If this is the case, the second estimate in 
Lemma \ref{integral kernel of Y(lambda)(I-Y(lambda))^{-1}} and 
(\ref{estimate of the kernel of Y(lambda)}) implies that there 
exist constants $C > 0$, $\mu_1 > 0$ and $0 < \delta_1 \leq 1$ 
such that for any $0 < \delta \leq \delta_1$ and  
$\lambda \in \C_{\delta_0}, \mu \geq \mu_1\delta^{-3}$, 
\begin{align*}
\Big\vert A^{ijj}(\xi, \zeta; \lambda) &- \int_{\partial{D_j}}
{}^tY_{22}^{ij}(\xi, \eta; \lambda)M_{D_j}(\eta, \zeta; \lambda)dS_\eta \Big\vert 
\\&
= 
\Big\vert \int_{\partial{D_j}}{}^tY_{22}^{ij}(\xi, \eta; \lambda)
\big(M_{jj}(\eta, \zeta; \lambda) - M_{D_j}(\eta, \zeta; \lambda)\big)dS_\eta \Big\vert 
\\&\leq C
\int_{\partial{D_j}}
\left(\mu+\frac{1}{\vert \xi - \eta \vert}\right)
e^{-\mu\vert \xi - \eta \vert}
\delta^{-2}e^{-{\delta}d_1\mu}e^{-\mu{d_3}}e^{-\mu\vert \eta - \zeta \vert}
dS_\eta
\\&
\leq C\delta^{-2}e^{-{\delta}d_1\mu}e^{-\mu\vert \xi - \zeta \vert}e^{-\mu{d_3}}
\int_{\partial{D_j}}\left(\mu+\frac{1}{\vert \xi - \eta \vert}\right)dS_\eta.
\end{align*}
From (\ref{estimate of integral of 1/v{x - z} on partial{D}_p})
and 
${\mu}e^{-\mu{d_3}} \leq {d_3}^{-1}$, there exist constants
$C > 0$ and $0 < \delta_1 \leq 1$ such that
\begin{align}
\vert A^{ijj}(\xi, \zeta; \lambda) - \int_{\partial{D_j}}&
{}^tY_{22}^{ij}(\xi, \eta; \lambda)M_{D_j}(\eta, \zeta; \lambda)dS_\eta \vert \leq 
C\delta^{-2}e^{-{\delta}d_1\mu}e^{-\mu\vert \xi - \zeta \vert}
\label{for proof of integral kernel of Y(lambda)(I-Y(lambda))^{-1} no 2}
\\&
\quad(\xi \in \partial{D_i}, \zeta \in \partial{D_j}, 
\lambda \in \C_{\delta_0}, \mu \geq \mu_1\delta^{-3}, 0 < \delta \leq \delta_1).
\nonumber 
\end{align}

\par\noindent
(iv) The case $j = p$ and $i \neq j$:
If this is the case, (\ref{form of tilde{W}_{ij}(x, y, lambda)}) and 
(\ref{for proof of integral kernel of Y(lambda)(I-Y(lambda))^{-1} no 2}) yield 
$$
\vert A^{ijj}(\xi, \zeta; \lambda) - \tilde{W}^{ij}(\xi, \zeta; \lambda)\vert
\leq C\delta^{-2}e^{-\delta{d_1}\mu}e^{-\mu\v{\xi - \zeta}}.
$$
Hence, estimate (\ref{estimate of tilde{W}_{ij}(x, y, lambda) (i neq j)}) 
for $\tilde{W}^{ij}(\xi, \zeta; \lambda)$ $(i \neq j)$ implies that
\begin{align}
\vert A^{ijj}(\xi, \zeta; \lambda) \vert &\leq 
C\delta^{-2}e^{-{\delta}d_1\mu}e^{-(1-\delta)\mu\vert \xi - \zeta \vert}
\label{for proof of integral kernel of Y(lambda)(I-Y(lambda))^{-1} no 3}
\\&
\quad(\xi \in \partial{D_i}, \zeta \in \partial{D_j}, i \neq j, 
\lambda \in \C_{\delta_0}, \mu \geq \mu_1\delta^{-3}, 0 < \delta \leq \delta_1).
\nonumber 
\end{align}

\par

Now we are in the position to give the estimate of $M^{(1)}_{ij}(\xi, \zeta; \lambda)$ 
for the case $i \neq j$. 
For (\ref{representation of M^{(1)}_{ij}(x, y, lambda)}), it follows that
$$
M^{(1), ij}(\xi, \zeta; \lambda) =\frac{1}{2\pi}e^{-\lambda\vert \xi - \zeta \vert}
H_1(\xi, \zeta)
+ \sum_{p = 1}^NA_{ijp}(\xi, \zeta; \lambda).
$$
From (\ref{for proof of integral kernel of Y(lambda)(I-Y(lambda))^{-1} no 1}), 
(\ref{for proof of integral kernel of Y(lambda)(I-Y(lambda))^{-1} no 3}) and 
the argument for getting 
(\ref{estimate of Y_{ij}(x, y, lambda) for i neq j}) implies 
\begin{align*}
\vert M^{(1)}_{ij}(\xi, \zeta; \lambda) \vert
&\leq 
C\delta^{-4}e^{-{\delta}d_1\mu}e^{-(1-\delta)\mu\vert \xi - \zeta \vert}
\\&
\quad(\xi \in \partial{D_i}, \zeta \in \partial{D_j}, i \neq j, 
\lambda \in \C_{\delta_0}, \mu \geq \mu_1\delta^{-3}, 0 < \delta \leq \delta_1).
\end{align*}

\par
Next is the case $i = j$. The representations 
(\ref{representation of M^{(1)}_{ij}(x, y, lambda)}) and 
(\ref{representation of M^{(1)}_{D_j}(x, y, lambda)})
yield
\begin{align*}
M^{(1), jj}(\xi, \zeta; \lambda) &- M^{(1)}_{D_j}(\xi, \zeta; \lambda)
\\&
= \sum_{p \neq j}A^{jjp}(\xi, \zeta; \lambda)
+ A^{jjj}(\xi, \zeta; \lambda) - \int_{\partial{D_j}}
{}^tY_{22}^{ij}(\xi, \eta; \lambda)M_{D_j}(\eta, \zeta; \lambda)dS_\eta.
\end{align*}
This equality, (\ref{for proof of integral kernel of Y(lambda)(I-Y(lambda))^{-1} no 1'}) and
(\ref{for proof of integral kernel of Y(lambda)(I-Y(lambda))^{-1} no 2}) imply that
\begin{align*}
\vert M^{(1), jj}(\xi, \zeta; \lambda) - M^{(1)}_{D_j}(\xi, \zeta; \lambda) \vert
&\leq 
C\delta^{-4}e^{-{\delta}d_1\mu}e^{-\mu\vert \xi - \zeta \vert}
\\&
\quad(\xi, \zeta \in \partial{D_j}, 
\lambda \in \C_{\delta_0}, \mu \geq \mu_1\delta^{-3}, 0 < \delta \leq \delta_1),
\end{align*}
which completes the proof of Proposition \ref{integral kernel of M^{(1)}}.
\hfill$\blacksquare$
\vskip1pc\noindent
\par
\setcounter{equation}{0}
\section{Estimate of some Laplace integrals}
\label{Estimate of some Laplace integrals}%

Let $U \subset \R^n$ be a bounded open set, and
$S(\sigma)$ is a $C^2$ function in $U$, and
$h(\sigma; \lambda)$ be a continuous function 
in $\sigma \in U$ with a parameter $\lambda \in \C_{\delta_0}$ 
for some $\delta_0 > 0$. For $S$ and $h$, assume that 
\par\noindent
(S.1) $\tau_{-\infty} = \inf_{\sigma \in U}S(\sigma)$ exists and $\tau_{-\infty} = S(0)$, 
\par\noindent
(S.2) there exists a constant $C_0 > 0$ such that
$\tau_{-\infty} \leq S(\sigma) \leq \tau_{-\infty}+C_0\v{\sigma}^2$ 
$(\sigma \in U)$,
\par\noindent
(H.1) $h$ is of the form: 
$h(\sigma; \lambda) = h_1(\sigma; \lambda)+\lambda^{-1}\tilde{h}_1(\sigma; \lambda)$
$(\sigma \in U, \lambda \in \C_{\delta_0})$,
\par\noindent
(H.2) there exist constants $C_1 > 0$, $C_1' > 0$ and $\mu_0 > 0$ such that
$$
{\rm Re}\,[h_1(\sigma; \lambda)] \geq C_1, \qquad
\v{h_1(\sigma; \lambda)}+\v{\tilde{h}_1(\sigma; \lambda)} \leq C_1'
\qquad(\sigma \in U, \lambda \in \C_{\delta_0}, \mu \geq \mu_0).
$$
For the functions $S$ and $h$, and a cutoff function $\varphi \in C_0^2(U)$ with
$0 \leq \varphi \leq 1$ and $\varphi(0) = 1$, we introduce a Laplace integral
$I(\lambda)$ of the form:
\begin{align}
I(\lambda) = \int_{U}e^{-{\lambda}S(\sigma)}\varphi(\sigma)h(\sigma; \lambda)d\sigma.
\label{Laplace integral}
\end{align}

\begin{Prop}\label{lower estimate for some Laplace integral}
For integral (\ref{Laplace integral}), assume that $S$ and $h$ satisfy
(S.1), (S.2), (H.1) and (H.2) in the above. Then there exist constants 
$0 < \delta_1 < \delta_0$, $\mu_1 > 0$ and $C > 0$ such that
$$
{\rm Re}\,[e^{\lambda\tau_{-\infty}}I(\lambda)] \geq C\mu^{-n/2}
\quad(\lambda \in \Lambda_{\delta_1}, \mu \geq \mu_1).
$$
\end{Prop}
Proof: We put $\tau_{\infty} = \sup_{\sigma \in U}S(\sigma)$, 
$E_\tau = \{\sigma \in U \,\vert\, S(\sigma) \leq \tau \}$ for $\tau \in \R$, and
$$
\beta_{\lambda}(\tau) = \int_{E_\tau}\varphi(\sigma)h(\sigma; \lambda)d\,\sigma
\qquad(\tau \in \R).
$$
Note that $\beta_{\lambda}(\tau)$ is a function of bounded variation, 
$\beta_{\lambda}(\tau) = 0$ for $\tau < \tau_{-\infty}$ and
$\beta_{\lambda}(\tau) = \beta_{\lambda}(\tau_{\infty})$ 
for $\tau \geq \tau_{\infty}$.
Note also that $\beta_{\lambda}$ is a right continuous function in $\tau \in \R$ 
since for any $\tau_0 \in \R$, $\lim_{\tau \to \tau_0+0}\chi_{E_\tau}(\sigma) = \chi_{E_{\tau_0}}(\sigma)$,
where $\chi_{E_\tau}(\sigma)$ is the characteristic function of the set $E_\tau$.
From Stieltjes integral with respect to $\beta_{\lambda}$, 
for any $\tilde\tau_{-\infty} < \tau_{-\infty}$, it follows that
\begin{align}
I(\lambda) &= \int_{\tilde\tau_{-\infty}}^{\tau_{\infty}}
e^{-\lambda\tau}d\beta_{\lambda}(\tau)
= 
e^{-\lambda\tau_{\infty}}\beta_{\lambda}(\tau_{\infty})
+\lambda\int_{\tau_{-\infty}}^{\tau_{\infty}}
e^{-\tau\lambda}\beta_{\lambda}(\tau)d\tau.
\label{for estimate from below 1}
\end{align}

We put
$$
\beta_{\lambda, 0}(\tau) = \int_{E_\tau}\varphi(\sigma)
h_1(\sigma; \lambda)d\,\sigma
\qquad(\tau \in \R).
$$
From (H.2), it follows that
\begin{align}
\v{\beta_{\lambda, 0}(\tau)-\beta_{\lambda, 0}(\tau_{-\infty})} \leq
C_1'\int_{(E_\tau\setminus{E_{\tau_{-\infty}})}\cup({E_{\tau_{-\infty}}\setminus{E_\tau}})}
&\varphi(\sigma)d\sigma
\leq C_1'\V{\varphi}_{L^1(U)}
\label{for estimate from below 2}
\\&
\quad(\tau \in \R, \lambda \in \C_{\delta_0}, \mu \geq \mu_0). 
\nonumber
\end{align}
Note also that (H.1) and (H.2) yield that
\begin{align}
\v{\beta_{\lambda}(\tau)-\beta_{\lambda, 0}(\tau)} &\leq
C_1'C_1^{-1}\v{\lambda}^{-1}\int_{E_\tau}\varphi(\sigma)
{\rm Re}\,[h_1(\sigma; \lambda)]d\sigma 
\nonumber\\&
\leq C_2\v{\lambda}^{-1}{\rm Re}\beta_{\lambda, 0}(\tau)
\leq C_2C_1'\v{\lambda}^{-1}\V{\varphi}_{L^1(U)}
\label{for estimate from below 3}
\\&
\hskip30mm
\quad(\tau \in \R, \lambda \in \C_{\delta_0}, \mu \geq \mu_0),
\nonumber
\end{align}
where  $C_2 = C_1'/C_1 > 0$. 
A similar argument for getting 
(\ref{for estimate from below 3}) implies
\begin{align}
\v{\beta_{\lambda}(\tau)} \leq C_1'\V{\varphi}_{L^1(U)},
\quad
\v{{\rm Im}\beta_{\lambda, 0}(\tau)} \leq C_2{\rm Re}\beta_{\lambda, 0}(\tau)
\quad(\tau \in \R, \lambda \in \C_{\delta_0}, \mu \geq \mu_0),
\label{for estimates of beta and beta_0}
\end{align}
which yields 
\begin{align}
{\rm Re}\, \Big(e^{\lambda\tau_{-\infty}}\lambda
&\int_{\tau_{-\infty}}^{\tau_{-\infty}}e^{-\tau\lambda}
\beta_\lambda(\tau)d\tau\Big)
\geq J_\delta(\lambda) 
-C_1'\V{\varphi}_{L^1(U)}\frac{\v{\lambda}}{\mu}e^{-\mu\delta}
\label{estimate of 2-nd term}
\end{align}
for any $0 \leq \delta \leq \tau_0$, 
where $\tau_0 = \tau_{\infty} - \tau_{-\infty}$ and
$$
J_\delta(\lambda) = 
{\rm Re}\, \Big(e^{\lambda\tau_{-\infty}}
\lambda\int_{\tau_{-\infty}}^{\tau_{-\infty}+\delta}e^{-\tau\lambda}
\beta_{\lambda, 0}(\tau)d\tau \Big).
$$

Further, (H.2) implies an estimate of ${\rm Re}\beta_{\lambda, 0}$ from below:
\begin{align}
{\rm Re}\beta_{\lambda, 0}(\tau) \geq C_1\gamma(\tau)
\quad\quad(\tau \in \R, \lambda \in \C_{\delta_0}, \mu \geq \mu_0),
\label{for estimate from below 4}
\end{align}
where 
$$
\gamma(\tau) = \int_{E_\tau}\varphi(\sigma)d\sigma
\qquad(\tau \in \R).
$$

\par

We can divide the following three cases:
Case 1: $\tau_{-\infty} = \tau_{\infty}$, 
Case 2: $\tau_{-\infty} < \tau_{\infty}$ 
and $\gamma(\tau_{-\infty}) > 0$, 
%
Case 3: $\tau_{-\infty} < \tau_{\infty}$ and 
$\gamma(\tau_{-\infty}) = 0$.
\par
Case 1: In this case, $E_{\tau_{-\infty}} = U$. This and assumption (H.2)
imply 
$$
{\rm Re}\,\beta_{\lambda, 0}(\tau_{-\infty}) \geq C_1\int_{U}\varphi(\sigma)d\sigma
= C_1\V{\varphi}_{L^1(U)}
\qquad(\lambda \in \C_{\delta_0}, \mu \geq \mu_0).
$$
This estimate, (\ref{for estimate from below 1}) and 
(\ref{for estimate from below 3}) yield 
$$
{\rm Re}\,[e^{\lambda\tau_{-\infty}}I(\lambda)] 
= {\rm Re}\,\beta_{\lambda}(\tau_{-\infty})
\geq (1-C_2\v{\lambda}^{-1}){\rm Re}\,\beta_{\lambda, 0}(\tau_{-\infty})
\geq C_1(1-C_2\v{\lambda}^{-1})\V{\varphi}_{L^1(U)}, 
$$
which implies 
${\rm Re}\,[e^{\lambda\tau_{-\infty}}I(\lambda)] \geq C$ for some 
constant $C > 0$ for large $\v{\lambda}$ in $\lambda \in \C_{\delta_0}$
uniformly.

\par
Case 2: In this case,  we put $C_3 = C_1\gamma(\tau_{-\infty}) > 0$.
Since $E_\tau \supset E_{\tau_{-\infty}}$ for $\tau \geq \tau_{-\infty}$
and $\lim_{\tau \to \tau_{-\infty}+0}
\chi_{E_\tau\setminus{E_{\tau_{-\infty}}}} = 0$, 
(\ref{for estimate from below 2}) implies that 
there exists a constant $\delta > 0$ such that
$$
\v{\beta_{\lambda, 0}(\tau)-\beta_{\lambda, 0}(\tau_{-\infty})} 
\leq \frac{C_3}{2(1+\delta_0)}
\quad(\tau_{-\infty} \leq \tau \leq \tau_{-\infty}+\delta).
$$
From the above estimate, (\ref{for estimate from below 3}) and (\ref{for estimate from below 4}), 
it follows that
\begin{align*}
J_\delta(\lambda) &\geq 
{\rm Re}\, \Big(e^{\lambda\tau_{-\infty}}
\lambda\int_{\tau_{-\infty}}^{\tau_{-\infty}+\delta}
e^{-\tau\lambda}\beta_{\lambda, 0}(\tau_{-\infty})d\tau \Big)
\\&\hskip20mm
-\v{\lambda}\int_{\tau_{-\infty}}^{\tau_{-\infty}+\delta}
e^{-(\tau-\tau_{-\infty})\mu}
\big(\v{\beta_\lambda(\tau)-\beta_{\lambda, 0}(\tau)}
+ \v{\beta_{\lambda, 0}(\tau)-\beta_{\lambda, 0}(\tau_{-\infty})}
\big)d\tau
\\&
\geq {\rm Re}\,\big[\beta_{\lambda, 0}(\tau_{-\infty})\big(1
-e^{-\lambda\delta}\big)\big]
-\frac{C_1'C_2}{\mu}\V{\varphi}_{L^1(U)}
- \frac{C_3\v{\lambda}}{2(1+\delta_0)\mu}
\\&
\geq \frac{C_3}{2}-C\big(\mu^{-1}+e^{-\mu\delta}\big)
\qquad(\lambda \in \C_{\delta_0}).
\end{align*}
Combining the above estimate and (\ref{estimate of 2-nd term}) with
(\ref{for estimate from below 1}), 
we obtain ${\rm Re}\,[e^{\lambda\tau_{-\infty}}I(\lambda)] \geq C$ 
for some constant $C > 0$ for large $\v{\lambda}$ in 
$\lambda \in \C_{\delta_0}$ uniformly.

\par
Case 3: In this case, take $r_1 > 0$ with 
$\overline{B(0, r_1)} \subset U$ and
$\varphi(\sigma) \geq \varphi(0)/2$ $(\v{\sigma} \leq r_1)$, 
where $B(0, r_1)$ is the open ball with the center $0$ and
the radius $r_1$. Note that 
(S.2) implies that
$\overline{B(0, \sqrt{(\tau-\tau_{-\infty})/C_0})} \subset E_\tau$
for $\tau_{-\infty} \leq \tau \leq \tau_{-\infty}+C_0r_1^2$, 
which yields 
\begin{align*}
\gamma(\tau) 
\geq \frac{\varphi(0)}{2C_0^{n/2}}{\rm Vol}(B(0, 1))
(\tau-\tau_{-\infty})^{n/2}
\quad\quad(\tau_{-\infty} \leq \tau \leq \tau_{-\infty}+C_0r_1^2).
\end{align*}
Hence, taking 
$C_4 = 2^{-1}\varphi(0){\rm Vol}(B(0, 1))\min\{C_0^{-n/2},
\frac{r_1^n}{(\tau_{\infty}-\tau_{-\infty})^{n/2}}\}$,
we obtain
\begin{align}
\gamma(\tau) 
\geq C_4(\tau-\tau_{-\infty})^{n/2}
\quad\quad(\tau_{-\infty} \leq \tau \leq \tau_{\infty}).
\label{for estimate from below 5}
\end{align}

\par
From (\ref{for estimates of beta and beta_0}), for any $0 < \delta \leq \tau_0$, 
it follows that
\begin{align*}
J_\delta(\lambda) &\geq 
{\rm Re}\, \Big(\lambda
\int_{0}^{\delta}e^{-\tau\lambda}
\beta_{\lambda, 0}(\tau+\tau_{-\infty})d\tau\Big)
-C_2\int_{0}^{\delta}\big\vert{{\rm Im }\big({\lambda}e^{-\tau\lambda}\big)}\big\vert
{\rm Re}\beta_{\lambda, 0}(\tau+\tau_{-\infty})d\tau. 
\end{align*}
From this estimate, (\ref{estimate of 2-nd term}) 
and (\ref{for estimate from below 1}), 
for any $0 < \delta \leq \tau_0$ and $\lambda \in \Lambda_{\delta_1}$, we obtain
\begin{align*}
{\rm Re}\,[e^{\tau_{-\infty}\lambda}I(\lambda)]
\geq \mu\int_{0}^{\delta}e^{-\tau\mu}\Phi(\tau, \lambda)
{\rm Re}\beta_{\lambda, 0}(\tau+\tau_{-\infty})d\tau
-C_1'\V{\varphi}_{L^1(U)}(e^{-\tau_{0}\mu}+\frac{\v{\lambda}}{\mu}e^{-\delta\mu}), 
\end{align*}
where 
$$
\Phi(\tau, \lambda) = \cos({\rm Im}\,\lambda\tau)-
C_2\v{\sin({\rm Im}\,\lambda\tau)}
+\frac{{\rm Im}\,\lambda}{\mu}\sin({\rm Im}\,\lambda\tau)
-\frac{C_2\v{{\rm Im}\,\lambda}}{\mu}
\v{\cos({\rm Im}\,\lambda\tau)}.
$$
\par
We take constants $0 < c_0 < 1$ and $0 < \theta_0 < \pi/2$ satisfying
$\cos{x} - C_2\v{\sin x} \geq 2c_0$ for $\v{x} \leq \theta_0$, 
and choose $\delta = \min\{\theta_0/\v{{\rm Im}\,\lambda}, \tau_0\}$ and
$\mu_1 = e^{\delta_1(C_2+1)/c_0}$.
Since  
\begin{align*}
\Phi(\tau, \lambda) \geq \cos({\rm Im}\,\lambda\tau)
-C_2\v{\sin({\rm Im}\,\lambda\tau)}
-\frac{\delta_1(C_2+1)}{\log\mu}
\qquad(\lambda \in \Lambda_{\delta_1}), 
\end{align*}
it follows that $\Phi(\tau, \lambda) \geq c_0$ for $ \lambda \in \Lambda_{\delta_1}$,
$\mu \geq \mu_1$ and $0 \leq \tau \leq \delta$.
From this fact, (\ref{for estimate from below 4}) and
(\ref{for estimate from below 5}), and $\v{\lambda}/\mu \leq 1+\delta_1(\log\mu)^{-1} \leq 2$
for $\lambda \in \Lambda_{\delta_1}$ with $\mu \geq \mu_1$ and $0 < \delta_1 \leq 1$, 
it follows that
\begin{align*}
{\rm Re}\,[e^{\tau_{-\infty}\lambda}I(\lambda)]
\geq c_0C_1C_4\mu^{-n/2}\int_{0}^{\delta\mu}e^{-\tau}\tau^{n/2}d\tau
-3C_1'\V{\varphi}_{L^1(U)}e^{-\delta\mu}.
\end{align*}
If $\v{{\rm Im}\,\lambda} \leq \theta_0/\tau_0$, 
$\delta\mu = \tau_0\mu \geq \tau_0e$, and if 
$\v{{\rm Im}\,\lambda} \geq \theta_0/\tau_0$,
$\delta\mu = \theta_0\mu/\v{{\rm Im}\,\lambda}
\geq \theta_0\delta_1^{-1}\log\mu \geq \theta_0$ ($0 < \delta_1 \leq 1$
and $\lambda \in \Lambda_{\delta_1}$). 
Hence, in any case, we obtain
\begin{align*}
{\rm Re}\,[e^{\tau_{-\infty}\lambda}I(\lambda)]
\geq c_0C_1C_4\mu^{-n/2}\int_{0}^{\min\{\tau_0e, \theta_0\}}
e^{-\tau}\tau^{n/2}d\tau
&-3C_1'\V{\varphi}_{L^1(U)}\big(\mu^{-\theta_0\delta_1^{-1}}+e^{-\tau_0\mu}\big)
\\&
\quad(\lambda \in \Lambda_{\delta_1}, \mu \geq \mu_1).
\end{align*}
This implies 
${\rm Re}\,[e^{\lambda\tau_{-\infty}}I(\lambda)] \geq C\mu^{-n/2}$ for some 
constant $C > 0$ for large $\v{\lambda}$ in $\lambda \in \Lambda_{\delta_1}$
uniformly if we take $\delta_1$ sufficiently small.
This completes the proof of 
Proposition \ref{lower estimate for some Laplace integral}.
\par
\hfill$\blacksquare$
\vskip1pc\noindent
\par

\par

Next, we treat the non-degenerate case, i.e. 
\par\noindent
(S.3) $\nabla_{\sigma}S(0) = 0$, ${\rm Hess }S(0) > 0$ and 
$S(\sigma) > \tau_{-\infty}$ ($0 \neq \sigma \in U$) 
\par\noindent
is assumed. For the amplitude function $h(\sigma; \lambda)$, we also assume 
\par\noindent
(H.3) there exists a constant $\mu_0 > 0$ such that 
$\lim_{\sigma \to 0}h(\sigma; \lambda) = h(0; \lambda)$ 
uniformly in $\lambda \in \C_{\delta_0}$ with $\mu \geq \mu_0$,
\par\noindent
(H.4) $h(\sigma; \lambda)$ is bounded for $\sigma \in U$ and 
$\lambda \in \C_{\delta_0}$.

\begin{Prop}\label{lower estimate for some Laplace integral for non-degenerate case}
Assume that $S(\sigma)$ satisfies (S.1) and (S.3). If $h(\sigma; \lambda)$ 
($\sigma \in U$, $\lambda \in \C_{\delta_0}$) is continuous in $\sigma \in U$, 
then
there exists a constant $C > 0$ such that
$$
\v{I(\lambda)} \leq Ce^{-\mu\tau_{-\infty}}\mu^{-n/2}\V{\varphi(\cdot)h(\cdot\,; \lambda)}_{C(U)}
\quad(\lambda \in \C_{\delta_0}),
$$
where $I(\lambda)$ is given by (\ref{Laplace integral}). Further, assume also that
$h$ satisfies (H.3) and (H.4). Then the following asymptotic formula holds:
$$
I(\lambda) = \frac{e^{-\lambda\tau_{-\infty}}}{\sqrt{{\rm Hess }S(0)}}\lambda^{-n/2}
\Big(h(0; \lambda)+o(1)\Big)
\qquad(\text{as } \mu \to \infty \text{ uniformly in } \lambda \in \C_{\delta_0}).
$$
\end{Prop}
Proof: Take a cutoff function $\psi \in C^\infty_0(U)$ satisfies $\psi = 1$ near ${\rm supp }\varphi$ and
$0 \leq \psi \leq 1$, and decompose the integral $I(\lambda)$ in (\ref{Laplace integral}) 
as follows:
\begin{align}
I(\lambda) = \varphi(0)h(0; \lambda)\int_{U}e^{-{\lambda}S(\sigma)}\psi(\sigma)d\sigma
+ \int_{U}e^{-{\lambda}S(\sigma)}\psi(\sigma)\tilde{h}(\sigma; \lambda)d\sigma, 
\label{decomposition for non-degenerate case}
\end{align}
where
$\tilde{h}(\sigma; \lambda) = \varphi(\sigma)h(\sigma; \lambda) - \varphi(0)h(0; \lambda)$. 
We write the first and second terms of the right side of 
(\ref{decomposition for non-degenerate case}) as 
$I_1(\lambda)$ and $I_2(\lambda)$ respectively. 
From the usual Laplace method, $I_1(\lambda)$ is expanded as
\begin{align}
I_1(\lambda) = h(0; \lambda)\frac{e^{-\lambda\tau_{-\infty}}}{\sqrt{{\rm Hess }S(0)}}\lambda^{-n/2}
\Big(1+O(\lambda^{-1})\Big)
\quad(\lambda \in \C_{\delta_0}, {\rm Re }\lambda \to \infty).
\label{for the main part}
\end{align}

\par
From (S.3), it follows that there exists a constant $C_0' > 0$ such that
$S(\sigma) \geq \tau_{-\infty}+C_0'\v{\sigma}^2$ $(\sigma \in U)$,
which yields
\begin{align}
\v{I_2(\lambda)} \leq 
\mu^{-n/2}e^{-\mu\tau_{-\infty}}\int_{\R^n}e^{-C_0'\v{\sigma}^2}
\v{\psi(\mu^{-1/2}\sigma)\tilde{h}(\mu^{-1/2}\sigma; \lambda)}d\sigma. 
\label{estimate of I_2}
\end{align}
Put $M = \sup_{\sigma \in U, \lambda \in \C_{\delta_0}}\v{h(\sigma; \lambda)} < \infty$
for (H.4). There exists a constant $C > 0$ such that 
\begin{align}
\v{\tilde{h}(\sigma; \lambda)} 
&\leq \v{\varphi(\sigma)}\v{h(\sigma; \lambda) - h(0; \lambda)}
+\v{h(0; \lambda)}\v{\varphi(\sigma) - \varphi(0)}
\nonumber
\\& 
\leq C\{\v{h(\sigma; \lambda) - h(0; \lambda)}+ M\v{\sigma}\}
\quad(\sigma \in U, \lambda \in \C_{\delta_0}).
\label{estimate of tilde{h}}
\end{align}
For any $\eta_0 > 0$, it follows that 
\begin{align*}
\v{\psi(\mu^{-1/2}\sigma)(h(\mu^{-1/2}\sigma; \lambda) - h(0; \lambda))}
&\leq \sup_{\v{\sigma} \leq \eta_0\mu^{-1/2}}
\v{h(\sigma; \lambda) - h(0; \lambda)}
\,\, (\v{\sigma} \leq \eta_0, \lambda \in \C_{\delta_0})
\intertext{and}
\v{\psi(\mu^{-1/2}\sigma)(h(\mu^{-1/2}\sigma; \lambda) - h(0; \lambda))}
&\leq 2M
\,\, (\v{\sigma} \geq \eta_0, \lambda \in \C_{\delta_0}). 
\end{align*}
These estimates and (\ref{estimate of I_2}) imply that there exists a constant $C > 0$ independent 
of $\eta_0 > 0$ such that
\begin{align*}
\v{I_2(\lambda)} \leq C\mu^{-n/2}e^{-\mu\tau_{-\infty}}
&\Big(\mu^{-1/2}M
+\sup_{\v{\sigma} \leq \eta_0\mu^{-1/2}}\v{h(\sigma; \lambda) - h(0; \lambda)}
\\
&\hskip20mm
+ M\int_{\v{\sigma} \geq \eta_0}e^{-C_0'\v{\sigma}^2}d\sigma\Big)
\quad(\lambda \in \C_{\delta_0}).
\nonumber
\end{align*}
Hence, taking $\eta_0 = \mu^{1/4}$ in the above estimate, and noting (H.3), 
(\ref{for the main part}) and (\ref{decomposition for non-degenerate case}),
we obtain
the asymptotic behavior of $I(\lambda)$ in 
Proposition \ref{lower estimate for some Laplace integral for non-degenerate case}. 
\par

Similarly to (\ref{estimate of I_2}), we have
\begin{align*}
\v{I(\lambda)} &\leq 
\mu^{-n/2}e^{-\mu\tau_{-\infty}}\int_{\R^n}e^{-C_0'\v{\sigma}^2}
\v{\varphi(\mu^{-1/2}\sigma)h(\mu^{-1/2}\sigma; \lambda)}d\sigma
\\& 
\leq \mu^{-n/2}e^{-\mu\tau_{-\infty}}\V{\varphi(\cdot)h(\cdot\,; \lambda)}_{C(U)}
\int_{\R^n}e^{-C_0'\v{\sigma}^2}d\sigma, 
\end{align*}
which shows the estimate of $I(\lambda)$ in 
Proposition \ref{lower estimate for some Laplace integral for non-degenerate case}. 
\hfill$\blacksquare$
\vskip1pc\noindent
\par

\begin{Remark}\label{the case of Horder continuous boundary}
Instead of (H.3) and (H.4), assume that $h(\cdot; \lambda)$ is H\"order continuous in $\sigma \in U$
of order $0 < \alpha_0 < 1 $. 
In this case, from (\ref{estimate of tilde{h}}), it follows that
$$
\v{\psi(\mu^{-1/2}\sigma)\tilde{h}(\mu^{-1/2}\sigma; \lambda)}
\leq C\mu^{-\alpha_0/2}\V{h(\cdot; \lambda)}_{C^{0, \alpha_0}(U)}
\quad(\lambda \in \C_{\delta_0}).
$$
Hence, there exist a constant $C > 0$ and a neighborhood $V$ 
of $0$ with $\overline{V} \subset U$ such that 
$\v{I_2(\lambda)} \leq C\mu^{-n/2-\alpha_0/2}e^{-\mu\tau_{-\infty}}
\V{h(\cdot; \lambda)}_{C^{0, \alpha_0}(V)}$ ($\lambda \in \C_{\delta_0}$). 
This estimate, (\ref{for the main part}) and (\ref{decomposition for non-degenerate case})
imply 
\begin{align*}
I(\lambda) = \frac{e^{-\lambda\tau_{-\infty}}}{\sqrt{{\rm Hess }S(0)}}\lambda^{-n/2}
\Big(h(0; \lambda)&+O(\lambda^{-\alpha_0/2})\V{h(\cdot; \lambda)}_{C^{0, \alpha_0}(V)}\Big)
\\& 
\qquad(\text{as } \mu \to \infty \text{ uniformly in } \lambda \in \C_{\delta_0}). 
\end{align*}
\end{Remark}


\setcounter{equation}{0}
\appendix
\renewcommand{\theequation}{A.\arabic{equation}}

\section{The case of one strictly convex cavity with $C^2$ boundary}

We discuss reducing regularities of $\partial{D}$ to obtain
the estimates of $M^{(1)}_{D_j}(\xi, \zeta; \lambda)$ 
in Proposition \ref{s-convex case estimates}. Since this estimate is
for the case of one strictly convex boundary, from now on, we assume that
$\partial{D}$ is a strictly convex $C^2$ surface. As described in Remark \ref{C^2 boundaries are OK},
this estimate is given for $C^{2, \alpha_0}$ boundary with some $\alpha_0 \in (0, 1)$.
In \cite{Ikehata and Kawashita}, for any $\xi \in \partial{D}$, standard local coordinates
\begin{align}
U_\xi \ni \,\sigma=(\sigma_1,\sigma_2)\mapsto 
\xi + \sigma_1 e_1+\sigma_2e_2-g_{\xi}(\sigma_1,\sigma_2)\nu_\xi \in \partial D\cap B(\xi,2r_0)
\label{expression of standard local coordinates}
\end{align}
are used to show the estimate of the integral kernels. In this case, $g_\xi$ can be 
extended as $ g_\xi \in {\mathcal B}^{2, \alpha_0}(\R^2)$ (i.e. $g_\xi \in {\mathcal B}^{2}(\R^2)$
and each derivative $\partial_\sigma^{\alpha}g_{\xi}$ for $\v{\alpha} = 2$ is uniform H\"order continuous
in $\R^2$). Since $g_{\xi}$ is uniformly bounded in ${\mathcal B}^{2, \alpha_0}(\R^2)$ with respect
to $\xi \in \partial{D}$, there exists a constant $C > 0$ such that 
$\v{\partial_\sigma^{\alpha}g_\xi(\sigma') - \partial_\sigma^{\alpha}g_\xi(\sigma)} 
\leq C\v{\sigma' - \sigma}^{\alpha_0}$
for any $\sigma, \sigma' \in \R^2$, $\v{\alpha} = 2$ and $\xi \in \partial{D}$. 
Thus, we can use perturbation arguments. 
When $\partial{D}$ is $C^2$, more delicate arguments than that in \cite{Ikehata and Kawashita}
are necessary since we only have $g_\xi \in {\mathcal B}^2(\R^2)$. 

\par
For $C^2$ class boundary, we need to show the following properties:
\begin{Lemma}\label{equi-continuous of g_xi}
All derivatives $\partial_{\sigma}^{\alpha}g_\xi \in {\mathcal B}(\R^2)$ for $\v{\alpha} \leq 2$ 
of the functions $g_\xi \in {\mathcal B}^2(\R^2)$ for $\xi \in \partial{D}$ given 
in Lemma \ref{standard coordinates} are equi-continuous, that is, for any $\varepsilon > 0$, 
there exists $\delta_\varepsilon > 0$ such
that $\v{\partial_{\sigma}^{\alpha}g_\xi(\tilde\sigma) - \partial_{\sigma}^{\alpha}g_\xi(\sigma)} 
< \varepsilon$ holds for $\v{\tilde\sigma - \sigma} < \delta_\varepsilon$ and $\xi \in \partial{D}$.
\end{Lemma}
A proof of Lemma \ref{equi-continuous of g_xi} is given later. We proceed to show 
how to treat the $C^2$ boundary case.

\par

Take any $\xi \in \partial{D}$ and a standard local coordinate 
(\ref{expression of standard local coordinates}) around $\xi$. Note that
we can choose $r_0 > 0$ in Lemma \ref{standard coordinates} sufficiently small enough.
In what follows, we change $r_0 > 0$ to be small several finite times. 
Since $\partial{D}$ is strictly convex and compact, and (\ref{well known estimate for surface 1})
holds for any $C^2$ surface, there exist constants $M_1 > M_0 > 0$ independent of $r_0$ such that
$M_1\v{\zeta - \xi}^2 \geq -\nu_{\xi}\cdot(\zeta - \xi) \geq M_0\v{\zeta - \xi}^2$ 
($\xi \in \partial{D}$ and $\zeta \in \partial{D}{\cap}B(\xi, 2r_0)$).
Choose $r_0 > 0$ satisfying $M_1r_0 < 1/2 $. For $\sigma \in U_\xi$, we put
$\zeta = \xi + \sigma_1 e_1+\sigma_2e_2-g_{\xi}(\sigma)\nu_\xi \in \partial{D}{\cap}B(\xi, 2r_0))$. 
From 
$$
M_0\v{\sigma}^2 \leq g_\xi(\sigma) \leq M_1\v{\zeta - \xi}^2 \leq \frac{\v{\zeta - \xi}}{2}
\leq \frac{\v{\sigma}+\v{g_\xi(\sigma)}}{2} 
\quad(\sigma \in U_\xi),
$$
$M_0\v{\sigma}^2 \leq g_\xi(\sigma) \leq \v{\sigma}$ holds.
Since $\v{\zeta - \xi}^2 = \v{\sigma}^2+\v{g_\xi(\sigma)}^2 \leq 2\v{\sigma}^2$, we obtain
\begin{equation}
M_0\v{\sigma}^2 \leq g_\xi(\sigma) \leq 2M_1\v{\sigma}^2 
\quad (\sigma \in U_\xi, \xi \in \partial{D}). 
\label{estimate of g_xi}
\end{equation}
We put $r_1 = 2r_0/\sqrt{1+16M_1^2r_0^2} < 2r_0 $. For $\sigma \in U_\xi$, 
it follows that 
$(2r_0)^2 > \v{\zeta - \xi}^2 \geq \v{\sigma}^2$ and 
$\v{\zeta - \xi}^2 = \v{\sigma}^2+\v{g_\xi(\sigma)}^2 \leq \v{\sigma}^2(1+4M_1^2\v{\sigma}^2)$, 
which imply 
\begin{align}
\v{\sigma} \leq \v{\zeta - \xi} &< \sqrt{1+16M_1^2r_0^2}\v{\sigma} \quad
\label{equivalence of length and local coordinate expression}
\\&
(\zeta = \xi+\sigma_1e_1+\sigma_2e_2 - g_\xi(\sigma)\nu_\xi \in \partial{D}{\cap}B(\xi, 2r_0)).
\nonumber 
\end{align}
\par

\par
Take any $\eta \in \partial{D}{\cap}B(x,2r_0)$ with $\xi \neq \eta$ and fixed. 
Choose $\{e_1,\,e_2\}$ in the standard system of local coordinates 
(\ref{expression of standard local coordinates}) around $\xi$ in such a way that $\eta - \xi$ is
perpendicular to $e_2$ and $(\eta - \xi){\cdot}e_1 > 0 $. Thus, 
one can write
$$
\eta = \xi + \sigma_1^0 e_1-g_\xi(\sigma_1^0,0)\nu_\xi
$$
with $(\sigma_1^0)^2+g_\xi(\sigma_1^0,0)^2<(2r_0)^2$ and $\sigma_1^0 > 0 $.

\begin{Prop}\label{estimates of lengths with a relaxed assumption than IK}
Assume that $\partial D$ is of class $C^{2}$ and strictly convex. 
\par\noindent
(i) It follows that 
$$
\v{\xi - \zeta}+\v{\zeta - \eta} \geq \v{\xi - \eta}+\frac{1}{2}
\frac{\sigma_2^2}{\vert \zeta - \xi\vert}
\qquad(\zeta \in \partial D\cap B(\xi, 2r_0)).
$$ 
\par\noindent
(ii)  If $r_0$ is chosen small enough, it follows that 
$$
\v{\xi - \zeta}+\v{\zeta - \eta} \geq \v{\xi - \eta}
+\frac{c_0}{\vert \zeta - \xi \vert}((\sigma_1^0)^2\sigma_1^2+\sigma_2^2)
$$
for all $\sigma=(\sigma_1,\sigma_2)$ and 
$\sigma^0=(\sigma_1^0, 0)$ with
$\sigma_1 < 2\sigma_1^0/3$, $\vert\sigma\vert < r_1$ and $\vert\sigma^0\vert < r_1$,
where $r_1 ={2r_0}/{\sqrt{1+16M_1^2r_0^2}} $ 
and $c_0$ is a positive constant depending only on $\partial D$.
\end{Prop}
Proof: For $\zeta = \xi +\sigma_1e_1+\sigma_2e_2-g(\sigma)\nu_{\xi} 
\in \partial D\cap B(\xi, 2r_0) $,
we put $\rho_0 = \v{\eta - \xi}$, $\rho = \v{\zeta - \xi}$, and denote by $\theta$
the angle made by the line segments $\xi\zeta$ and $\xi\eta$. The cosine theorem implies
$ \v{\zeta - \eta} = \sqrt{\rho_0^2-2\rho_0\rho\cos\theta+\rho^2}
\geq \rho_0-{\rho}\cos\theta$, 
which yields
\begin{equation}
\v{\xi - \zeta} + \v{\zeta - \eta} \geq \rho_0+\rho(1-\cos\theta)
= \rho_0+\frac{\rho\sin^2\theta}{1+\cos\theta}
\geq \rho_0+\frac{\rho}{2}\sin^2\theta.
\label{yogennteiriniyoru}
\end{equation}
Since $\rho_0\rho\cos\theta = \sigma_1^0\sigma_1+g_\xi(\sigma_0, 0)g_\xi(\sigma)$, it follows that
\begin{equation}
\v{\sin\theta}^2 = \frac{\rho_0^2\sigma_2^2  
+ (g_\xi(\sigma)\sigma^0_1-\sigma_1g_\xi(\sigma^0_1, 0))^2}{\rho_0^2\rho^{2}}
\geq \frac{\sigma_2^2}{\rho^{2}},
\label{kakudonosinnohyouka}
\end{equation}
which implies (i) of Proposition \ref{estimates of lengths with a relaxed assumption than IK}. 
\par
We put $r = \v{\sigma}$ and $\omega_j = \sigma_j/r$ $(j = 1, 2)$. Take any $0 < \epsilon < 1/2$
fixed later. For $\omega_1 \leq 1-\epsilon$, $\omega_2^2 \geq 1 - (1-\epsilon)^2 > \epsilon$ holds,
which yields $\sigma_2^2 = r^2\omega_2^2 \geq \epsilon\v{\sigma}^2$. Thus, we get
\begin{align*}
\v{\xi - \zeta} + \v{\zeta - \eta} \geq \rho_0+\frac{\sigma_2^2}{2\rho}
&\geq 
\rho_0+\frac{\epsilon}{2\rho}\v{\sigma^2}
\geq
\rho_0+\frac{\epsilon}{2\rho}(\frac{(\sigma^0_1)^2}{(2r_0)^2}\sigma_1^2+\sigma_2^2)
\quad
\\
&\hskip7mm
(\zeta \in \partial{D}\cap{B(\xi, 2r_0)}, \omega_1 \leq 1-\epsilon).
\nonumber
\end{align*}
Hence, to obtain (ii) of Proposition \ref{estimates of lengths with a relaxed assumption than IK}, 
from (\ref{yogennteiriniyoru}) and (\ref{kakudonosinnohyouka}), and 
$$
\frac{(\sigma_0^1)^2}{\rho_0^2} \geq \frac{1}{1+16M_1^2r_0^2}
$$ 
given by 
(\ref{equivalence of length and local coordinate expression}),
it suffices to show that 
\begin{align}
\frac{g_\xi(\sigma^0_1, 0)}{\sigma^0_1} - \frac{g_\xi(\sigma)}{\sigma_1}
\geq \frac{M_0}{12}\sigma^0_1
\quad(\v{\sigma^0_1} < r_1, \v{\sigma} < r_1, \sigma_1 < \frac{2}{3}\sigma^0_1, \omega_1 \geq 1 - \epsilon)
\label{key estimate for (ii)}
\end{align}
if we choose $0 < \epsilon < 1$ sufficiently small. 

\par
Since $\partial D$ is $C^{2}$ and $\partial D$ is strictly convex, 
$g_\xi$ is expressed by 
\begin{align*}
g_\xi(\sigma) = \sum_{ij = 1}^2a^{ij}_\xi(\sigma)\sigma_i\sigma_j 
\quad\text{and}\quad
a^{ij}_\xi(\sigma) = \int_0^1(1 - \theta)\partial_{\sigma_i}\partial_{\sigma_j}g_\xi(\theta\sigma)d\theta
\quad(i, j = 1, 2). 
\end{align*}
Note that each $a^{ij}_\xi \in C(U_\xi)$ is uniformly bounded for $\v{\sigma} \leq r_1$. Hence, 
there exists a constant $M_2 > 0$ such that $\v{a_\xi^{ij}(\sigma)} \leq M_2$ for $\v{\sigma} \leq r_1$. 
Note that this constant $M_2 > 0$ does not depend on $\xi \in \partial{D}$ and $r_1 > 0$.
\par
From (\ref{estimate of g_xi}), $a^{11}_\xi(\sigma^0_1, 0) \geq M_0$ ($(\sigma^0_1, 0) \in U_\xi$).
For this $M_0 > 0$, 
$\v{a_\xi^{11}(\sigma) - a_\xi^{11}(0)} < M_0/8$ ($\v{\sigma} \leq r_1$) if we take
$r_1 > 0$ sufficiently small. 
Note that this $r_1 > 0$ (and $r_0 > 0 $ also) can be chosen as 
a constant independent of $\xi \in \partial{D}$ 
since Lemma \ref{equi-continuous of g_xi} implies that 
$g_\xi$ is equi-continuous with respect to $\xi \in \partial{D}$. 
Hence, it follows that 
\begin{align*}
\v{a_\xi^{11}(\sigma^0_1, 0) - a_\xi^{11}(\sigma)} \leq \v{a_\xi^{11}(\sigma^0_1, 0) - a_\xi^{11}(0, 0)}
+ \v{a_\xi^{11}(0, 0) - a_\xi^{11}(\sigma)} \leq M_0/4, 
\end{align*}
which yields 
\begin{align*}
a^{11}_\xi(\sigma^0_1, 0) - \frac{2}{3}a^{11}_\xi(\sigma) 
\geq \frac{1}{3}a_\xi^{11}(\sigma^0_1, 0) - \frac{2}{3}\v{a_\xi^{11}(\sigma^0_1, 0) - a_\xi^{11}(\sigma)}
\geq \frac{1}{6}M_0,
\end{align*}
for $\v{\sigma^0_1} < r_1$ and $\v{\sigma} < r_1$, and 
$$
a^{11}_\xi(\sigma) \geq a^{11}_\xi(\sigma^0_1, 0) - \v{a^{11}_\xi(\sigma^0_1, 0) - a^{11}_\xi(\sigma)}
\geq M_0 - M_0/4 > 0
$$
for $\v{\sigma} < r_1$. 
\par
When $\omega_1 \geq 1 - \epsilon$ and $0 < \epsilon \leq 1/2$, 
$\v{\omega_2} \leq \sqrt{2\epsilon}$ holds, which yields 
$\v{\omega_2}/\omega_1 \leq \sqrt{2\epsilon}/(1-\epsilon) \leq 2\sqrt{2\epsilon} \leq 2$.
Hence, for any $\v{\sigma^0_1} < r_1$ and $\v{\sigma} < r_1$
with $0 < r\omega_1 = \sigma_1 < 2\sigma^0_1/3$, it follows that
\begin{align*}
\frac{g_\xi(\sigma^0_1, 0)}{\sigma^0_1} - \frac{g_\xi(\sigma)}{\sigma_1}
&
\geq a^{11}_\xi(\sigma^0_1, 0)\sigma^0_1 - a^{11}_\xi(\sigma)\sigma_1
- 2M_2\sigma_1\frac{\v{\omega_2}}{\omega_1} - M_2\frac{\omega_2^2}{\omega_1^2}\sigma_1
\\&
\geq a^{11}_\xi(\sigma^0_1, 0)\sigma^0_1 - \frac{2}{3}a^{11}_\xi(\sigma)\sigma^0_1
- 2M_2\frac{2}{3}\sigma^0_12\sqrt{2\epsilon}
- M_2\frac{2}{3}\sigma^0_14\sqrt{2\epsilon}
\\&
\geq \Big(\frac{1}{6}M_0 - \frac{16}{3}M_2\sqrt{2\epsilon}\Big)\sigma_1^0
\geq \frac{M_0}{6}\Big(1 - \frac{32M_2\sqrt{2\epsilon}}{M_0}\Big)\sigma^0_1,
\end{align*}
which implies (\ref{key estimate for (ii)}) if we choose 
$\epsilon = \min\{1/2, M_0^2/(2(64M_2)^2)\}$.
This completes the proof of Proposition \ref{estimates of lengths with a relaxed assumption than IK}.
\hfill$\blacksquare$
\vskip1pc\noindent
\par

\par

Last, we show Lemma \ref{equi-continuous of g_xi} used to show 
Proposition \ref{estimates of lengths with a relaxed assumption than IK}.
\par\noindent
Proof of Lemma \ref{equi-continuous of g_xi}: 
Since $\partial{D}$ is $C^2$ class, for any $\xi \in \partial{D}$, there exist a constant $r_{\xi} > 0 $,  
an open neighborhood $U_\xi$ of the origin $0$ in $\R^2$ 
and a function $g_\xi \in {\mathcal B}^{2}(\R^2)$ with $g_{\xi}(0) = 0$ and 
$\nabla g_{\xi}(0) = 0$ such that 
$$
U_\xi \ni \,\sigma = (\sigma_1, \sigma_2)\mapsto 
s_\xi(\sigma) = \xi + \sigma_1 e_1(\xi)+\sigma_2e_2(\xi) 
- g_\xi(\sigma)\nu_\xi \in \partial D\cap B(\xi, r_{\xi}),
$$
where $\{e_1(\xi), e_2(\xi)\}$ is an orthogonal basis for $T_{\xi}(\partial D)$. 
Take any $\varepsilon_1$ with $0 < \varepsilon_1 \leq 1/4$ fixed later.
We can also assume that 
$\nu_\zeta\cdot\nu_\eta \geq 1-\varepsilon_1$ holds for any $\xi \in \partial{D}$ and 
$\zeta, \eta \in \partial{D}{\cap}B(\xi, r_{\xi})$ since for $C^2$
class surfaces, it is well known that 
there exists a constant $C > 0$
such that $\v{\nu_\xi - \nu_\zeta} \leq C\v{\xi - \zeta}$
($\xi, \zeta \in \partial{D}$). In what follows, we write $e_3(\xi) = -\nu_\xi$.
\par
Since $\nu_{s_\xi(\sigma)}\cdot\nu_{\xi} = 1/\sqrt{1+\v{\nabla_{\sigma}g_\xi(\sigma)}^2}$, 
$\nu_{s_\xi(\sigma)}\cdot\nu_{\xi} \geq 1-\varepsilon_1$ implies $\v{\nabla_{\sigma}g_\xi(\sigma)}^2
\leq 1/(1 - \varepsilon_1)^2-1 \leq 2\epsilon_1/(1 - \varepsilon_1)^2$, which yields
$ \v{\partial_{\sigma_k}g_\xi(\sigma)} \leq 2\sqrt{\varepsilon_1}$ 
($\sigma \in U_{\xi}$, $k = 1, 2$) since $0 < \varepsilon_1 \leq 1/4$.
\par
From compactness of $\partial{D}$, we can choose finitely many points $\xi^{(j)}$ 
$(j = 1, 2, \ldots, N)$ satisfying 
$\partial{D} \subset \cup_{j = 1}^{N}B(\xi^{(j)}, r_{\xi^{(j)}}/4)$.
Put $r_0 = \min_{j = 1, 2, \ldots, N}r_{\xi^{(j)}}/8 > 0$. Note that 
\begin{align}
\partial{D} = \cup_{j = 1}^N\{\zeta \in \partial{D} \vert 
B(\zeta, 2r_0) \subset B(\xi^{(j)}, r^{(j)}/2) \}.
\label{coverling of partial{D}}
\end{align} 
Indeed, for any $\zeta \in \partial{D}$, there exists some $\xi^{(j)} \in \partial{D}$ 
satisfying $\zeta \in B(\xi^{(j)}, r_{\xi^{(j)}}/4)$.
For this $\xi^{(j)}$ and $ z \in B(\zeta, 2r_0)$, 
$\v{z - \xi^{(j)}} \leq \v{z - \zeta}+ \v{\zeta - \xi^{(j)}} < 2r_0 + r_{\xi^{(j)}}/4 
\leq r_{\xi^{(j)}}/2 $, which yields $ B(\zeta, 2r_0) \subset B(\xi^{(j)}, r_{\xi^{(j)}}/2)$.

\par
We take any $j \in \{1, 2, \ldots, N \}$ and 
$\zeta \in \partial{D}$ satisfying $B(\zeta, 2r_0) \subset B(\xi^{(j)}, r^{(j)}/2)$.
We define $V_\zeta \subset \R^2$ by 
$
V_\zeta = \{ \tau = (\tau_1, \tau_2) \in \R^2 \vert 
\zeta + \tau_1e_1(\zeta)+\tau_2e_2(\zeta)+\tau_3e_3(\zeta) \in \partial{D}{\cap}B(\zeta, 2r_0) 
\text{ for some }
\tau_3 \in \R \}.
$ 
Note that for any $\tau \in V_\zeta$, there exists a unique $\tau_3 \in \R$ satisfying 
$\eta = \zeta + \tau_1e_1(\zeta)+\tau_2e_2(\zeta)+\tau_3e_3(\zeta) \in \partial{D}{\cap}B(\zeta, 2r_0)$.
Hence $\tau_3$ is a function in $\tau$, which is written by $\tau_3 = h_\zeta(\tau)$.
This fact is shown as follows:
Assume that there exists different $\tilde\tau_3 \in \R$ from $\tau_3$ satisfying 
$\tilde\eta = \zeta + \tau_1e_1(\zeta)+\tau_2e_2(\zeta)+\tilde\tau_3e_3(\zeta) 
\in \partial{D}{\cap}B(\zeta, 2r_0)$. From 
$\partial{D}{\cap}B(\zeta, 2r_0) \subset \partial{D}{\cap}B(\xi^{(j)}, r_{\xi^{(j)}}/2)$, 
$\eta$ and $\tilde\eta$ are written as 
$\eta = \xi^{(j)}+\sigma_1e_1(\xi^{(j)})+\sigma_2e_2(\xi^{(j)})+ g_{\xi^{(j)}}(\sigma)e_3(\xi^{(j)}) $ and 
$\tilde\eta = \xi^{(j)}+\tilde\sigma_1e_1(\xi^{(j)})+\tilde\sigma_2e_2(\xi^{(j)}) 
+ g_{\xi^{(j)}}(\tilde\sigma)e_3(\xi^{(j)}) $ 
by taking some $\sigma$ and $\tilde\sigma \in U_{\xi^{(j)}}$, respectively. 
Put $\eta_t = ((1-t)\sigma_1+t\tilde\sigma_1)e_1(\xi^{(j)})+((1-t)\sigma_2+t\tilde\sigma_2)e_2(\xi^{(j)})
+ g_{\xi^{(j)}}((1-t)\sigma+t\tilde\sigma)e_3(\xi^{(j)}) \in 
\partial{D}{\cap}B(\xi^{(j)}, r_{\xi^{(j)}})$ ($0 \leq t \leq 1$). 
From mean value theorem, it follows that 
$g_{\xi^{(j)}}(\tilde\sigma) - g_{\xi^{(j)}}(\sigma) 
= (\tilde\sigma - \sigma)\cdot\partial_{\sigma}g_{\xi^{(j)}}(\sigma^{(0)})$
where $\sigma^{(0)} = (1-t_0)\sigma+t_0\tilde\sigma$ for some $0 < t_0 < 1$.
Hence, we obtain
\begin{align*}
\tilde\eta - \eta &= (\tilde\sigma_1 - \sigma_1)e_1(\xi^{(j)})
+ (\tilde\sigma_2 - \sigma_2)e_2(\xi^{(j)})
+ (g_{\xi^{(j)}}(\tilde\sigma) - g_{\xi^{(j)}}(\sigma))e_3(\xi^{(j)})
\\
&= \sum_{k = 1}^2(\tilde\sigma_k - \sigma_k)(e_k(\xi^{(j)})
+\partial_{\sigma_k}g_{\xi^{(j)}}(\sigma^{(0)})e_3(\xi^{(j)}))
\in T_{\eta_{t_0}}(\partial{D}).
\end{align*}
Thus $(\tilde\eta - \eta)\cdot{\nu_{\eta_{t_0}}} = 0$, which yields 
$(\tilde\tau_3-\tau_3)\nu_{\zeta}\cdot{\nu_{\eta_{t_0}}} = 0$. 
This gives a contradiction since 
$\nu_{\zeta}\cdot{\nu_{\eta_{t_0}}} \geq 1 - \varepsilon_1 > 0$ holds.
Hence, $\tau_3$ is uniquely determined.

\par
From the above argument, the map
$V_\zeta \in \tau \mapsto \zeta+\tau_1e_1(\zeta)+\tau_2e_2(\zeta)+h_\zeta(\tau)e_3(\zeta)
\in \partial{D}{\cap}B(\zeta, 2r_0)$ is bijective, and 
the function $h_\zeta$ is related to $g_{\xi^{(j)}}$ by the equality
$ \xi^{(j)} + \sigma_1e_1(\xi^{(j)})+\sigma_2e_2(\xi^{(j)})+g_{\xi^{(j)}}(\sigma)e_3(\xi^{(j)})
= \zeta + \tau_1e_1(\zeta)+\tau_2e_2(\zeta)+h_\zeta(\tau)e_3(\zeta)$,
which is equivalent to the following equalities:
\begin{align*}
h_\zeta(\tau) &= e_3(\zeta)\cdot(\xi^{(j)} - \zeta + \sigma_1e_1(\xi^{(j)})+\sigma_2e_2(\xi^{(j)})
+g_{\xi^{(j)}}(\sigma)e_3(\xi^{(j)}))
\\
\tau_k &= e_k(\zeta)\cdot(\xi^{(j)} - \zeta + \sigma_1e_1(\xi^{(j)})+\sigma_2e_2(\xi^{(j)})
+g_{\xi^{(j)}}(\sigma)e_3(\xi^{(j)}))
\quad(k = 1, 2).
\end{align*}
We put $\tau = \Phi_{\zeta, \xi^{(j)}}(\sigma)$, which has the inverse
$\sigma = \Psi_{\xi^{(j)}, \zeta}(\tau)$ for $\tau \in V_\zeta $. 
Since 
$\{e_1(\zeta), e_2(\zeta), e_3(\zeta)\}$ and $\{e_1(\xi^{(j)}), e_2(\xi^{(j)}), e_3(\xi^{(j)})\}$
are orthogonal basis, 
it follows that
\begin{align*}
\Big\vert{\rm det}\big(\frac{\partial\tau}{\partial\sigma}\big) \Big\vert
&\geq 
\Big\vert{\rm det}\begin{pmatrix}
e_1(\zeta){\cdot}e_1(\xi^{(j)}) & e_1(\zeta){\cdot}e_2(\xi^{(j)}) 
\\
e_2(\zeta){\cdot}e_1(\xi^{(j)}) & e_2(\zeta){\cdot}e_2(\xi^{(j)}) 
\end{pmatrix}
\Big\vert
- 2(\v{\partial_{\sigma_1}g_{\xi^{(j)}}(\sigma)}+\v{\partial_{\sigma_2}g_{\xi^{(j)}}(\sigma)})
\\&
\qquad
(\sigma \in U_{\xi^{(j)}}, s_{\xi^{(j)}}(\sigma) \in \partial{D}{\cap}B(\zeta, 2r_0))
\text{ and } \zeta \in \partial{D}{\cap}B(\zeta, 2r_0))
\end{align*}
and 
\begin{align*}
1 &\leq \Big\vert{\rm det}\begin{pmatrix}
e_1(\zeta){\cdot}e_1(\xi^{(j)}) & e_1(\zeta){\cdot}e_2(\xi^{(j)}) 
\\
e_2(\zeta){\cdot}e_1(\xi^{(j)}) & e_2(\zeta){\cdot}e_2(\xi^{(j)}) 
\end{pmatrix}
\Big\vert
\v{e_3(\zeta){\cdot}e_3(\xi^{(j)})}
+ 2\sqrt{2}\sum_{k = 1}^2\v{e_k(\zeta){\cdot}e_3(\xi^{(j)})}^2
\\
&\leq
\Big\vert{\rm det}\begin{pmatrix}
e_1(\zeta){\cdot}e_1(\xi^{(j)}) & e_1(\zeta){\cdot}e_2(\xi^{(j)}) 
\\
e_2(\zeta){\cdot}e_1(\xi^{(j)}) & e_2(\zeta){\cdot}e_2(\xi^{(j)}) 
\end{pmatrix}
\Big\vert
+ 2\sqrt{2}(1-\v{e_3(\zeta){\cdot}e_3(\xi^{(j)})}^2).
\end{align*}
From these estimates and 
$e_3(\zeta){\cdot}e_3(\xi) = \nu_{\zeta}{\cdot}\nu_{\xi} \geq 1-\varepsilon_1$, we obtain
\begin{align*}
\Big\vert{\rm det}\big(\frac{\partial\tau}{\partial\sigma}\big) \Big\vert
&\geq 1-4\sqrt{2}\varepsilon_1
- 2(\v{\partial_{\sigma_1}g_{\xi^{(j)}}(\sigma)}+\v{\partial_{\sigma_2}g_{\xi^{(j)}}(\sigma)})
\geq 1-4\sqrt{2}\varepsilon_1 -8\sqrt{\varepsilon_1}
\\&\quad
(\sigma \in U_{\xi^{(j)}}, s_{\xi^{(j)}}(\sigma) \in \partial{D}{\cap}B(\zeta, 2r_0))
\text{ and } \zeta \in \partial{D}{\cap}B(\zeta, 2r_0)). 
\end{align*}
\par
From now on, take $\varepsilon_1 = 1/1024$ to be 
$\big\vert{\rm det}\big(\frac{\partial\tau}{\partial\sigma}\big) \big\vert \geq 1/2$
for $\sigma \in U_{\xi^{(j)}}$, $s_{\xi^{(j)}}(\sigma) \in \partial{D}{\cap}B(\zeta, 2r_0)$
and $\zeta \in \partial{D}{\cap}B(\zeta, 2r_0)$). 
Thus, the implicit function theorem implies that $\Psi_{\xi^{(j)}, \zeta} \in C^2(V_\zeta)$ and
$\frac{\partial{\Psi_{\xi^{(j)}, \zeta}}}{\partial\tau}(\tau) = 
\big(\frac{\partial{\Phi_{\zeta, \xi^{(j)}}}}{\partial\sigma}(\sigma)\big)^{-1}$
($\tau \in V_\zeta$).  From these facts and 
$h_\zeta(\tau) = g_{\xi^{(j)}}(\Psi_{\xi^{(j)}, \zeta}(\tau))$, we can see that for any $\alpha$ with
$\v{\alpha} \leq 2$, the function $\partial_{\tau}^{\alpha}h_\zeta(\tau)$ 
is equi-continuous with respect to $\zeta$ and $j = 1, 2, \ldots, N$. Thus, we obtain
Lemma \ref{equi-continuous of g_xi} if we note (\ref{coverling of partial{D}}).
\hfill$\blacksquare$
\vskip1pc
%

%
%
%

%
%
%
%
\par
%
%
%
%
%
%
\par
\footnotesize
{\sc Kawashita, Mishio}
\par
{\sc Department of Mathematics, Graduate School of Science,
 Hiroshima University, }
\par {\sc Higashi-Hiroshima 739-8526 Japan}

\end{document}